\documentclass[envcountsame,envcountsect,leqno]{svmult}

\usepackage{a4,amsfonts,exscale}

\usepackage{amsmath}
\usepackage{amssymb}
\usepackage{amscd}
\usepackage{textcomp}

\usepackage{xy}
\xyoption{all}

\newdir^{ (}{{}*!/-3pt/\dir^{(}}    
\newdir^{  }{{}*!/-3pt/\dir^{}}    
\newdir_{ (}{{}*!/-3pt/\dir_{(}}    
\newdir_{  }{{}*!/-3pt/\dir_{}}

\newcounter{zahl}

\newcommand{\DS}{\displaystyle}
\newcommand{\TS}{\textstyle}

\newcommand{\SSC}{\scriptscriptstyle}

\DeclareMathOperator{\Aut}{Aut}

\DeclareMathOperator{\End}{End}

\DeclareMathOperator{\Frob}{Frob}
\DeclareMathOperator{\Gal}{Gal}
\DeclareMathOperator{\GL}{GL}
\DeclareMathOperator{\Koh}{H}

\DeclareMathOperator{\Hom}{Hom}
\newcommand{\CHom}{{\cal H}om}
\DeclareMathOperator{\Id}{Id}
\DeclareMathOperator{\Isom}{Isom}
\DeclareMathOperator{\Lie}{Lie}

\DeclareMathOperator{\PGL}{PGL}

\newcommand{\Rkoh}{{\rm R}}

\DeclareMathOperator{\Spec}{Spec}
\DeclareMathOperator{\Spf}{Spf}
\DeclareMathOperator{\Stab}{Stab}

\DeclareMathOperator{\Tor}{Tor}

\DeclareMathOperator{\Var}{V}

\newcommand{\alg}{{\rm alg}}

\DeclareMathOperator{\coker}{coker}
\DeclareMathOperator{\di}{div}

\DeclareMathOperator{\id}{\,id}
\DeclareMathOperator{\im}{im}

\newcommand{\red}{{\rm red}}

\newcommand{\rig}{{\rm rig}}
\DeclareMathOperator{\rk}{rk}

\DeclareMathOperator{\Sch}{\CS \!{\it ch}}
\DeclareMathOperator{\Nilp}{\CN \!{\it ilp}}
\DeclareMathOperator{\AbSh}{\CA {\it b}-\CS {\it h}}

\DeclareMathOperator{\DrMod}{\CalD {\it r}-\CM {\it od}}
\DeclareMathOperator{\Pol}{\CP {\it ol}}
\DeclareMathOperator{\Vect}{\CV {\it ec}}
\DeclareMathOperator{\Seq}{\CS {\it eq}}
\DeclareMathOperator{\Hecke}{\CH {\it ecke}}
\DeclareMathOperator{\Quot}{\CQ {\it uot}}

\renewcommand{\phi}{\varphi}
\renewcommand{\epsilon}{\varepsilon}

\usepackage{amsfonts}
\newcommand{\BOne}{\hbox{\rm1\kern-2.7pt l\kern.9pt}}

\newcommand{\BA}{{\mathbb{A}}}

\newcommand{\BE}{{\mathbb{E}}}
\newcommand{\BF}{{\mathbb{F}}}
\newcommand{\BG}{{\mathbb{G}}}

\newcommand{\BM}{{\mathbb{M}}}

\newcommand{\BP}{{\mathbb{P}}}

\newcommand{\BZ}{{\mathbb{Z}}}

\newcommand{\CA}{{\cal{A}}}

\newcommand{\CalD}{{\cal{D}}}

\newcommand{\CF}{{\cal{F}}}
\newcommand{\CG}{{\cal{G}}}
\newcommand{\CH}{{\cal{H}}}
\newcommand{\CI}{{\cal{I}}}
\newcommand{\CJ}{{\cal{J}}}
\newcommand{\CK}{{\cal{K}}}

\newcommand{\CM}{{\cal{M}}}
\newcommand{\CN}{{\cal{N}}}
\newcommand{\CO}{{\cal{O}}}
\newcommand{\CP}{{\cal{P}}}
\newcommand{\CQ}{{\cal{Q}}}
\newcommand{\CR}{{\cal{R}}}
\newcommand{\CS}{{\cal{S}}}
\newcommand{\CT}{{\cal{T}}}
\newcommand{\CU}{{\cal{U}}}
\newcommand{\CV}{{\cal{V}}}

\newcommand{\CX}{{\cal{X}}}
\newcommand{\CY}{{\cal{Y}}}
\newcommand{\CZ}{{\cal{Z}}}

\renewcommand{\Gamma}{\varGamma}

\newcommand{\rsurj}{\mbox{\mathsurround=0pt \;$\longrightarrow \hspace{-0.7em} \to$\;}}
\newcommand{\lsurj}{\mbox{\mathsurround=0pt \;$\gets \hspace{-0.7em} \longleftarrow$\;}}

\newcommand{\rbij}{\mbox{\mathsurround=0pt \;$-\hspace{-0.75em}\stackrel{\sim\quad}{\longrightarrow}$\;}}

\let\setminus\smallsetminus

\newcommand{\es}{\enspace}

\newcommand{\dual}{^{\SSC\lor}}

\newcommand{\ul}[1]{{\underline{#1}}}
\newcommand{\ol}[1]{{\overline{#1}}}
\newcommand{\wh}[1]{{\widehat{#1}}}
\newcommand{\wt}[1]{{\widetilde{#1}}}

\usepackage{ifthen}

\newcommand{\invlim}[1][]{\ifthenelse{\equal{#1}{}}
{\DS \lim_{\longleftarrow}}
{\DS \lim_{\underset{#1}{\longleftarrow}}}
}

\newcommand{\dirlim}[1][]{\ifthenelse{\equal{#1}{}}
{\DS \lim_{\longrightarrow}}
{\DS \lim_{\underset{#1}{\longrightarrow}}}
}

\newcommand{\dbl}{{\mathchoice{\mbox{\rm [\hspace{-0.15em}[}}
                              {\mbox{\rm [\hspace{-0.15em}[}}
                              {\mbox{\scriptsize\rm [\hspace{-0.15em}[}}
                              {\mbox{\tiny\rm [\hspace{-0.15em}[}}}}
\newcommand{\dbr}{{\mathchoice{\mbox{\rm ]\hspace{-0.15em}]}}
                              {\mbox{\rm ]\hspace{-0.15em}]}}
                              {\mbox{\scriptsize\rm ]\hspace{-0.15em}]}}
                              {\mbox{\tiny\rm ]\hspace{-0.15em}]}}}}
\newcommand{\dpl}{{\mathchoice{\mbox{\rm (\hspace{-0.15em}(}}
                              {\mbox{\rm (\hspace{-0.15em}(}}
                              {\mbox{\scriptsize\rm (\hspace{-0.15em}(}}
                              {\mbox{\tiny\rm (\hspace{-0.15em}(}}}}
\newcommand{\dpr}{{\mathchoice{\mbox{\rm )\hspace{-0.15em})}}
                              {\mbox{\rm )\hspace{-0.15em})}}
                              {\mbox{\scriptsize\rm )\hspace{-0.15em})}}
                              {\mbox{\tiny\rm )\hspace{-0.15em})}}}}

\def\?{\ ???\ \immediate\write16{}%
\immediate\write16{Warning: There was still a question mark . . . }%
\immediate\write16{}}

\long\def\forget#1{}

\spnewtheorem{bigexample}[theorem]{Example}{\rm\bf}{\rm}
\spnewtheorem{construction}[theorem]{Construction}{\rm\bf}{\rm}
\spnewtheorem{notation}[theorem]{Notation}{\rm\bf}{\rm}
\spnewtheorem{point}[theorem]{}{\rm\bf}{\rm}
\spnewtheorem{uniftheorem}[theorem]{Uniformization Theorem}{\rm\bf}{\it}

\smartqed

\def\longto{\longrightarrow}
\def\into{\hookrightarrow}

\def\isoto{\arrover{\sim}}

\newbox\mybox
\def\arrover#1{\mathrel{
       \setbox\mybox=\hbox spread 1.4em{\hfil$\scriptstyle#1$\hfil}
       \vbox{\offinterlineskip\copy\mybox
             \hbox to\wd\mybox{\rightarrowfill}}}}

\begin{document}


\author{Urs Hartl}

\institute{University of Freiburg,
Institute of Mathematics,
Eckerstr.~1,
D-79104 Freiburg,
Germany,
E-mail: 
urs.hartl@math.uni-freiburg.de}

\date{June 9th, 2004}

\title*{Uniformizing the Stacks of Abelian Sheaves}

\maketitle
\begin{abstract}
Elliptic sheaves (which are related to Drinfeld modules) were introduced by Drinfeld and further studied by Laumon--Rapoport--Stuhler and others. They can be viewed as function field analogues of elliptic curves and hence are objects ``of dimension~$1$''. Their higher dimensional generalizations are called abelian sheaves.
In the analogy between function fields and number fields, abelian sheaves are counterparts of abelian varieties. In this article we study the moduli spaces of abelian sheaves and prove that they are algebraic stacks. We further transfer results of \v{C}erednik--Drinfeld and Rapoport--Zink on the uniformization of Shimura varieties to the setting of abelian sheaves. Actually the analogy of the \v{C}erednik--Drinfeld uniformization is nothing but the uniformization of the moduli schemes of Drinfeld modules by the Drinfeld upper half space. Our results generalize this uniformization. The proof closely follows the ideas of Rapoport--Zink. In particular, analogies of $p$-divisible groups play an important role. As a crucial intermediate step we prove that in a family of abelian sheaves with good reduction at infinity, the set of points where the abelian sheaf is uniformizable in the sense of Anderson, is formally closed.

\end{abstract}


\renewcommand{\theequation}{\thesection.\arabic{equation}}

%
%

\section*{Introduction}  \label{Intro}

In arithmetic algebraic geometry the moduli spaces of abelian varieties are of great importance. For instance they have played a major role in Faltings' proof of the Mordell conjecture \cite{Faltings}, the proof of Fermat's Last Theorem \cite{FLT}, and the proof of Langlands reciprocity for $\GL_n$ over non-archimedean local fields of characteristic zero by Harris--Taylor~\cite{Harris}. Therefore their structure and especially their reduction at bad primes is intensively studied. One way to investigate their reduction is through $p$-adic uniformization. This was begun by \v{C}erednik~\cite{Cerednik} and Drinfeld~\cite{Drinfeld2} and continued by Rapoport--Zink~\cite{RZ}. \v{C}erednik--Drinfeld obtained the uniformization of certain Shimura curves of EL-type by a formal scheme whose associated rigid-analytic space is Drinfeld's $p$-adic upper half plane. This formal scheme can be viewed as a moduli space for $p$-divisible groups which are isogenous to a fixed supersingular $p$-divisible group. See Boutot--Carayol~\cite{BC} for a detailed account. Rapoport--Zink generalized these results to the (partial) uniformization of higher dimensional Shimura varieties by more general moduli spaces for $p$-divisible groups.

In this article we study \emph{abelian sheaves} as positive characteristic analogues of abelian varieties. We investigate their moduli spaces and prove that these are algebraic stacks. Then our aim is to transfer the above uniformization results to the case of positive characteristic. For the case considered by \v{C}erednik--Drinfeld this was accomplished already by Drinfeld~\cite{Drinfeld}; see below. We use this as a guide line to transfer the results of Rapoport--Zink. We also obtain a partial uniformization of the moduli stacks of abelian sheaves. In the Hilbert-Blumenthal situation a similar uniformization result was obtained by Stuhler~\cite{Stuhler} using different methods.

\smallskip

Let us explain what abelian sheaves are by first going back 30 years to Drinfeld's elliptic sheaves.
Exploiting the analogy between number fields and function fields, Drinfeld~\cite{Drinfeld, Drinfeld3} invented the notions of \emph{elliptic modules} (today called \emph{Drinfeld modules}) and the dual notion of \emph{elliptic sheaves}. These structures are analogues of elliptic curves for characteristic $p$ in the following sense. Their endomorphism rings are rings of integers in global function fields of positive characteristic or orders in central division algebras over the later. 
On the other hand, the moduli spaces are varieties over smooth curves over a finite field. Through these two aspects in which global function fields of positive characteristic come into play, elliptic sheaves and variants of them proved to be fruitful for establishing the Langlands correspondence for $\GL_n$ over local and global function fields of positive characteristic. See the work of Drinfeld~\cite{Drinfeld, Drinfeld5, Drinfeld4}, Laumon--Rapoport--Stuhler~\cite{LRS}, and Lafforgue~\cite{Lafforgue}. Beyond this the analogy between elliptic modules and elliptic curves is abundant.

In this spirit, Anderson~\cite{Anderson} introduced higher dimensional generalizations of Drinfeld's elliptic modules and called them \emph{abelian $t$-modules}. The concept of \emph{abelian sheaves} is a higher dimensional generalization of elliptic sheaves. Both serve as characteristic $p$ analogues of abelian varieties. Abelian sheaves were studied in various special instances in the past. In this article we intend to give a systematic treatment. The definition of abelian sheaves is as follows. Let $C$ be a smooth projective curve over $\BF_q$ and let $\infty\in C(\BF_q)$ be a fixed point. For every $\BF_q$-scheme $S$ we denote by $\sigma$ the endomorphism of $C\times_{\BF_q} S$ that acts as the identity on the coordinates of $C$ and as $b\mapsto b^q$ on the sections $b\in \CO_S$. Now an \emph{abelian sheaf of rank $r$ and dimension $d$ over $S$} consists of the following data: a collection of locally free sheaves $\CF_i$ of rank $r$ on $C\times_{\BF_q} S$ satisfying a certain periodicity condition. These sheaves are connected by two commuting sets of morphisms $\Pi_i:\CF_i\to\CF_{i+1}$ and $\tau_i:\sigma^\ast\CF_i\to\CF_{i+1}$ such that $\coker \Pi_i$ and $\coker\tau_i$ are locally free $\CO_S$-modules of rank $d$, supported respectively on $\infty\times S$ and on the graph of a morphism $c:S\to C$ called the \emph{characteristic} of the abelian sheaf.
An abelian sheaf of dimension $1$ is the same as an elliptic sheaf. In this sense abelian sheaves are higher dimensional elliptic sheaves. The notion of abelian sheaf is dual to the notion of abelian $t$-module and related to Anderson's $t$-motives. In fact, if the characteristic is different from $\infty$ an abelian sheaf over a field is nothing but a pure $t$-motive equipped with additional structure at infinity (Section~\ref{SectPolTMotives}). Therefore we like to view abelian sheaves as characteristic $p$ analogues of polarized abelian varieties.

Abelian sheaves with appropriately defined level structure possess moduli spaces which are algebraic stacks locally of finite type over the curve $C$. The morphism to $C$ is given by assigning to an abelian sheaf its characteristic. We denote the algebraic stacks of abelian sheaves of rank $r$ and dimension $d$ with $H$-level structures by $\AbSh^{r,d}_H$. Here $H\subset\GL_r(\BA_f)$ is a compact open subgroup and $\BA_f$ are the finite adeles of $C$.
It should be noted that opposed to the case of elliptic sheaves, the stacks of abelian sheaves will in general not be schemes, not even if we add high level structures. This is due to the fact that for every level there are abelian sheaves having non-trivial automorphisms (see Remark~\ref{RemStacksNotSmooth}). The non-representability also reflects in the uniformization; see below.

\smallskip

Then our aim is to study the uniformization of these moduli stacks at $\infty$. Let $z$ be a uniformizing parameter of $C$ at $\infty$. In the case of elliptic sheaves, Drinfeld~\cite{Drinfeld, Drinfeld2} showed that the moduli stacks are in fact smooth affine schemes which can be uniformized by a formal scheme $\Omega^{(r)}$. This formal scheme is the characteristic $p$ version of the one used by \v{C}erednik and Drinfeld to uniformize Shimura curves. Correspondingly it is a moduli space for certain formal groups on which multiplication with $z$ is an isogeny. All this was worked out in detail by Genestier~\cite{Genestier}. We like to call these formal groups ``$z$-divisible groups''. They play an important role also in our uniformization results.

So let us next explain some facts about $z$-divisible groups. Naturally these groups are of most use over schemes on which $z$ is not a unit. Therefore we will from now on work over schemes $S$ in $\Nilp_{\BF_q\dbl z\dbr}$, the category of $\BF_q\dbl z\dbr$-schemes on which $z$ is locally nilpotent. However, since it is important to separate the two roles played by $z$ as a uniformizing parameter at $\infty$ and as an element of $\CO_S$, we use the symbol $z$ only for the first and we denote the image of $z$ in $\CO_S$ by $\zeta$. Then $S$ belongs to $\Nilp_{\BF_q\dbl\zeta\dbr}$.

Classically $p$-divisible groups may be studied via their Dieudonn\'e modules. There is a corresponding notion for $z$-divisible groups. 
A \emph{Dieudonn\'e $\BF_q\dbl z\dbr$-module over $S$} is a finite locally free $\CO_S\dbl z\dbr$-module $\wh{\CF}$ with a $\sigma$-linear endomorphism $F$ such that
\begin{enumerate}
\item 
$\coker F$ is locally free as an $\CO_S$-module,
\item
$(z-\zeta)$ is nilpotent on $\coker F$,
\end{enumerate}
The theory of $z$-divisible groups and their Dieudonn\'e $\BF_q\dbl z\dbr$-modules resembles many facets of the theory of $p$-divisible groups; see \cite{Crystals}. 

Also $z$-divisible groups are related to abelian sheaves through their Dieudonn\'e $\BF_q\dbl z\dbr$-modules. Namely, the completion of an abelian sheaf over $S$ at $\infty\times S$ is a Dieudonn\'e $\BF_q\dbl z\dbr$-module. The connection between abelian sheaves and $z$-divisible groups parallels the situation for abelian varieties. In particular, there is an analogue of the Serre-Tate-Theorem relating the deformation theory of abelian sheaves to the deformation theory of their $z$-divisible groups.

\smallskip

Finally we come to the uniformization of $\AbSh^{r,d}_H$ at infinity. We begin by describing the uniformizing spaces.
Let $r,d,k,\ell$ be positive integers with $\frac{d}{r}=\frac{k}{\ell}$ and $k$ and $\ell$ relatively prime. Let $\CO_{\!\Delta}$ be the ring of integers in the central skew field over $\BF_q\dpl z\dpr$ of invariant $k/\ell$. A \emph{special $z$-divisible $\CO_{\!\Delta}$-module over $S\in\Nilp_{\BF_{q^\ell}\dbl\zeta\dbr}$} is a $z$-divisible group $E$ of height $r\ell$ and dimension $d\ell$ with an action of $\CO_{\!\Delta}$ prolonging the action of $\BF_q\dbl z\dbr$, such that the inclusion $\BF_{q^\ell}\subset\CO_{\!\Delta}$ makes $\Lie E$ into a locally free $\BF_{q^\ell}\otimes_{\BF_q}\CO_S$-module of rank $d$. The $z$-divisible groups associated to abelian sheaves of rank $r$ and dimension $d$ are special $z$-divisible $\CO_{\!\Delta}$-modules.
Let $\BE$ be a special $z$-divisible $\CO_{\!\Delta}$-module over $\Spec\BF_{q^\ell}$. Then the moduli problem of special $z$-divisible $\CO_{\!\Delta}$-modules which are isogenous to $\BE$ is solved by a formal scheme $G$ locally formally of finite type over $\Spf\BF_{q^\ell}\dbl\zeta\dbr$. The later means that the reduced closed subscheme $G_\red$ of $G$ is locally of finite type over $\Spec\BF_{q^\ell}$.

We fix an abelian sheaf $\ol{\BM}$ of rank $r$ and dimension $d$ over $\Spec\BF_{q^\ell}$ whose restriction to $C\setminus\infty$ satisfies $\tau_i=\Id_r\cdot\sigma^\ast$ and we let $\BE$ be its $z$-divisible group. The Newton polygon of $\BE$ is a straight line. Let $Z$ be the set of points $s$ of $\AbSh^{r,d}_H \times_C\,\infty$ such that the universal abelian sheaf $\ul\CF_s$ over $s$ is isogenous to $\ol{\BM}$ over an algebraic closure of $\kappa(s)$. It is an important step to show that $Z$ is the set of points over which the Newton polygons of $\BE$ and of the $z$-divisible group associated to $\ul\CF_s$ coincide. This implies that $Z$ is a closed subset. We consider the formal completion $\AbSh^{r,d}_H{}_{\!\!/Z}$ of $\AbSh^{r,d}_H$ along $Z$. It is no longer an algebraic stack, but it is a \emph{formal algebraic stack over $\Spf\BF_{q^\ell}\dbl\zeta\dbr$}. Formal algebraic stacks are generalizations of algebraic stacks in the same sense as formal schemes generalize usual schemes. (The relevant facts on formal algebraic stacks are collected in an appendix.) Now we can uniformize $\AbSh^{r,d}_H{}_{\!\!/Z}$ as follows. Being isogenous to the $z$-divisible group $\BE$ of $\ol{\BM}$, the universal special $z$-divisible $\CO_{\!\Delta}$-module on $G$ gives rise to an abelian sheaf on $G$ which we call its \emph{algebraization}. Let $J(Q)$ be the group of quasi-isogenies of $\ol{\BM}$. There are natural embeddings of the group $J(Q)$ into $\GL_r(\BA_f)$ and into the group of quasi-isogenies of $\BE$. The later group acts on the formal scheme $G$. We let $J(Q)$ act diagonally on $G\times\GL_r(\BA_f)$. Taking into account level structures we obtain the following 

\medskip

\noindent
{\bf Uniformization Theorem~\ref{ThmUnifOfAbSh}.} {\it There is a canonical 1-iso\-mor\-phism of formal algebraic $\Spf\BF_{q^\ell}\dbl\zeta\dbr$-stacks}
\[
\Theta:\es\raisebox{-0.2mm}{$J(Q)$}\backslash \raisebox{0mm}{$G\times\GL_r(\BA_f)$}/\raisebox{-0.2mm}{$H$} \es\isoto\es\AbSh^{r,d}_H{}_{\!\!/Z}\times_{\Spf\BF_q\dbl\zeta\dbr}\;\Spf\BF_{q^\ell}\dbl\zeta\dbr
\]

\medskip

At this point note that the non-representability of the algebraic stacks $\AbSh^{r,d}_H$ also reflects in the uniformization. Namely, since all unipotent subgroups of $J(Q)$ are torsion, in general a discrete subgroup of $J(Q)$ cannot act fixed point free on $G$. So the quotients $\raisebox{-0.2mm}{$J(Q)$}\backslash \raisebox{0mm}{$G\times\GL_r(\BA_f)$}/\raisebox{-0.2mm}{$H$}$ can only be formal algebraic stacks and not formal algebraic spaces. This phenomenon does not occur in the $p$-adic uniformization of Shimura-varieties.

\smallskip

We like to mention an interesting aspect of the proof that is also related to the uniformizability of $t$-motives. Namely by work of Gardeyn~\cite{Gardeyn} the $t$-motive associated to an abelian sheaf over a complete field extension $K$ of $\BF_q\dpl\zeta\dpr$ is \emph{uniformizable} in the sense of Anderson~\cite{Anderson} if and only if firstly the abelian sheaf extends to an abelian sheaf over the valuation ring $R$ of a finite extension of $K$, and secondly its reduction modulo the maximal ideal of $R$ is isogenous to $\ol\BM$. We show that the second condition is closed. More precisely if $\ul\CF$ is an abelian sheaf or rank $r$ and dimension $d$ over $S\in\Nilp_{\BF_q\dbl\zeta\dbr}$, then the set of points in $S$ over which $\ul\CF$ is isogenous to $\ol\BM$ is closed. This is the key ingredient in the proof of the Uniformization Theorem. It is proved in Section~\ref{SectUnifClosed} where we avoid the language of stacks and proceed in more down-to-earth terms.

\smallskip

Let us end by explaining the relation of our Uniformization Theorem to the results of Rapoport--Zink~\cite{RZ} and Drinfeld~\cite{Drinfeld}.
In the case of elliptic sheaves we have
\[
\AbSh^{r,d}_H{}_{\!\!/Z}\es=\es\AbSh^{r,d}_H\times_C \;\Spf\BF_{q}\dbl\zeta\dbr\,.
\]
This follows from the fact that there is only one polygon between the points $(0,0)$ and $(r,1)$ with non-negative slopes and integral break points, namely the straight line. So all special $z$-divisible $\CO_{\!\Delta}$-modules have the same Newton polygon as $\BE$ and therefore $Z$ is all of $\AbSh^{r,d}_H\times_C\;\infty$.
Moreover in this case $G$ is the formal scheme $\Omega^{(r)}$ introduced above and $J(Q)\cong\GL_r(Q)$. So we recover Drinfeld's uniformization theorem
\[
\raisebox{-0.2mm}{$\GL_r(Q)$}\backslash \raisebox{0mm}{$\Omega^{(r)}\times\GL_r(\BA_f)$}/\raisebox{-0.2mm}{$H$} \es\isoto\es M^r_H\times_C \;\Spf\BF_{q^\ell}\dbl\zeta\dbr
\]
where $M^r_H$ is the moduli scheme of Drinfeld modules of rank $r$ with $H$-level structure.

Compared to \cite{RZ}, our uniformization theorem is analogous to the uniformization of the formal completion of Shimura varieties along the most supersingular isogeny class. In this sense uniformizable abelian sheaves correspond to supersingular abelian varieties. 
There is no doubt that the more general uniformization in \cite{RZ} of an arbitrary isogeny class of abelian varieties also carries over to the setting of abelian sheaves. Furthermore, we have only described the uniformization at $\infty$ in this article. But the analogous uniformization results at other places of $C$ should likewise hold. For example in the case of $\CalD$-elliptic sheaves these were described by Hausberger~\cite{Hausberger}.

%
%

\section*{Table of Contents}
\setcounter{minitocdepth}{1}
\dominitoc

%
%

\section*{Notation} 

Throughout this article we will denote by

\begin{tabbing}
$A=\Gamma(C',\CO_C)$ \es\=\kill
$\BF_q$ \> the finite field having $q$ elements and characteristic $p$, \\[1mm]
$C$ \> a smooth projective geometrically irreducible curve over $\BF_q$, \\[1mm]
$\infty\in C(\BF_q)$ \>a fixed point,\\[1mm]
$C'=C\setminus\infty$ \\[1mm]
$A=\Gamma(C',\CO_{C'})$ \> the ring of regular functions on $C'$, \\[1mm]
$Q =\BF_q(C)$ \> the function field of $C$, viz.\ the field of fractions of $A$,\\[1mm]
$Q_\infty$ \> the completion of $Q$ at $\infty$,\\[1mm]
$A_\infty$ \> the ring of integers in $Q_\infty$,\\[1mm]
$v\in M_f$ \> the finite places of $C$, i.e.\ the points of $C'$, \\[1mm]
$A_v$ \> the completion of $A$ at the finite place $v$ of $C$,\\[1mm]
$\DS\wh{A} = \prod_{v\in M_f} \,A_v$\\[1mm]
$\BA_f = Q\otimes_A \wh{A}$ \> the finite adeles of $C$,\\[1mm]
$d,r,k,\ell$ \> positive integers with $\frac{d}{r}=\frac{k}{\ell}$ and $k$ and $\ell$ relatively prime,\\[1mm]
$\Delta$ \> the central skew field over $Q_\infty$ of invariant $k/\ell$,\\[1mm]
$\CO_{\!\Delta}$ \> its ring of integers.
\end{tabbing}

\noindent
All schemes, as well as their products and morphisms between them, are supposed to be over $\Spec\BF_q$. If $X$ is a scheme we let $\Sch_X$ be the category of $X$-schemes. For two schemes $X$ and $Y$ we write $X\times Y$ for their product over $\Spec\BF_q$. A similar notation will be employed to the tensor product over $\BF_q$.
If $i:Y\into X$ is a closed immersion of schemes and $\CF$ is a quasi-coherent sheaf on $X$ we denote the restriction $i^\ast\CF$ by $\CF|_Y$.
As is customary, we will use the term vector bundle for a locally free coherent sheaf on a scheme. 

Starting with Section~\ref{SectZDivGps} we denote by 
\begin{tabbing}
$A=\Gamma(C',\CO_C)$ \es\=\kill
$\zeta$ \> an indeterminant over $\BF_q$,\\[1mm]
$\BF_q\dbl\zeta\dbr$ \> the ring of formal power series in $\zeta$,\\[1mm]
$\Spf\BF_q\dbl\zeta\dbr$ \> the formal scheme which is the formal spectrum of $\BF_q\dbl\zeta\dbr$,\\[1mm]
$\Nilp_{\BF_q\dbl\zeta\dbr}$ \> \parbox[t]{.79\textwidth}{the category of schemes over $\Spf\BF_q\dbl\zeta\dbr$, viz.\ the category of schemes over $\Spec\BF_q\dbl\zeta\dbr$ on which $\zeta$ is locally nilpotent.}
\end{tabbing}
From Section~\ref{SectZDivGps} on all schemes will  be in $\Nilp_{\BF_q\dbl\zeta\dbr}$.

Let $S$ be a scheme. We denote by
\begin{tabbing}
$A=\Gamma(C',\CO_C)$ \es\=\kill
$\sigma_S:S\to S$ \> \parbox[t]{.79\textwidth}{its Frobenius endomorphism which acts as the identity on points and as the $q$-power map on the structure sheaf,}\\[1mm]
$C_S = C\times S$\\[1mm]
$\sigma = \id_C\times\sigma_S$ \> \parbox[t]{.79\textwidth}{the endomorphism of $C_S$ that acts as the identity on the coordinates of $C$ and as $b\mapsto b^q$ on the elements $b\in \CO_S$.}
\end{tabbing}

\noindent
For a divisor $D$ on $C$ we denote by $\CO_{C_S}(D)$ the invertible sheaf on $C_S$ whose sections have divisor $\geq -D$. If $\CF$ is a coherent sheaf on $C_S$ we set $\CF(D) := \CF\otimes_{\CO_{C_S}}\CO_{C_S}(D)$. This notation applies in particular to the divisor $D=n\cdot\infty$ for an integer $n$.

%
%

\section*{{\Large Part One \es Abelian Sheaves}}  \label{PartI}

\mtaddtocont{\protect\contentsline {mtchap}{\protect\large Part One \es Abelian Sheaves}{\thepage}}

%
%

\section{Definition of Abelian Sheaves}  \label{SectAbelianSheaves}
\setcounter{equation}{0}

Let $S$ be a scheme and fix a morphism $c:S\to C$. Let $\CJ$ be the ideal sheaf on $C_S$ of the graph of $c$.

\begin{definition} \label{DefAbelianSheaf}
An \emph{abelian sheaf $\ul{\CF}=(\CF_i,\Pi_i,\tau_i)$ of rank $r$, dimension $d$, and characteristic $c$ over $S$\/} is a ladder of vector bundles $\CF_i$ on $C_S$ of rank $r$ and injective homomorphisms $\Pi_i$, $\tau_i$ of $\CO_{C_S}$-modules ($i\in \BZ$) of the form
\[
\begin{CD}
\cdots & @>>> & \CF_{i-1} & @>{\Pi_{i-1}}>> & \CF_i & @>{\Pi_i}>> & \CF_{i+1} & @>{\Pi_{i+1}}>> & \cdots \\
& & & & @AA{\tau_{i-2}}A & & @AA{\tau_{i-1}}A & & @AA{\tau_i}A & & \\
\cdots & @>>> & \sigma^\ast\CF_{i-2} & @>{\sigma^\ast\Pi_{i-2}}>> & \sigma^\ast\CF_{i-1} & @>{\sigma^\ast\Pi_{i-1}}>> & \sigma^\ast\CF_i & @>{\sigma^\ast\Pi_i}>> & \cdots 
\end{CD}
\]
subject to the following conditions (for all $i\in \BZ$):
\begin{enumerate}
\item \label{DefAbelianSheafCond1}
the above diagram is commutative,
\item \label{DefAbelianSheafCond2}
the morphism $\Pi_{i+\ell-1}\circ\ldots\circ\Pi_i$ identifies $\CF_i$ with the subsheaf $\CF_{i+\ell}(-k\cdot\infty)$ of $\CF_{i+\ell}$,
\item\label{DefAbelianSheafCond3}
the cokernel of $\Pi_i$ is a locally free $\CO_S$-module of rank $d$,
\item\label{DefAbelianSheafCond4}
the cokernel of $\tau_i$ is a locally free $\CO_S$-module of rank $d$ and annihilated by $\CJ^d$.
\end{enumerate}
A \emph{morphism} between two abelian sheaves $(\CF_i,\Pi_i,\tau_i)$ and $(\CF'_i,\Pi'_i,\tau'_i)$ is a collection of morphisms $\CF_i\to\CF'_i$ which commute with the $\Pi$'s and the $\tau$'s. 
\end{definition}

Let us make a few remarks. By condition~\ref{DefAbelianSheafCond2} the cokernel of $\Pi_i$ is supported at $\infty$. 
Moreover, due to the periodicity condition~\ref{DefAbelianSheafCond2} we have $\tau_{i+\ell n}=\tau_i\otimes\id_{\CO_{C_S}(kn)}$ for all $n\in\BZ$. 
Finally, the reader should be aware that we allow that the ideal sheaf $\CJ$ acts non-trivially on $\coker\tau_i$. In this respect our abelian sheaves are more general than the abelian sheaves studied so far in the literature. As such we like to mention elliptic sheaves \cite{Drinfeld3,Blum-Stuhler} which we discuss later, and $\CalD$-elliptic sheaves \cite{LRS} which are abelian sheaves equipped with an action of an order $\CalD$ in a central division algebra over $Q$.

\smallskip

\begin{definition} \label{DefStack1}
We denote by $\AbSh^{r,d}(S)$ the category whose objects are the abelian sheaves of rank $r$ and dimension $d$ over $S$ and whose morphisms are the isomorphisms of abelian sheaves. If $S'\to S$ is a morphism of schemes the pullback of an abelian sheaf over $S$ is an abelian sheaf over $S'$. This defines a fibered category $\AbSh^{r,d}$ over the category of $\BF_q$-schemes, which is a stack for the $fppf$-topology%
. The functor which assigns to an object of $\AbSh^{r,d}(S)$ the characteristic $c:S\to C$ defines a 1-morphism of stacks
\[
\AbSh^{r,d} \to C\,.
\]
\end{definition}

\medskip

Next we introduce level structures on abelian sheaves. Let $I\subset C'=C\setminus\infty$ be a finite closed subscheme and let $\ul{\CF}=(\CF_i,\Pi_i,\tau_i)$ be an abelian sheaf of rank $r$ over $S$. Then the restrictions $\CF_i|_{I\times S}$ are all isomorphic via the morphisms $\Pi_i$. We call this restriction $\ul{\CF}|_{I\times S}$. The same holds for the morphisms $\tau_i$. So we obtain a morphism 
\[
\tau|_{I\times S}:\es \sigma^\ast\ul{\CF}|_{I\times S} \longto \ul{\CF}|_{I\times S}
\]
which we consider as a $\sigma$-linear map of $\ul{\CF}|_{I\times S}$ to itself. In this article we always assume that the characteristic $c(S)$ of $(\CF_i,\Pi_i,\tau_i)$ is disjoint from $I$. 
Due to this assumption, $\tau|_{I\times S}$ is an isomorphism. We consider the \emph{functor of $\tau$-invariants of $\ul{\CF}|_{I\times S}$\/}
\[
\begin{array}{cccc}
(\ul{\CF}|_I)^\tau : & \es \Sch_S \es & \longto & \CO_I-\text{modules} \\
 & T/S & \longmapsto & \es \ker \Koh^0(I\times T, \, \tau|_{I\times T} - \id_{\ul{\CF}|_{I\times T}})\,.
\end{array}
\]
In B{\"o}ckle--Hartl \cite[Theorem 2.5]{BH} the following fact is proved.

\begin{proposition} \label{TauInvRepresentable}
The functor $(\ul{\CF}|_I)^\tau$ is representable by a finite {\'e}tale scheme over $S$ which is a $\GL_r(\CO_I)$-torsor.
\end{proposition}

\begin{definition} \label{DefLevelIStructure}
An \emph{$I$-level structure on $(\CF_i,\Pi_i,\tau_i)$ over $S$\/} is an isomorphism 
\[
{\bar\eta}:\es\bigl(\ul{\CF}|_{I\times S},\tau|_{I\times S}\bigr)\es\isoto\es \bigl(\CO_{I\times S}^r,\Id_r\cdot\sigma^\ast\bigr)
\]
from $\ul{\CF}|_{I\times S}$ to $\CO_{I\times S}^r$ that commutes with the $\sigma$-linear endomorphisms $\tau|_{I\times S}$ on one and $\Id_r\cdot\sigma^\ast$ on the other side.
\end{definition}

If the characteristic $c(S)$ meets $I$ then $\ul{\CF}$ does not possess any $I$-level structure. 

\begin{proposition} \label{PropAutomGpFinite}
Let $\ul\CF$ be an abelian sheaf over a field. Then the automorphism group of $\ul\CF$ is finite.
\end{proposition}

\begin{proof}
We may assume that the base field is algebraically closed. Let $\bar\eta$ be an $I$-level structure on $\ul\CF=(\CF_i,\Pi_i,\tau_i)$. It induces a group homomorphism
\[
\alpha_I:\Aut(\ul\CF)\;\to\; \Aut\bigl(\CO_{I\times S}^r,\Id_r\cdot\sigma^\ast\bigr) \;=\;\GL_r(\CO_I)\,.
\]
The later group is finite. We claim that $\alpha_I$ is injective for some sufficiently large finite subscheme $I\subset C'$. From this the proposition will follow. To establish the claim let $(f_i:\CF_i\to\CF_i)_i$ be an automorphism of $\ul\CF$. Note that if $f_0$ is the identity then $f_i$ must also be the identity for all $i$. We now consider a finite flat morphism $\pi:C\to\BP^1_{\BF_q}$. Since the map $\pi_\ast:\Aut_C(\CF_0)\to\Aut_{\BP^1}(\pi_\ast\CF_0)$ is injective we may assume $C=\BP^1$. Then the vector bundle $\pi_\ast\CF_0$ decomposes
\[
\pi_\ast\CF_0 \es = \es \bigoplus_{i=1}^s \CO_{\BP^1}(n_i)
\]
for uniquely determined integers $n_1\ge\ldots\ge n_s$\,. We let $P\in\BP^1$ be a point and set $I= (n_1-n_s+1)\cdot P$. Then $\alpha_I$ is injective.
\qed
\end{proof}

\medskip

Next we want to give a different definition of $I$-level structures. Note that via the natural isomorphism $(\ul{\CF}|_I)^\tau \otimes_{\CO_I}\CO_{I\times S} \isoto \ul{\CF}|_{I\times S}$ the $I$-level structures ${\bar\eta}$ on $\ul{\CF}$ over a connected $S$ correspond bijectively to the isomorphisms of $\CO_I$-modules
\[
{\bar\eta}':\es (\ul{\CF}|_I)^\tau(S) \rbij \CO_I^r\,.
\]

We use this observation to define more general level structures by introducing the adelic point of view. Let $\ul{\CF}$ be an abelian sheaf of rank $r$ over $S$. We define the functor
\[
\begin{array}{cccc}
(\ul{\CF}|_{\wh{A}})^\tau : & \es \Sch_S \es & \longto & \es \wh{A}-\text{modules} \\[2mm]
 & T/S & \longmapsto & \es \invlim[I] \,(\ul{\CF}|_I)^\tau (T)\,,
\end{array}
\]
where the limit is taken over all finite closed subschemes $I\subset C'$.
Assume that $S$ is connected and choose an algebraically closed base point $\iota:s\to S$. Due to Proposition~\ref{TauInvRepresentable} we may view $(\ul{\CF}|_{\wh{A}})^\tau$ as the $\wh{A}[\pi_1(S,s)]$-module $(\iota^\ast\ul{\CF}|_{\wh{A}})^\tau(s)$.

We consider the set $\Isom_{\wh{A}}\bigl((\ul{\CF}|_{\wh{A}})^\tau,\wh{A}^{\,r}\bigr):=\Isom_{\wh{A}}\bigl((\iota^\ast\ul{\CF}|_{\wh{A}})^\tau(s),\wh{A}^{\,r}\bigr)$ of isomorphisms of $\wh{A}$-modules. Via its natural action on $\wh{A}^{\,r}$ the group $\GL_r(\wh{A})$ acts on this set from the left. Via its action on $(\iota^\ast\ul{\CF}|_\wh{A})^\tau(s)$ the group $\pi_1(S,s)$ acts on it from the right.

\begin{definition} \label{DefHLevelStructure}
Let $H\subset \GL_r(\wh{A})$ be a compact open subgroup. An \emph{$H$-level structure on $\ul{\CF}$ over $S$\/} is an $H$-orbit in $\Isom_{\wh{A}}\bigl((\ul{\CF}|_{\wh{A}})^\tau,\wh{A}^{\,r}\bigr)$ which is fixed by $\pi_1(S,s)$. (Because of this later condition the notion of level structure is independent of the chosen base point.)
\end{definition}

In particular if $H=H_I=\ker\bigl(\GL_r(\wh{A})\to\GL_r(\CO_I)\bigr)$ an $H$-level structure is nothing else than an $I$-level structure. Note that as before an $H$-level structure can only exist if the characteristic $c(S)$ does not meet the set of places $v$ of $C$ for which $H_v \neq \GL_r(A_v)$.

\begin{definition} \label{DefStack2}
Let $\AbSh^{r,d}_H(S)$ be the category whose objects are the abelian sheaves of rank $r$ and dimension $d$ together with an $H$-level structure over $S$ and whose morphisms are the isomorphisms of abelian sheaves which respect the level structures. Analogous to Definition~\ref{DefStack1}, this defines a stack $\AbSh^{r,d}_H$ over $C$.
\end{definition}

For $H=H_\emptyset=\GL_r(\wh{A})$ the definition of $H$-level structure is vacuous. Therefore we have $\AbSh^{r,d}_{H_\emptyset}=\AbSh^{r,d}$. We will show in Section~\ref{SectAlg} that these stacks are algebraic over $C$.

There is a free action of the group $\BZ$ on these stacks given on objects by the map
\[
[n]: \es (\CF_i,\Pi_i,\tau_i) \es \mapsto \es (\CF_{i+n},\Pi_{i+n},\tau_{i+n})\,.
\]

\begin{bigexample}[Drinfeld Modules and Elliptic Sheaves] \label{ExDrinfeldModules}
Let $d=1$ and let $H=H_I$. Then an abelian sheaf is what was called an \emph{elliptic sheaf} in Blum--Stuhler~\cite{Blum-Stuhler}. We consider the open substack
\[
\AbSh^{r,1}_H\times_C \;C'
\]
of abelian sheaves with characteristic disjoint from $\infty$. It is shown in \cite[Theorem 3.2.1]{Blum-Stuhler} that there is a 1-isomorphism between the stack $\DrMod^r_I$ of Drinfeld $A$-modules of rank $r$ with $I$-level structure and the open and closed substack of $\AbSh^{r,1}_H\times_C \;C'$ consisting of those $(\CF_i,\Pi_i,\tau_i)$ with $\deg(\CF_0|_{C_s}) = 1-r$ for each algebraically closed point $s$ of $S$.

On the other hand in $\AbSh^{r,1}_H$ we always have $\deg(\CF_i|_{C_s})= \deg(\CF_0|_{C_s})+i$. 
Using the free action of $\BZ$ on $\AbSh^{r,1}_H$ we thus obtain a 1-isomorphism of stacks
\[
\AbSh^{r,1}_H\times_C \;C' \es \cong\es \BZ \times \DrMod^r_I\,.
\]
The first factor gives the degree of $\CF_0$.

In fact if $I\neq\emptyset$ the stack $\DrMod^r_I$ is a smooth affine scheme of finite type over $\BF_q$. See for example \cite[Theorem 2.3.8]{Blum-Stuhler}. So $\AbSh^{r,1}_H\times_C \;C'$ is a smooth scheme locally of finite type in this case.

\smallskip

We will give another example in Section~\ref{SectExample} which shows that for $d>1$ one can neither expect that $\AbSh^{r,d}_H$ is a scheme, nor that it is smooth over $C'$.
\end{bigexample}

%
%

\section{Relation to Anderson's $t$-Motives}  \label{SectPolTMotives}
\setcounter{equation}{0}

We will show that abelian sheaves with characteristic disjoint from $\infty$ are the same as polarized $A$-motives, a variant of Anderson's~\cite{Anderson} $t$-motives. Let $S$ be a scheme and fix a characteristic morphism $c:S\to C'$ disjoint from $\infty$. Let $\Gamma(c)\subset C_S$ be the graph of $c$ and let $\CJ$ be the ideal sheaf defining $\Gamma(c)$.

\begin{definition} \label{DefPolTMotive}
A \emph{(pure) polarized $A$-motive $(\CF,\tau)$ of rank $r$, dimension $d$, and characteristic $c$ over $S$\/} consists of a vector bundle $\CF$ on $C_S$ of rank $r$ and a morphism of coherent sheaves $\tau:\sigma^\ast \CF \to \CF(k\cdot \infty)$ such that
\begin{enumerate}
\item \label{1.CondPolTMotive}
the cokernel of $\tau$ is supported on $\Gamma(c) \cup \infty$, the part supported on $\Gamma(c)$ is annihilated by $\CJ^d$,
\item \label{2.CondPolTMotive}
for every $i=1,\ldots,\ell$ the image of $\tau^i:\sigma^{i\,\ast}\CF \to \CF(i k\cdot\infty)$ lies in $\CF(k\cdot\infty)$ and
$\tau^\ell:\sigma^{\ell\,\ast}\CF \to \CF(k\cdot\infty)$ is a local isomorphism at $\infty$,
\item \label{3.CondPolTMotive}
locally at $\infty$, $\CF$ is contained in the image $\tau\CF$.
\end{enumerate}
\end{definition}

\noindent
Here we denote by $\tau^i$ the composition $\bigl(\tau \otimes \id_{\CO\left((i-1)k\cdot\infty\right)}\bigr)\circ \ldots\circ\sigma^{(i-1)\,\ast}\tau$ mapping
\[
\sigma^{i\,\ast}\CF \longto \sigma^{(i-1)\,\ast}\CF(k\cdot\infty) \longto \ldots \longto \sigma^\ast\CF\bigl((i-1)k\cdot\infty\bigr) \longto\CF(i k\cdot\infty)\,.
\]


\begin{lemma} \label{LemCokerLocFree}
Conditions~\ref{1.CondPolTMotive} and \ref{2.CondPolTMotive} imply that $\tau$ is injective and that the part of $\coker \tau$ which is supported on $\Gamma(c)$ (respectively on $\infty$) is a locally free $\CO_S$-module of rank $d$ (respectively $(\ell-1)d$). 
\end{lemma}

\begin{proof}
Consider the exact sequences of $\CO_S$-modules induced by $\tau$
\[ \xymatrix @R=1.5pc @C=1pc {
0\ar[r] &\ker\tau\ar[r] &\sigma^\ast\CF \ar[r]\ar[d]^\tau & \im\tau\ar[r] &0\es\\
0\ar[r] &\im\tau \ar[r] &\CF(k\cdot\infty)\ar[r] & \coker\tau\ar[r]&0\,.
}
\]
For a point $s\in S$ we tensor them with the residue field $\kappa(s)$ at $s$ to obtain
\[ \xymatrix @R=1.5pc @C=1.5pc {
0\ar[r] &\Tor^{\CO_S}_1\bigl(\im\tau,\kappa(s)\bigr)\ar[r] &(\ker\tau)_s\ar[r] &(\sigma^\ast\CF)_s \ar[r]^\beta\ar[d]^{\tau\otimes\id} & (\im\tau)_s\ar[r] &0\es\\
0\ar[r] &\Tor^{\CO_S}_1\bigl(\coker\tau,\kappa(s)\bigr)\ar[r] &(\im\tau)_s \ar[r]^{\alpha\quad} &\CF(k\cdot\infty)_s\ar[r] & (\coker\tau)_s\ar[r]&0\,.
}
\]
We claim that $\tau\otimes\id$ is injective. Indeed consider a point on $C_s$ and its local ring which is a PID. Then over this PID $\tau\otimes\id$ is a morphism between finite free modules of the same rank with torsion cokernel. Hence by the elementary divisor theorem $\tau\otimes\id=\alpha\beta$ must be injective.

Now the surjectivity of $\beta$ shows that $\alpha$ is injective. Therefore $\Tor^{\CO_S}_1\bigl(\coker\tau,\kappa(s)\bigr)$ vanishes and by the local criterion for flatness \cite[6.8]{Eisenbud}, $\coker\tau$ is a flat $\CO_S$-module. This in turn implies the flatness of $\im\tau$. Hence $(\ker\tau)_s$ is identified with the kernel of $\beta$. However, $\beta$ is an isomorphism and thus $(\ker\tau)_s$ is zero. By Nakayama's lemma $\ker\tau$ is zero and $\tau$ is injective. Finally let $\CG$ be the part of $\coker \tau$ which is supported on $\Gamma(c)$. Since the characteristic is disjoint from $\infty$, $\CG$ is a direct summand of $\coker\tau$. Then condition~\ref{2.CondPolTMotive} implies that $\coker\tau^\ell = \bigoplus_{i=0}^{\ell-1}\sigma^{i\,\ast}\CG$ and the lemma follows.
\qed
\end{proof}

If $S=\Spec K$ is the spectrum of a field, $A=\BF_q[t]$, and $(\CF,\tau)$ is a polarized $A$-motive over $S$ then the $K[t,\tau]$-module $\Gamma(C'_S, \CF)$ is a \emph{pure $t$-motive of weight $d/r$\/} as defined by Anderson \cite{Anderson}. Compared to Anderson's definition however, our \emph{polarized\/} $A$-motive contains additional data at infinity which rigidifies the structure of the $A$-motive.
Since Anderson's $t$-motives serve as characteristic $p$ analogues of abelian varieties this may justify our terminology.

As in Section~\ref{SectAbelianSheaves} we can define $H$-level structures on polarized $A$-motives for compact open subgroups $H\subset\GL_r(\wh{A})$. Correspondingly we obtain the stack $\Pol^{r,d}_H$ of polarized $A$-motives with $H$-level structure. This is a stack over $C'$.

\begin{theorem} \label{ThmIsomOfStacks}
There is a 1-isomorphism of stacks 
\[
\AbSh^{r,d}_H\times_C \;C' \es \cong \es \Pol^{r,d}_H\,.
\]
\end{theorem}

\begin{proof}
We construct two mutually 2-inverse 1-morphisms $T$ and $T'$ between the stacks in question. First consider the 1-morphism $T:\AbSh^{r,d}_H\times_C \;C' \to \Pol^{r,d}_H$ which assigns to an abelian sheaf $(\CF_i,\Pi_i,\tau_i)$ over $S$ with characteristic disjoint from $\infty$ the polarized $A$-motive of rank $r$ and dimension $d$ consisting of
\[
\CF := \CF_0\,,\qquad \tau :=\Pi_{\ell-1}\circ\ldots\circ\Pi_1\circ\tau_0: \es\sigma^\ast \CF \to \CF(k\cdot\infty)\,.
\]
Clearly this construction is also compatible with level structures.

Conversely let $(\CF,\tau)$ be a polarized $A$-motive of rank $r$ and dimension $d$ over $S$. For $0\leq i \leq \ell$ we set
\[
\CF_i \es := \es \CF + \ldots + \tau^i\CF \es \subset \es \CF(k\cdot\infty)\,.
\]
Then $\CF_i$ is equal to $\CF$ outside $\infty$ due to the definition of $\tau$. And locally at $\infty$  it is isomorphic to $\sigma^{i\,\ast}\CF$ due to condition~\ref{3.CondPolTMotive} of Definition~\ref{DefPolTMotive}. Therefore $\CF_i$ is locally free on $C_S$. 
Let $\Pi_i:\CF_i\to\CF_{i+1}$ be the inclusion and let $\tau_i:\sigma^\ast\CF_i\to\CF_{i+1}$ be the morphism $\tau$. Again by condition~\ref{3.CondPolTMotive}, $\coker\Pi_i$ is supported at $\infty$ and $\coker\tau_i$ is supported on $\Gamma(c)$ and annihilated by $\CJ^d$. By Lemma~\ref{LemCokerLocFree} these cokernels are locally free $\CO_S$-modules of rank $d$. Furthermore, by condition~\ref{2.CondPolTMotive} we have $\CF_\ell = \CF(k\cdot\infty)$. For arbitrary $i\in \BZ$ we take $n\in \BZ$ such that $0\leq i-n\ell <\ell$ and define
\[
\CF_i := \CF_{i-n\ell} \otimes_{\CO_{C_S}} \CO_{C_S}(nk\cdot\infty)\,,\quad \Pi_i := \Pi_{i-n\ell} \otimes \id\,,\quad \tau_i :=\tau_{i-n\ell} \otimes \id\,.
\]
Then $(\CF_i,\Pi_i,\tau_i)$ is an abelian sheaf of rank $r$ and dimension $d$ over $S$. This construction is also compatible with level structures and defines a 1-morphism $T': \Pol^{r,d}_H \to \AbSh^{r,d}_H\times_C \;C'$. One easily proves that $T$ and $T'$ are mutually 2-inverse.
\qed
\end{proof}

\begin{remark}
Example~\ref{ExDrinfeldModules} together with Theorem~\ref{ThmIsomOfStacks} shows that every Drinfeld module carries a canonical polarization. This parallels the situation for elliptic curves.
\end{remark}

%
%

\section{Algebraicity of the Stacks of Abelian Sheaves}  \label{SectAlg}
\setcounter{equation}{0}

In this section we will prove the following theorem.

\begin{theorem} \label{ThmAlgebraicStacks}
$\AbSh^{r,d}_H$ is an algebraic stack in the sense of Deligne--Mumford \cite{DM}, locally of finite type over $C$.
\end{theorem}
%

Note that the example in Section~\ref{SectExample} below shows that it will in general not be smooth over $C$. Nevertheless, if the characteristic is disjoint from $\infty$ and one stratifies the stack according to the isomorphy type of the $\CO_C$-modules $\coker\Pi_i$ and $\coker\tau_i$ then each stratum will be smooth over $C'$. 

More precisely, we let $z$ be a uniformizing parameter on $C$ at $\infty$ and we fix a flag $\CG_\bullet$ of $\BF_q[z]/z^k$-submodules
\begin{equation} \label{EqFixedFlag}
(0)\subset \CG_1\subset\ldots\subset \CG_{\ell-1} \subset \CG_\ell = (\BF_q[z]/z^k)^{\oplus r}
\end{equation}
such that the successive quotients all have dimension $d$ over $\BF_q$.
We say that an abelian sheaf $\ul\CF$ of rank $r$ and dimension $d$ over $S$ has \emph{isomorphy type $\CG_\bullet$ of $\Pi$} if for every point $s\in S$ the flag of $\CO_{C_s}$-modules
\[
(0)\subset \CF_1/\CF_0\subset\ldots\subset\CF_{\ell-1}/\CF_0 \subset \CF_{\ell}/\CF_0 = \CF_0(k\cdot\infty)/\CF_0
\]
is isomorphic to the flag $\CG_\bullet\otimes_{\BF_q}\kappa(s)$.

To fix the isomorphy type of $\tau$ let $S$ be a scheme with characteristic morphism $c:S\to C$. We denote by $\CJ\subset\CO_{C_S}$ the ideal sheaf defining the graph of $c$, by $\Gamma_d$ the closed subscheme of $C_S$ defined by $\CJ^d$, and by $\CO_d$ the structure sheaf $\CO_{C_S}/\CJ^d$ of $\Gamma_d$. We fix integers 
\begin{equation} \label{EqElemDiv}
( e_1\leq e_2\leq \ldots\leq e_r)\es=\es\ul{e}
\end{equation}
between $0$ and $d$ with $e_1+\ldots+e_r=d$. 
We say that an abelian sheaf $\ul\CF$ of rank $r$ and dimension $d$ over $S$ has \emph{has isomorphy type $\ul{e}$ of $\tau$} if for every point $s\in S$ the $\CO_{C_s}$-module $\coker\tau_{-1}$ is isomorphic to
\[
\bigoplus_{\nu=1}^r \CO_{C_s}/\,\CJ^{e_\nu}\CO_{C_s}\,.
\]
Note that if the characteristic is disjoint from $\infty$ then the cokernels of all $\tau_i$ are isomorphic.

We consider the locally closed subset (see Laumon--Moret-Bailly~\cite[{\S}5]{LaumonMB}) of $\AbSh^{r,d}_H$ consisting of the points over which the universal abelian sheaf has isomorphy types $\CG_\bullet$ of $\Pi$ and $\ul e$ of $\tau$.
We give this subset the reduced induced structure (\cite[4.10]{LaumonMB}) and obtain a substack $\AbSh^{r,d,\CG_\bullet,\ul{e}}_H$ of $\AbSh^{r,d}_H$.

\begin{theorem} \label{ThmSmoothStratification}
The substack $\AbSh^{r,d,\CG_\bullet,\ul{e}}_H \times_C\,C'$
is smooth over $C'$ of relative dimension $\sum_{\mu>\nu}(e_\mu-e_\nu)$ if non-empty.
\end{theorem}
%

The proof of these theorems follows Laumon--Rapoport--Stuhler \cite{LRS}. For a given compact open subgroup $H\subset\GL_r(\wh A)$ consider a closed subscheme $I\subset C'$ which is supported on the places $v$ for which $H_v\ne\GL_r(A_v)$ and satisfies 
\[
H\supset H_I:=\ker\bigl(\GL_r(\wh{A})\to\GL_r(\CO_I)\bigr)\,.
\]
Then $\AbSh^{r,d}_{H_I}$ is finite {\'e}tale over $\AbSh^{r,d}_H$ by Proposition~\ref{TauInvRepresentable}. In fact it is even a $H/H_I$-torsor. Hence $\AbSh^{r,d}_H$ is a quotient of $\AbSh^{r,d}_{H_I}$ in the sense of stacks; cf.\ \cite[4.6.1]{LaumonMB}. Therefore it suffices to prove the two theorems for $H=H_I$. We assume this situation from now on.

We will cover $\AbSh^{r,d}_H$ by open substacks corresponding to stable vector bundles with additional level structure. To be precise we proceed as follows.
Let $\CF$ be a locally free sheaf of rank $r$ on $C_S$. An \emph{$I$-level structure on $\CF$} is an isomorphism ${\bar\eta}:\CF|_{I\times S} \isoto \CO_{I\times S}^r$ of $\CO_{I\times S}$-modules. 

\begin{definition} \label{DefStableVb}
We say that the pair $(\CF,{\bar\eta})$ is \emph{stable} if for all algebraically closed points $s\in S$ and all locally free $\CO_{C_s}$-modules $\CG$ properly contained in $\CF_s$, we have
\[
\frac{\deg(\CG) - \deg(I)}{\rk(\CG)} \es < \es \frac{\deg(\CF_s) - \deg(I)}{\rk(\CF_s)}
\]
(compare Seshadri \cite[4.I. D{\'e}finition 2]{Seshadri}). 
\end{definition}

We denote by $\AbSh^{r,d}_{H,st}$ the open substack of $\AbSh^{r,d}_H$ of those abelian sheaves for which $(\CF_0,{\bar\eta})$ is stable. We will show that $\AbSh^{r,d}_{H,st}$ is representable by a disjoint union of quasi-projective schemes of relative dimension $d(r-1)$ over $C\setminus I$ if $I\neq\emptyset$. For this purpose we need to introduce some additional stacks.

We denote by $\Vect^r_I$ the stack classifying locally free sheaves of rank $r$ on $C$ with $I$-level structure and by $\Vect^r_{I,st}$ the substack of such locally free sheaves which are stable. Seshadri \cite[4.III]{Seshadri} proves the following fact (see also \cite[4.3]{LRS}). Note that there is a typing error in the formula for the dimension in \cite{Seshadri}.

\begin{proposition}
If $I\neq\emptyset$ the stack $\Vect^r_{I,st}$ is representable by a disjoint union of quasi-projective schemes over $\BF_q$ which are smooth of dimension $r^2(g-1+\deg I)$. Here $g$ is the genus of $C$.
\end{proposition}

\begin{definition} \label{DefStackVec}
Let $S$ be a scheme and let $(\CF_i,\Pi_i)$ be a sequence of locally free sheaves on $C_S$ as in the first row of the Definition~\ref{DefAbelianSheaf} of an abelian sheaf. Suppose that $(\CF_i,\Pi_i)$ satisfies the conditions~\ref{DefAbelianSheafCond2} and \ref{DefAbelianSheafCond3} of Definition~\ref{DefAbelianSheaf}. An $I$-level structure on $(\CF_i,\Pi_i)$ is a collection of $I$-level structures ${\bar\eta}_i$ on the $\CF_i$ that are compatible with the morphisms $\Pi_i$.
We denote by $\Seq^{r,d}_I$ the stack classifying sequences $(\CF_i,\Pi_i)$ as above together with $I$-level structures ${\bar\eta}_i$ and by $\Seq^{r,d}_{I,st}$ the open substack of those sequences for which $\CF_0$ with its level structure is stable.
\end{definition}

\begin{lemma} \label{LemReprSeq}
The natural 1-morphism 
\[
\Seq^{r,d}_I \to\Vect^r_I\,,\quad(\CF_i,\Pi_i,{\bar\eta}_i)\mapsto (\CF_0,{\bar\eta}_0)
\]
is representable by a closed subscheme of a flag variety.
The strata $\Seq^{r,d,\CG_\bullet}_I$ with fixed isomorphy type $\CG_\bullet$ of $\Pi$ are smooth over $\Vect^r_I$ if non-empty. In particular $\Seq^{r,d}_{I,st}$ and $\Seq^{r,d,\CG_\bullet}_{I,st}$ are representable by a disjoint union of quasi-projective schemes (respectively quasi-projective and smooth schemes) over $\BF_q$ if $I\neq\emptyset$.
\end{lemma}

\begin{proof}
Let $\CF_0$ on $C_S$ be given corresponding to a 1-morphism $S\to\Vect^r_I$. Due to the periodicity condition~\ref{DefAbelianSheafCond2} the sequence $(\CF_i,\Pi_i)$ corresponds to a flag of length $\ell$ of $\CO_{C_S}$-submodules
\begin{equation} \label{EqChainAtInfty}
(0)\subset \CF_1/\CF_0\subset\ldots\subset\CF_{\ell-1}/\CF_0 \subset \CF_{\ell}/\CF_0 = \CF_0(k\cdot\infty)/\CF_0
\end{equation}
such that the successive quotients are all locally free $\CO_S$-modules of rank $d$. Hence the first assertion follows.

To prove the statement about $\Seq^{r,d,\CG_\bullet}_I$ note that locally on $S$ the sheaf $\CF_\ell/\CF_0$ is isomorphic to $\CO_S[z]/z^k$. By the elementary divisor theorem a point $\Spec K$ of $\Seq^{r,d}_I$ belongs to $\Seq^{r,d,\CG_\bullet}_I$ if and only if the flag (\ref{EqChainAtInfty}) is conjugate to $\CG_\bullet$ under $\GL_r\bigl(K[z]/z^k\bigr)$. Therefore $\Seq^{r,d,\CG_\bullet}_I$ is relatively representable over $\Vect^r_I$ by the homogeneous space
\[
\GL_r(\BF_q[z]/z^k)/\Stab(\CG_\bullet)\,.
\]
The group $\GL_r(\BF_q[z]/z^k)$ is the Weil restriction $\CR_{(\BF_q[z]/z^k)/\BF_q}\GL_r$ and hence a smooth connected algebraic group over $\BF_q$. The stabilizer of $\CG_\bullet$ corresponds to a closed algebraic subgroup defined over $\BF_q$. Thus the above homogeneous space is a smooth algebraic variety over $\BF_q$. From this the lemma follows.
\qed
\end{proof}

\begin{definition}\label{DefStackHecke}
We let $\Hecke^{r,d}_I$ be the stack classifying the commutative diagrams with $I$-level structures
\[
\begin{CD}
\cdots & @>>> & \CF_{i-1} & @>{\Pi_{i-1}}>> & \CF_i & @>{\Pi_i}>> & \CF_{i+1} & @>{\Pi_{i+1}}>> & \cdots \\
& & & & @AA{t_{i-2}}A & & @AA{t_{i-1}}A & & @AA{t_i}A & & \\
\cdots & @>>> & \CF'_{i-2} & @>{\Pi'_{i-2}}>> & \CF'_{i-1} & @>{\Pi'_{i-1}}>> & \CF'_i & @>{\Pi'_i}>> & \cdots 
\end{CD}
\]
such that the sequences $(\CF_i,\Pi_i)$ and $(\CF'_i,\Pi'_i)$ with their $I$-level structures belong to $\Seq^{r,d}_I$ and such that the $t_i$ satisfy conditions~\ref{DefAbelianSheafCond1} and \ref{DefAbelianSheafCond4} of Definition~\ref{DefAbelianSheaf} and respect the $I$-level structures. Assigning to such a diagram over $S$ the morphism $S\to C\setminus I$ on whose graph the cokernels of the $t_i$ are supported, defines a 1-morphism of stacks $\Hecke^{r,d}_I\to C\setminus I$.
\end{definition}

The above stacks fit into the following 2-cartesian diagram of stacks
\begin{equation} \label{DiagramStacks}
\begin{CD}
\AbSh^{r,d}_H & @>>> & \Seq^{r,d}_I \\
@VVV & & @VV(\id,\sigma_{\Seq})V \\
\Hecke^{r,d}_I & @>(\rm 1^{st}\,row,\,2^{nd}\,row)>> & \Seq^{r,d}_I \times \Seq^{r,d}_I\\
@VVV \\
C\setminus I
\end{CD}  
\end{equation}

On the stack $\Hecke^{r,d}_I$ the cokernel of $t_{-1}$ is a quotient of $\CF_0$ which is locally free of rank $d$ over the base and supported on the graph of the characteristic morphism $c$. We analyze this property. Let $\CT$ be the stack
\[
(C\setminus I)\times\Seq^{r,d}_I\,.
\]
On $C\times\CT$ consider the ideal sheaf $\CJ$ defining the graph $\Gamma(c)$ of the characteristic morphism $c:\CT\to C\setminus I\subset C$ (see \cite[{\S}12]{LaumonMB}). We denote by $\Gamma_d$ the closed substack of $C\times\CT$ defined by $\CJ^d$ and by $\CO_d$ the sheaf $\CO_{C\times\CT}/\CJ^d$. The quotients of $\CF_0$ that are supported on $\Gamma(c)$ and are locally free over $\CT$ of rank $d$ are classified by  Grothendieck's Quot-scheme relative to $\CT$ \cite[$\rm n^o$ 221, Th{\'e}or{\`e}me 3.1]{FGA} 
\[
\Quot^d_{\CF_0\otimes \CO_d/\Gamma_d/\CT}\,.
\]
It is a stack projective over $\CT$.

\begin{lemma} \label{LemHecke}
The 1-morphism
\[
\Hecke^{r,d}_I\to \Quot^d_{\CF_0\otimes \CO_d/\Gamma_d/\CT} \times \Seq^{r,d}_I
\]
given by the cokernel of $t_{-1}$ and the second row, is representable by a closed immersion. Over $C\setminus(I\cup\infty)$ the 1-morphism
\[
\Hecke^{r,d}_I\to \Quot^d_{\CF_0\otimes \CO_d/\Gamma_d/\CT} 
\]
obtained by projection onto the first factor is a 1-isomorphism. In particular, the 1-morphism $\Hecke^{r,d}_I\to \Seq^{r,d}_I \times \Seq^{r,d}_I$ from (\ref{DiagramStacks}) is representable by a quasi-projective morphism.
\end{lemma}

\begin{proof}
The substack $\Hecke^{r,d}_I$ of $\Quot^d_{\CF_0\otimes \CO_d/\Gamma_d/\CT} \times \Seq^{r,d}_I$ is defined by the following conditions 
\begin{enumerate}
\item \label{LemHeckeCond1}
$\CF'_{-1}$ equals the kernel of the morphism from $\CF_0$ to the universal quotient,
\item \label{LemHeckeCond2}
for each $i= -\ell,\ldots,-2$ the sheaf $\CF'_i$ is contained in the intersection of $\CF_{i+1}$ and $\CF'_{-1}$ which we view as subsheaves of $\CF_0$ via $\Pi_{-1}\circ\ldots\circ\Pi_{i+1}$ and $t_{-1}$,
\item \label{LemHeckeCond3}
if we let $t_i$ be the inclusion $\CF'_i\subset\CF_{i+1}$ then $\coker t_i$ is annihilated by $\CJ^d$,
\item \label{LemHeckeCond4}
$t_{-1}$ is compatible with the $I$-level structures on $\CF'_{-1}$ and $\CF_0$.
\end{enumerate}
Namely by descending induction on $i$ the short exact sequences of $\CO_S$-modules
\[ \xymatrix @R=1pc {
0\ar[r] &\coker t_{i-1}\ar[r] &\coker(\Pi_i\circ t_{i-1})\ar[r]\ar@{=}[d] &\coker\Pi_i\ar[r] & 0\\
0\ar[r] &\coker \Pi'_{i-1}\ar[r] &\coker(t_i\circ \Pi'_{i-1})\ar[r] &\coker t_i\ar[r] & 0
}
\]
imply that $\coker t_i$ is a locally free $\CO_S$-module of rank $d$.
Now clearly the above conditions are represented by a closed immersion.

Over $C\setminus(I\cup\infty)$ defining $\CF'_i$ as the intersection $\CF_{i+1}\cap\CF'_{-1}$ in conditions~\ref{LemHeckeCond1} and \ref{LemHeckeCond2} automatically gives an object $(\CF'_i,\Pi'_i)$ in $\Seq^{r,d}_I$. This proves that the projection onto the first factor is a 1-isomorphism there.
\qed
\end{proof}

\smallskip

\begin{proof}[of Theorem~\ref{ThmAlgebraicStacks}]
Recall that we have assumed $H=H_I$.
Considering the diagram (\ref{DiagramStacks}) we conclude from the previous lemmas that $\AbSh^{r,d}_{H,st}$ is representable by a disjoint union of quasi-projective schemes over $C\setminus I$ if $I\ne\emptyset$.

Now we let $I\subset I'\subset C'$ be two finite closed subschemes with $I'\neq \emptyset$ and we set $H'=H_{I'}$. By restricting $I'$-level structures to $I$-level structures we obtain a 1-morphism of stacks
\[
r_{I',I}:\es\AbSh^{r,d}_{H'}\es\to\es\AbSh^{r,d}_H\,.
\]
Over $C\setminus I'$ this 1-morphism is a torsor under the finite group
\[
G_{I',I} \es :=\es\ker\bigl(\GL_r(\CO_{I'})\to\GL_r(\CO_I)\bigr)\es\cong\es H/H'
\]
due to Proposition~\ref{TauInvRepresentable}. Since $r_{I',I}^{-1}(\AbSh^{r,d}_{H,st})\subset\AbSh^{r,d}_{H',st}$, the open substack $\AbSh^{r,d}_{H',st}$ which is stable under $G_{I',I}$, gives as a quotient in the sense of stacks an open substack
\[
\AbSh^{r,d}_{H',st}\,/\,G_{I',I}\es\subset\es \AbSh^{r,d}_H
\]
that contains $\AbSh^{r,d}_{H,st} \times_{C\setminus I}\;C\setminus I'$. It is an algebraic stack in the sense of Deligne--Mumford. If we let $I'$ vary among finite closed subschemes of $C'$ containing $I$ these open substacks cover $\AbSh^{r,d}_H$, since every vector bundle becomes stable for a sufficiently high level structure. This proves Theorem~\ref{ThmAlgebraicStacks} except for the assertion on the dimension which follows from Theorem~\ref{ThmSmoothStratification}.
\qed
\end{proof}

\begin{proof}[of Theorem~\ref{ThmSmoothStratification}]
We denote by $\Quot^\ul{e}$ the reduced, locally closed substack of $\Quot^d_{\CF_0\otimes \CO_d/\Gamma_d/\CT}$ consisting of those points for which the universal quotient is isomorphic to
\[
\bigoplus_{\nu=1}^r \CO_{C_s}/\,\CJ^{e_\nu}\CO_{C_s}\,.
\]
We claim that $\Quot^\ul{e}$ is smooth over $\CT$ of relative dimension $\sum_{\mu>\nu}(e_\mu-e_\nu)$. 
Indeed, locally on $\CT$ the sheaf $\CF_0\otimes\CO_d$ is isomorphic to $\CO_d^r$. Let $\CH\subset \CO_d^r$ be the kernel of the morphism from $\CO_d^r$ to the universal quotient on $\Quot^\ul{e}$. The condition on the isomorphy type of the universal quotient implies that $\CH$ is conjugate to $\oplus_{\nu=1}^r \CJ^{e_\nu}/\CJ^d$ under $\GL_r(\CO_d)$. Therefore $\Quot^\ul{e}$ is locally isomorphic to the homogeneous space
\[
\GL_r(\CO_d)\,/\,\Stab(\oplus_{\nu=1}^r \CJ^{e_\nu}/\CJ^d)\,.
\]
As in the proof of Lemma~\ref{LemReprSeq} this homogeneous space is smooth and the claim follows. The theorem can now be deduced using Lemma~\ref{LemHecke} and applying \cite[Lemma 4.2]{LRS} to diagram (\ref{DiagramStacks}).
\qed
\end{proof}

%
%

\section{An Example} \label{SectExample}
\setcounter{equation}{0}

In this section let $C=\BP^1_{\BF_q}$ and $A=\BF_q[t]$. Let $I=\Var(t)\subset C$ and $H=H_I$ and set $z=\frac{1}{t}$. Then $C\setminus I=\Spec\BF_q[z]$.
We consider the case where $d=r=2$, $k=\ell=1$ and describe the algebraic stack $\AbSh^{2,2}_H$. It decomposes
\[
\AbSh^{2,2}_H\es = \es\coprod_{n\in\BZ}\;\AbSh^{2,2}_H(n)
\]
into the open and closed substacks on which the vector bundle $\CF_0$ has degree $n$. The shift by $1$ from \ref{DefStack2} yields a 1-isomorphism $\AbSh^{2,2}_H(n)\to\AbSh^{2,2}_H(n+2)$. So it suffices to describe $\AbSh^{2,2}_H(n)$ for $n=0,1$. We want to treat the case $n=0$ here.

Let $M_I^{2,2}$ be the scheme
\[
\Spec \BF_q[\zeta] [ a_{\mu\nu}: 1\leq \mu,\nu\leq2]\,/\,(a_{11}+a_{22}+2\zeta\,,\,a_{11}a_{22}-a_{12}a_{21}-\zeta^2)\,.
\]
We view $(a_{\mu\nu})$ as a $2\times 2$ matrix with trace $-2\zeta$ and determinant $\zeta^2$. Mapping $z$ to $\zeta$ defines a morphism $c:M_I^{2,2}\to C$.
On $S=M_I^{2,2}$ we set for $i\in \BZ$
\[
\CF_i=\CO_{C_S}(i\cdot\infty)^{\oplus 2} \quad\text{and}\quad \tau_i=\bigl(1+t(a_{\mu\nu})\bigr)\cdot\sigma^\ast:\;\sigma^\ast\CF_i\to\CF_{i+1}\quad\text{and}
\]
we let $\Pi_i:\CF_i\to\CF_{i+1}$ be the morphism induced by the inclusion $\CO_{C_S}\subset\CO_{C_S}(\infty)$.
Due to the trace and determinant condition on the matrix $(a_{\mu\nu})$, the data $(\CF_i,\Pi_i,\tau_i)$ is an abelian sheaf of rank $2$, dimension $2$, and characteristic $c$ over $S$. 

We want to define an $I$-level structure on $\ul{\CF}=(\CF_i,\Pi_i,\tau_i)$. Note that $\ul{\CF}|_{I\times S}$ is canonically isomorphic to $\CO_S^2$ with $\tau|_{I\times S} = \Id_2\cdot\sigma^\ast$. Hence
the identity morphism on $\CO_S^2$ defines an $I$-level structure ${\bar\eta}$ on $\ul{\CF}$.

\begin{proposition} \label{PropExr=d=2}
The 1-morphism $M^{2,2}_I\to\AbSh^{2,2}_H$ induced by $(\ul{\CF}, {\bar\eta})$ identifies $M^{2,2}_I$ with the (representable) open substack of $\AbSh^{2,2}_H$ on which the underlying vector bundle with level structure $(\CF_0,{\bar\eta})$ is stable (Definition~\ref{DefStableVb}) and has degree zero.
\end{proposition}
We will see below that $(\CF_0,{\bar\eta})$ is stable and of degree zero if and only if $\CF_0 \cong \CO_{C_S}^2$.

\begin{proof}
Let $T$ be a scheme together with a characteristic morphism $c':T \to \Spec \BF_q[z]$ disjoint from $I$. Denote the image $c'{}^\ast (z)$ in $\CO_T$ by $\zeta'$. Let $(\CF'_i,\Pi'_i,\tau'_i,{\bar\eta}')$ be an abelian sheaf of rank $2$, dimension $2$, and characteristic $c'$ over $T$ with $I$-level structure such that $(\CF'_0,{\bar\eta}')$ is stable and of degree zero. We have to exhibit a uniquely defined morphism $f:T\to M^{2,2}_I$ such that $f^\ast(\ul{\CF},{\bar\eta}) \cong (\CF'_i,\Pi'_i,\tau'_i,{\bar\eta}')$. Since the morphisms $\Pi'_i$ identify $\CF'_i$ with $\CF'_0(i\cdot\infty)$ it suffices to concentrate on $\CF'_0$ and $\tau'_0$.

We claim that the stability condition implies $\CF'_0\cong\CO_{C_T}^2$ globally on $T$. Indeed let $\pi:C_T\to T$ be the projection onto the second factor. We first show that $\pi_\ast\CF'_0$ is locally free of rank 2 on $T$ and that $\pi^\ast\pi_\ast\CF'_0\to\CF'_0$ is an isomorphism. Let $s\in T$ be an algebraically closed point. The stability implies that every invertible $\CO_{C_s}$-module $\CG\subset\CF'_0|_{C_s}$ has degree at most $0$. Hence $\CF'_0|_{C_s}\cong\CO_{C_s}^2$ and thus $H^1(C_s,\CF'_0|_{C_s})=(0)$. By the theorem on cohomology and base change \cite[III.12.11]{Hartshorne} this implies that $\Rkoh^1\pi_\ast\CF'_0$ vanishes and that $\pi_\ast\CF'_0$ is locally free of rank 2 on $T$. Moreover $\pi^\ast\pi_\ast\CF'_0\to\CF'_0$ is an isomorphism in the fiber over $s$ and hence on all of $T$ by Nakayama. Now the level structure $\bar\eta'$ induces an isomorphism $\pi_\ast\CF'_0\cong(\pi^\ast\pi_\ast\CF'_0)|_{I\times T}\isoto\CO_T^2$. From this our claim follows.

As $\tau'_0$ maps $\sigma^\ast\CF'_0$ into $\CF'_1=\CF'_0(\infty)$, it is represented with respect to a basis of $\CF'_0$ by a matrix
\[
\tau'_0 \es = \es (U_0 + tU_1)\cdot\sigma^\ast\qquad\text{with}\quad U_0\in \GL_2(\CO_T)\quad\text{and}\quad U_1\in M_2(\CO_T)\,.
\]
Identifying $\CF'_0|_{I\times T}$ with $\CO_T^2$, we can express $\bar\eta'$ by a matrix in $\GL_2(\CO_T)$. There is a uniquely defined change of basis of $\CF'_0$ such that this matrix becomes the identity. Then we also have $U_0=\Id$.
The condition on $\coker\tau'_0$ now implies that 
\[
\det(\Id + tU_1) \es = \es (1-\zeta' t)^2\,.
\]
Then the required morphism $f:T\to M^{2,2}_I$ is given by $f^\ast(a_{\mu\nu})=U_1$ and $f^\ast(\zeta)=\zeta'$.
\qed
\end{proof}

\begin{remark} \label{RemStacksNotSmooth}
From this example one sees that $\AbSh^{r,d}_H$ need not be smooth over $C$. Namely, $M^{2,2}_I$ is not smooth at the points with $a_{12}=a_{21}=0$, $a_{11}=a_{22}=-\zeta$. The reason for this is that at these points the $\CO_{C_s}$-module $\coker\tau_0$ is isomorphic to $\bigl(\kappa(s)[z]/(z-\zeta)\bigr)^{\oplus 2}$ whereas at all other points it is isomorphic to $\kappa(s)[z]/(z-\zeta)^2$. Compare with Theorem~\ref{ThmSmoothStratification}.

\smallskip

This example also shows that in general one cannot hope that the stacks $\AbSh^{r,d}_H$ are schemes. Namely, for any level $I$ there are abelian sheaves of rank $2$ and dimension $2$ with an $I$-level structure that have non-trivial automorphisms. Indeed, let $I=\Var(a)\subset C'$ for an $a\in\BF_q[t]$ with $\deg a=n$. Let $c:S\to C\setminus I$ be arbitrary and denote the image of $z$ in $\CO_S$ by $\zeta$. Let $f\in\CO_S[t]$ have degree $\leq n+1$. Then the abelian sheaf with
\[
\CF_i\es=\es \CO_{C_S}\bigl((i+n)\cdot\infty\bigr)\oplus\CO_{C_S}(i\cdot\infty)\,,\qquad\tau_i \es=\es\left(
\begin{array}{cc}
1-\zeta t & f \\
0 & 1-\zeta t
\end{array}\right)\cdot\sigma^\ast
\]
admits an $I$-level structure after a finite {\'e}tale extension of $S$. It has non-trivial automorphisms compatible with this level structure of the form
\[
\left(
\begin{array}{cc}
1 & xa \\
0 & 1
\end{array}\right)
\]
for $x\in\BF_q$.

\smallskip

As a consequence $\AbSh^{r,d}_H$ is not quasi-compact in general. Indeed recall from the proof of Theorem~\ref{ThmAlgebraicStacks} the covering of $\AbSh^{r,d}_H$ by open substacks 
\[
U_I\es :=\es \AbSh^{r,d}_{H_I,st}\,/\,(H/H_I)
\]
where $I\subset C'$ runs through all finite subschemes with $H_I\subset H$. Note that $U_I\subset U_{I'}$ if $I\subset I'$. Now if $\AbSh^{r,d}_H$ were quasi-compact this covering would have a finite refinement. So $\AbSh^{r,d}_H$ would equal a single $U_I$ for a large enough $I$. From this we could deduce that 
\[
\AbSh^{r,d}_{H_I}\es = \es\AbSh^{r,d}_{H_I,st}
\]
is in fact a scheme. Since for $r=d=2$ the later is not the case, $\AbSh^{2,2}_H$ is not quasi-compact.
\end{remark}

%
%

\section{Isomorphism Classes Versus Isogeny Classes}  \label{SectIsomVsIsog}
\setcounter{equation}{0}

The general yoga that mediates between isomorphism classes of abelian varieties over $\BZ_p$-schemes and prime-to-$p$ isogeny classes of such can be transfered to abelian sheaves. This will allow us to define $H$-level structures for arbitrary compact open subgroups $H\subset\GL_r(\BA_f)$. We begin by defining the notion of isogeny for abelian sheaves.

\begin{definition} \label{DefQIsogAbSh}
A morphism between abelian sheaves $(\CF_i,\Pi_i,\tau_i)$ and $(\CF'_i,\Pi'_i,\tau'_i)$ over $S$ is called an \emph{isogeny} if 
\begin{enumerate}
\item 
all morphisms $\CF_i\to\CF'_i$ are injective,
\item
$\coker(\CF_i\to\CF'_i)$ is supported on $D\times S$ for an effective divisor $D\subset C$, and
\item
$\coker(\CF_i\to\CF'_i)$ is locally free of finite rank as an $\CO_S$-module.
\end{enumerate}
The isogeny is called \emph{finite} (or \emph{prime to $\infty$}) if for all $i$ the support of $\coker(\CF_i\to\CF'_i)$ is disjoint from $\infty$. A \emph{(finite) quasi-isogeny} between $(\CF_i,\Pi_i,\tau_i)$ and $(\CF'_i,\Pi'_i,\tau'_i)$ is a (finite) isogeny between $(\CF_i,\Pi_i,\tau_i)$ and $(\CF'_i(D),\Pi'_i,\tau'_i)$ for some effective divisor $D\subset C$ (respectively $D\subset C\setminus\infty$).
\end{definition}

Note that an isogeny between two abelian sheaves can only exist if they both have the same rank and dimension. The following proposition is evident. It justifies our definition of quasi-isogenies.

\begin{proposition}
Let $\alpha:\ul{\CF}\to\ul{\CF}'$ be an isogeny of abelian sheaves over $S$. Then there exists an effective divisor $D\subset C$ and an isogeny $\alpha\dual:\ul{\CF}'\to\ul{\CF}(D)$ with $\alpha\circ\alpha\dual$ and $\alpha\dual\circ\alpha$ being the isogenies induced by the inclusion $\CO_C\subset\CO_C(D)$. If $D$ is chosen minimal then $D$ and $\alpha\dual$ are uniquely determined.
\end{proposition}

\begin{example}
Let $\ul\CF=(\CF_i,\Pi_i,\tau_i)$ be an abelian sheaf of rank $r$ and dimension $d$ over $S$. Then the collection of the $\Pi_i$ defines an isogeny
\[
(\Pi_i):\;\ul\CF[1]\to\ul\CF
\]
where $[1]$ denotes the shift by $1$ (cf.\ \ref{DefStack2}). Similarly let $P$ be a point on $C$ and let $S\in\Nilp_{A_P}$ via the characteristic morphism $c:S\to C$. Then the collection of the $\tau_i$ defines an isogeny
\[
(\tau_i):\;\sigma^\ast\ul\CF[1]\to\ul\CF\,.
\]
\end{example}

\medskip

For abelian varieties there is a general principle relating isomorphism classes of abelian varieties with level structure over $\BZ_p$-schemes to prime-to-$p$ isogeny classes of such. This principle also applies to abelian sheaves. We need to introduce the analogue for $\BZ_p$ in our situation. 

\begin{notation} \label{NotationLocalSituation}
The local ring $\CO_{C,\infty}$ is a discrete valuation ring. Let $z$ be a uniformizing parameter. We identify the completion of $\CO_{C,\infty}$ with $\BF_q\dbl z\dbr$ and $Q_\infty$ with $\BF_q\dpl z\dpr$.
Let $\zeta$ be an indeterminant over $\BF_q$ and denote by $\BF_q\dbl\zeta\dbr$ the ring of formal power series in $\zeta$. 
Fix the characteristic morphism $c:\Spec\BF_q\dbl\zeta\dbr\to C$ defined by $c^\ast(z)=\zeta$. Clearly this morphism identifies $\BF_q\dbl z\dbr$ with $\BF_q\dbl\zeta\dbr$. However, since $z-\zeta$ does not necessarily act trivially on $\coker\tau_i$ we use two different symbols to separate the two roles played by $z$ as a uniformizing parameter at $\infty$ and as an element of $\CO_S$. From now on all base schemes for abelian sheaves will be schemes over $\Spec\BF_q\dbl\zeta\dbr$. We consider the base change of our stacks
\[
\AbSh^{r,d}_H\times_C\;\Spec\BF_q\dbl\zeta\dbr\,.
\]
For the sake of brevity we denote them again by $\AbSh^{r,d}_H$. This should not cause confusion since from now on we work entirely in the local situation over $\Spec\BF_q\dbl\zeta\dbr$. 
\end{notation}

Working with isogeny classes instead of isomorphism classes we have to modify our definition of $H$-level structures. The new definition will have the additional advantage that it extends to arbitrary compact open subgroups $H\subset\GL_r(\BA_f)$ which are not necessarily contained in $\GL_r(\wh{A})$.
Let $\ul{\CF}$ be an abelian sheaf of rank $r$ over $S$. We define the functor
\[
\begin{array}{cccc}
(\ul{\CF}|_{\BA_f})^\tau : & \es \Sch_S \es & \longto & \BA_f-\text{modules} \\[2mm]
 & T/S & \longmapsto & \es \BA_f\otimes_\wh{A} \,(\ul{\CF}|_\wh{A})^\tau(T)\,.
\end{array}
\]
Assume that $S$ is connected. Choosing an algebraically closed base point $\iota:s\to S$ we may view $(\ul{\CF}|_{\BA_f})^\tau$ as the $\BA_f[\pi_1(S,s)]$-module $(\iota^\ast\ul{\CF}|_{\BA_f})^\tau(s)$.

We consider the set $\Isom_{\BA_f}\bigl((\ul{\CF}|_{\BA_f})^\tau,\BA_f^{\,r}\bigr):=\Isom_{\BA_f}\bigl((\iota^\ast\ul{\CF}|_{\BA_f})^\tau(s),\BA_f^{\,r}\bigr)$ of isomorphisms of $\BA_f$-modules. Via its natural action on $\BA_f^{\,r}$ the group $\GL_r(\BA_f)$ acts on this set from the left. Via its action on $(\iota^\ast\ul{\CF}|_{\BA_f})^\tau(s)$ the group $\pi_1(S,s)$ acts on it from the right.

\begin{definition} \label{DefRationalLevelStructure}
Let $H\subset \GL_r(\BA_f)$ be a compact open subgroup. A \emph{rational $H$-level structure on $\ul{\CF}$ over $S$} is an $H$-orbit in $\Isom_{\BA_f}\bigl((\ul{\CF}|_{\BA_f})^\tau,\BA_f^{\,r}\bigr)$ which is fixed by $\pi_1(S,s)$. (Again the later condition implies that the notion of level structure is independent of the chosen base point.)
\end{definition}

Every quasi-isogeny $\alpha:\ul{\CF}\to\ul{\CF}'$ induces an isomorphism
\[
(\alpha|_{\BA_f})^\tau: \es(\iota^\ast\ul{\CF}|_{\BA_f})^\tau(s) \es\isoto \es(\iota^\ast\ul{\CF}'|_{\BA_f})^\tau(s)
\]
and thus carries rational $H$-level structures on $\ul{\CF}$ to rational $H$-level structures on $\ul{\CF}'$.

\begin{theorem}
If $H\subset\GL_r(\wh{A})$ is a compact open subgroup then the stack $\AbSh^{r,d}_H$ is canonically 1-isomorphic to the stack $\CX$ whose category of $S$-valued points has 
\begin{description}
\item[as objects:\es] all pairs $(\ul{\CF}, {\bar\gamma})$ consisting of an abelian sheaf $\ul{\CF}$ of rank $r$ and dimension $d$ and a rational $H$-level structure ${\bar\gamma}$ on $\ul{\CF}$ over $S$  and
\item[as morphisms:\es] all finite quasi-isogenies that are compatible with the rational $H$-level structures.
\end{description}
\end{theorem}

\begin{proof}
Let $S$ be a $\Spec\BF_q\dbl\zeta\dbr$-scheme and let $(\ul{\CF},{\bar\eta})$ be an object of $\AbSh^{r,d}_H(S)$. I.e.\ $\ul{\CF}$ is an abelian sheaf and ${\bar\eta}$ is an $H$-level structure on $\ul{\CF}$ over $S$ (in the sense of Definition~\ref{DefHLevelStructure}). Then $\BA_f\otimes_\wh{A}{\bar\eta}$ is a rational $H$-level structure on $\ul{\CF}$. This defines a canonical 1-morphism of stacks $f:\AbSh^{r,d}_H \to \CX$.

For it to be a 1-isomorphism we have to show that $f(S):\AbSh^{r,d}_H(S) \to \CX(S)$ is an equivalence of categories for all $S$. Let $\bigl(\ul{\CF}=(\CF_i,\Pi,i,\tau_i),{\bar\gamma}\bigr)$ be an object of $\CX(S)$. Choose an algebraically closed base point $\iota:s\to S$ and a representative
\[
\gamma:(\iota^\ast\ul{\CF}|_{\BA_f})^\tau(s)\isoto\BA_f^{\,r}\,.
\]
of ${\bar\gamma}$. There is an element $a\in A$ with $\gamma^{-1}(a\wh{A}^r)\subset (\iota^\ast\ul{\CF}|_\wh{A})^\tau(s)$. Then the sheaf 
\[
(\iota^\ast\ul{\CF}|_\wh{A})^\tau(s)\;/\;\gamma^{-1}(a\wh{A}^r)
\]
on $C_s$ has finite length and support disjoint from $\infty$. Via the identification of $(\ul{\CF}|_\wh{A})^\tau$ with $(\iota^\ast\ul{\CF}|_\wh{A})^\tau(s)$ this sheaf can be viewed as a quotient sheaf of $\CF_i$ for all $i$. Its support is of the form $D\times S$ where $D\subset C'$ is the divisor of zeros of $a$. Note that $D\times S$ is disjoint from the characteristic of $\ul{\CF}$. Therefore the kernel of this quotient map is an abelian sheaf $(\CF'_i,\Pi'_i,\tau'_i)$ of rank $r$ and dimension $d$ over $S$ which is isogenous to $\ul{\CF}$. Consider the abelian sheaf $\ul{\CF}'(D):=(\CF'_i(D),\Pi'_i,\tau'_i)$ and the induced finite quasi-isogeny $\alpha$ from $\ul{\CF}'(D)$ to $\ul{\CF}$. By construction $\gamma\circ(\alpha|_{\BA_f})^\tau$ induces an isomorphism
\[
{\bar\eta}: (\ul{\CF}'(D)|_\wh{A})^\tau \isoto \wh{A}^r\,.
\]
This yields an $H$-level structure ${\bar\eta}$ on $\ul{\CF}'(D)$. Clearly this construction is independent of the choice of $\gamma$. 

The pair $(\ul{\CF}'(D),\BA_f\otimes_\wh{A}{\bar\eta})$ is isogenous to $(\ul{\CF},{\bar\gamma})$ by construction. This proves that the functor $f(S)$ is essentially surjective. Analyzing the above construction further shows that it is also fully faithful. Hence $f$ is a 1-isomorphism of stacks.
\qed
\end{proof}

The reader should note that if $v$ is a finite place of $C$, the analogous statement holds for abelian sheaves over $A_v$-schemes and prime-to-$v$ isogenies.

\smallskip

The above theorem enables us to define the stacks of abelian sheaves with $H$-level structure for arbitrary compact open subgroups $H\in\GL_r(\BA_f)$.

\begin{definition} \label{DefStack3}
$\AbSh^{r,d}_H$ is the stack over $\Spec\BF_q\dbl\zeta\dbr$ whose category of $S$-valued points has
\begin{description}
\item[as objects:\es] all pairs $(\ul{\CF}, {\bar\gamma})$ consisting of an abelian sheaf $\ul{\CF}$ of rank $r$ and dimension $d$ and a rational $H$-level structure ${\bar\gamma}$ on $\ul{\CF}$ over $S$  and
\item[as morphisms:\es] all finite quasi-isogenies that are compatible with the rational $H$-level structures.
\end{description}
\end{definition}

For varying $H$ the stacks $\AbSh^{r,d}_H$ form a projective system of stacks.
The transition maps $\AbSh^{r,d}_{H'}\to\AbSh^{r,d}_H$ for $H'\subset H$ are representable by finite {\'e}tale morphisms of schemes. We define a right action of $\GL_r(\BA_f)$ on this projective system by letting $g\in \GL_r(\BA_f)$ act through the 1-isomorphisms
\[
g:\es\AbSh^{r,d}_H\es\isoto\es\AbSh^{r,d}_{g H g^{-1}}
\]
which are defined by $(\ul{\CF},{\bar\gamma})\mapsto(\ul{\CF},\ol{g^{-1}\gamma})$ on $S$-valued points. Using these 1-isomorphisms it follows from Theorem~\ref{ThmAlgebraicStacks} that all the stacks $\AbSh^{r,d}_H$ are algebraic in the sense of Deligne--Mumford
and locally of finite type 
over $\Spec\BF_q\dbl\zeta\dbr$%
.

%

%
%

\bigskip

\section*{{\Large Part Two \es $z$-Divisible Groups}}  \label{PartII}
\setcounter{equation}{0}

\mtaddtocont{\protect\contentsline {mtchap}{\protect\large\protect\rule{0mm}{2.5ex}Part Two \es $z$-Divisible Groups}{\thepage}}

Our ultimate goal in this article is to study the uniformization of the stacks $\AbSh^{r,d}_H$ at $\infty$. Classically, this corresponds to the $p$-adic uniformization of moduli spaces of abelian varieties. For those uniformization questions an indispensable tool are the associated $p$-divisible groups. In the same manner we are thus lead to the idea of ``$z$-divisible groups''. These groups  were studied in detail in Hartl~\cite{Crystals}. But they already appeared in special cases in the work of Drinfeld~\cite{Drinfeld2}, Genestier~\cite{Genestier}, Laumon~\cite{Laumon}, Taguchi~\cite{Taguchi} and Rosen~\cite{Rosen}. For the uniformization of $\AbSh^{r,d}_H$ these $z$-divisible groups are of equal importance than $p$-divisible groups are for abelian varieties. Therefore the next few sections are devoted to them. We first review some facts from \cite{Crystals} in Sections~\ref{SectZDivGps} and \ref{SectCrystals}.

As $z$-divisible groups are of most use over schemes on which $z$ is not a unit we will from now on work over $\Spf\BF_q\dbl\zeta\dbr$ (see Notation~\ref{NotationLocalSituation}). Since we want to relate $z$-divisible groups to abelian sheaves, we should also consider abelian sheaves over $\Spf\BF_q\dbl\zeta\dbr$-schemes. Hence we right away introduce the base change
\[
\AbSh^{r,d}_H\times_C\;\Spf\BF_q\dbl\zeta\dbr\,.
\]
This is no longer an algebraic stack. But it is a \emph{formal algebraic $\Spf\BF_q\dbl\zeta\dbr$-stack} (\ref{DefFormalAlgebraicStacks}). Formal algebraic stacks are related to algebraic stacks in the same way as formal schemes are related to usual schemes. For the necessary background we refer the reader to the appendix.

%
%

\section{Definition of $z$-Divisible Groups}  \label{SectZDivGps}
\setcounter{equation}{0}

We continue to work in the local situation introduced in Notation~\ref{NotationLocalSituation}. In particular $z$ is a uniformizing parameter at $\infty$ and we identify $Q_\infty$ with $\BF_q\dpl z\dpr$. Let $\Delta$ be the central skew field over $Q_\infty$ of invariant $k/\ell$ and let $\CO_{\!\Delta}$ be its ring of integers. We identify $\CO_{\!\Delta}$ with the $\BF_q$-algebra $\BF_{q^\ell}\dbl z,\Pi\dbr$ of non-commutative power series subject to the relations
\[
\Pi^\ell \es = \es z^k\,,\quad z\,\Pi\es=\es \Pi\,z\,,\quad z\,\lambda \es = \es \lambda\,z\,,\quad \Pi\,\lambda^q \es = \es \lambda\,\Pi\quad \text{for all}\quad \lambda\in \BF_{q^\ell}\,.
\]
Denote by $\Nilp_{\BF_q\dbl\zeta\dbr}$ the category of schemes over $\Spec\BF_q\dbl\zeta\dbr$ on which $\zeta$ is locally nilpotent. 
From now on in this article the scheme $S$ will be in $\Nilp_{\BF_q\dbl\zeta\dbr}$.

\begin{definition}
Let $R$ be a unitary ring. We define an \emph{$R$-module scheme over $S$} to be a flat commutative $S$-group scheme $E$ together with a unitary ring homomorphism $R\to\End_S E$. It is called \emph{finite of order $d$} if it is so as an $S$-group scheme. A \emph{morphism of $R$-module schemes} is a morphism of the underlying $S$-group schemes which is compatible with the $R$-action.
\end{definition}

For an $R$-module scheme $E$ over $S$ we define its \emph{coLie module} $\Lie^\ast\! E$ as the $\CO_S$-module of invariant differentials. It is canonically isomorphic to $e^\ast\Omega_{E/S}$ where $e:S\to E$ is the zero section. We have $\Lie E=\CHom_S(\Lie^\ast\! E,\CO_S)$ as $\CO_S$-module.

\medskip

The additive group scheme $\BG_{a,S}$ is an example for an $\BF_q$-module scheme over $S$. Likewise every $S$-group scheme which locally on $S$ is isomorphic to $\BG_{a,S}^d$ for some integer $d\geq0$ is an $\BF_q$-module scheme. Such a scheme is called an \emph{$\BF_q$-vector group scheme of dimension $d$} over $S$. For every $a\in\BF_q$ the endomorphism induced on its coLie module equals the multiplication with $a$ viewed as an element of $\Gamma(S,\CO_S)$.

\begin{definition} \label{DefZDivisibleGroup}
Let $h,d\ge 1$ be integers. A \emph{$z$-divisible group of height $h$ and dimension $d$} over $S$ is an inductive system of finite $\BF_q\dbl z\dbr$-module schemes over $S$
\[
E \es = \es (E_1\xrightarrow{i_1}E_2\xrightarrow{i_2}E_3\xrightarrow{i_3}\ldots)
\]
such that for each integer $n\geq 1$
\begin{enumerate}
\item 
the $\BF_q$-module scheme $E_n$ can be embedded into an $\BF_q$-vector group scheme over $S$,
\item 
the order of $E_n$ is $q^{hn}$,
\item 
the following sequence of $\BF_q\dbl z\dbr$-module schemes over $S$ is exact
\[
0\to E_n\xrightarrow{i_n}E_{n+1}\xrightarrow{z^n}E_{n+1}\,,
\]
\item 
$(z-\zeta)^d =0$ on $\Lie^\ast\! E_n$\,.
\item
$d =\max\{\,\dim_{\kappa(s)}\bigl(\Lie^\ast\! E_n\otimes_{\CO_S}\kappa(s)\bigr):\, s\in S, n\ge1\,\}$\,.
\end{enumerate}
A \emph{morphism of $z$-divisible groups over $S$} is a morphism of inductive systems of $\BF_q\dbl z\dbr$-module schemes. 
\end{definition}

We set $\Lie^\ast\! E=\invlim\Lie^\ast\! E_n$. Conditions 1 to 4 imply that this is a locally free $\CO_S$-module. Condition 5 asserts that its rank is the dimension of $E$. 

The reader should observe that we do not require that $z-\zeta$ acts trivially on $\Lie^\ast\! E$. This is in conformity with the previous sections. In this respect our notion of $z$-divisible group is more general than the variants considered in \cite{Drinfeld2, Genestier, Taguchi, Rosen} and different from the classical case of $p$-divisible groups.

The group of morphisms $\Hom_S(E,E')$ between two $z$-divisible groups $E$ and $E'$ over $S$ is a torsion free $\BF_q\dbl z\dbr$-module. Let $\ul{\Hom}_S(E,E')$ denote the sheaf of germs of morphisms on $S$.

\begin{definition}
The category of \emph{$z$-divisible groups over $S$ up to isogeny} has as objects the $z$-divisible groups over $S$ and as morphisms from $E$ to $E'$ all global sections of the sheaf $\ul{\Hom}_S(E,E')\otimes_{\BF_q\dbl z\dbr}\BF_q\dpl z\dpr$ on $S$. An isomorphism in this new category is called a \emph{quasi-isogeny}. An \emph{isogeny between $z$-divisible groups} is a morphism between $z$-divisible groups which also is a quasi-isogeny.
\end{definition}

In particular, for every quasi-isogeny $\alpha$ there exists locally on $S$ an integer $n$ such that $z^n\alpha$ is an isogeny.
Quasi-isogenies have the following rigidity property.

\begin{proposition} \label{PropQIsogLiftZDiv}
Let $\iota:S'\into S$ be a closed subscheme defined by a sheaf of ideals which is locally nilpotent. Let $E$ and $E'$ be two $z$-divisible groups over $S$. Then every quasi-isogeny from $\iota^\ast E$ to $\iota^\ast E'$ lifts in a unique way to a quasi-isogeny from $E$ to $E'$.
\end{proposition}

\begin{proposition} \label{PropQIsogClosed}
Let $\alpha:E\to E'$ be a quasi-isogeny of $z$-divisible groups over $S$. Then the functor on $\Nilp_{\BF_q\dbl\zeta\dbr}$
\[
T\es\mapsto\es\{\,\phi\in\Hom_{\BF_q\dbl\zeta\dbr}(T,S): \phi^\ast\alpha \text{ is an isogeny}\,\}
\]
is representable by a closed subscheme of $S$.
\end{proposition}

\begin{definition} \label{DefODModule}
A \emph{$z$-divisible $\CO_{\!\Delta}$-module over $S$} is a $z$-divisible group $E$ over $S$ with an action $\CO_{\!\Delta}\to\End_S E$ of $\CO_{\!\Delta}$, which prolongs the natural action of $\BF_q\dbl z\dbr$. A morphism of $z$-divisible $\CO_{\!\Delta}$-modules which is an isogeny of $z$-divisible groups is called an \emph{isogeny}.
\end{definition}

\begin{definition}
If $S$ belongs to $\Nilp_{\Spf\BF_{q^\ell}\dbl\zeta\dbr}$ a $z$-divisible $\CO_{\!\Delta}$-module $E$ of height $r\ell$ and dimension $d\ell$ over $S$ is called \emph{special} if the action of $\CO_{\!\Delta}$ induced on $\Lie^\ast\! E$, makes $\Lie^\ast\! E$ into a locally free $\BF_{q^\ell}\otimes \CO_S$-module of rank $d$.
\end{definition}

At the end of the next section we will show that the later are precisely the $z$-divisible groups that arise from abelian sheaves.

\begin{proposition} \label{PropSpecialClosed}
For a $z$-divisible $\CO_{\!\Delta}$-module $E$ over $S$ the condition of being special is represented by an open and closed immersion into $S$.
\end{proposition}

\begin{proof}
Clearly $\Lie^\ast\! E$ decomposes into a direct sum $\sum_{i=0}^{\ell-1}(\Lie^\ast\! E)_i$ of components $(\Lie^\ast\! E)_i$ on which $\lambda\in\BF_{q^\ell}\subset\CO_{\!\Delta}$ acts via $\lambda^{q^i}\in\CO_S$. Then $E$ is special if and only if all $(\Lie^\ast\! E)_i$ are locally free $\CO_S$-modules of rank $d$. The proposition follows.
\qed
\end{proof}

%
%

\section{Dieudonn\'e modules of $z$-Divisible Groups}  \label{SectCrystals}
\setcounter{equation}{0}

We continue with the notation from Section~\ref{SectZDivGps}. Let $S$ be a scheme in $\Nilp_{\BF_q\dbl\zeta\dbr}$ and consider the completion of $C_S$ along the closed subscheme $\infty\times \Var(\zeta)$. Its structure sheaf is the sheaf $\CO_S\dbl z\dbr$ on $S$ of formal power series in $z$. 
We denote the sheaf $\CO_{\!\Delta}\otimes_{\BF_q\dbl z\dbr}\CO_S\dbl z\dbr$ on $S$ by $\CO_{\!\Delta}\wh{\otimes} \CO_S$.

Consider the additive group $\BG_{a,S}=\Spec\CO_S[\xi]$ over $S$. On $\BG_{a,S}$ we have the \emph{Frobenius isogeny} $Frob_q:\BG_{a,S}\to \BG_{a,S}$ defined by $Frob_q^\ast(\xi)=\xi^q$.
Let $E=(E_n,i_n)$ be a $z$-divisible group over $S$. We associate to $E$ the sheaf
\[
M_E \es=\es\invlim[n]\CH om_S(E_n,\BG_{a,S})
\]
on $S$. We make $M_E$ into a sheaf of $\CO_S\dbl z\dbr$-modules by letting $z$ act through the isogeny $z$ on $E$. The $\sigma_S$-linear multiplication with $\Frob_q$ on the left defines a morphism
\[
F_E:\sigma^\ast M_E\to M_E\,.
\]
If $E$ is moreover a $z$-divisible $\CO_{\!\Delta}$-module then $M_E$ becomes an $\CO_{\!\Delta}\wh{\otimes} \CO_S$-module through the action of $\CO_{\!\Delta}$ from the right. The module $M_E$ may be viewed as the analogue of the \emph{contravariant Dieudonn{\'e} module} associated to a $p$-divisible group. See \cite{Crystals} for a general discussion of this analogy. We recall the following facts.

\begin{definition}
A \emph{Dieudonn\'e $\BF_q\dbl z\dbr$-module over $S$ of dimension $d$ and rank $r$} is a sheaf $\wh\CF$ of $\CO_S\dbl z\dbr$-modules on $S$ equipped with a $\CO_S\dbl z\dbr$-module homomorphism $F:\sigma^\ast \wh\CF\to \wh\CF$ such that locally on $S$ in the Zariski topology 
\begin{enumerate}
\item 
$\wh\CF$ is free of rank $r$ as an $\CO_S\dbl z\dbr$-module,
\item 
$\coker F$ is free of rank $d$ as an $\CO_S$-module,
\item 
$(z-\zeta)^d=0$ on $\coker F$.
\end{enumerate}
A \emph{morphism} between Dieudonn\'e $\BF_q\dbl z\dbr$-modules is a morphism of sheaves of $\CO_S\dbl z\dbr$-modules which is compatible with $F$.
\end{definition}

Note that $F$ is automatically injective. As for $p$-divisible groups the $z$-divisible groups are classified by their Dieudonn\'e $\BF_q\dbl z\dbr$-modules.

\begin{theorem} \label{ThmGroupsVSCrystals}
The functor $E\mapsto (M_E,F_E)$ is an anti-equivalence between the category of $z$-divisible groups of height $r$ and dimension $d$ over $S$ and the category of Dieudonn\'e $\BF_q\dbl z\dbr$-modules of rank $r$ and dimension $d$ over $S$. 
There is a canonical isomorphism $\Lie^\ast\! E=\coker F_E$.
\end{theorem}

If one seeks a classification of $p$-divisible groups up to isogeny one works with isocrystals instead of crystals. For $z$-divisible groups we do the same.

\begin{definition}
A \emph{Dieudonn\'e $\BF_q\dpl z\dpr$-module} over $S$ is a finite locally free $\CO_S\dbl z\dbr[\frac{1}{z}]$-module $\CV$ together with an isomorphism $F:\sigma^\ast\CV\isoto\CV$.
\end{definition}

To every Dieudonn\'e $\BF_q\dbl z\dbr$-module $\ul{\wh\CF}=(\wh\CF,F)$ over $S$ we associate its Dieudonn\'e $\BF_q\dpl z\dpr$-module
\[
\ul{\wh\CF}\,[{\TS\frac{1}{z}}] \es:=\es \bigl(\wh\CF\otimes_{\CO_S\dbl z\dbr}\CO_S\dbl z\dbr[{\TS\frac{1}{z}}]\,,\,F\otimes\id\bigr)\,.
\]

\begin{definition}
A \emph{morphism} $\alpha:\ul{\wh\CF}\to\ul{\wh\CF}'$ of Dieudonn\'e $\BF_q\dbl z\dbr$-modules is called an \emph{isogeny} if $\alpha$ is injective and $\coker\alpha$ is a locally free $\CO_S$-module of finite rank.
We define a \emph{quasi-isogeny} between Dieudonn\'e $\BF_q\dbl z\dbr$-modules $\ul{\wh\CF}$ and $\ul{\wh\CF}'$ to be an isomorphism between $\ul{\wh\CF}\,[{\TS\frac{1}{z}}]$ and $\ul{\wh\CF}'[{\TS\frac{1}{z}}]$.
\end{definition}

\begin{proposition}
The functor $E\mapsto (M_E,F_E)$ maps isogenies to isogenies and quasi-isogenies to quasi-isogenies.
\end{proposition}

Now let $m/n$ be a rational number written in lowest terms with $n>0$. Then we define the Dieudonn\'e $\BF_q\dpl z\dpr$-module $\CV(m/n)$ over $\Spec \BF_q$ as follows
\[
\CV= \BF_q\dpl z\dpr^n\,,\qquad 
F=\left(\begin{array}{cccc}
0\, &  \ldots&&\!z^m \\[-2mm]
\raisebox{1mm}{$1$} & \ddots & \!& \raisebox{-1mm}{$\vdots$}\\[-3.4mm]
& \ddots & \hspace{-1.6mm}\ddots \\
& & 1 \!\!\!& 0
\end{array}\right)\cdot\sigma^\ast\!:\es\sigma^\ast\CV\to\CV\,.
\]
There is the following analogue of Dieudonn{\'e}'s Theorem \cite{Manin}.

\begin{theorem} \label{ThmClassZIsocrystals}
Let $K$ be an algebraically closed field with $\Spec K\in\Nilp_{\BF_q\dbl\zeta\dbr}$. Then every Dieudonn\'e $\BF_q\dpl z\dpr$-module over $\Spec K$ is isomorphic to a direct sum
\[
\bigoplus_i\CV(m_i/n_i)\otimes_{\BF_q\dpl z\dpr} K\dpl z\dpr
\]
for uniquely determined rational numbers $m_1/n_1 \,\leq\, m_2/n_2\,\leq\ldots$\,.
\end{theorem}

This result allows us to define the \emph{Newton polygon} of a Dieudonn\'e $\BF_q\dbl z\dbr$-module $\ul{\wh\CF}=(\wh\CF,F)$ over a field $K$ in $\Nilp_{\BF_q\dbl\zeta\dbr}$. Namely over an algebraically closed extension its Dieudonn\'e $\BF_q\dpl z\dpr$-module decomposes as in the theorem. Then the Newton polygon is the polygon which passes through the points 
\[
(\,n_1+\ldots+n_i\;,\;m_1+\ldots+m_i\,)
\]
for all $i$ and is extended linearly between them. It is independent of the chosen algebraically closed extension. We also define the \emph{Hodge polygon} as usual by the elementary divisors of the $K\dbl z\dbr$-module $\coker F$. The Hodge polygon lies below the Newton polygon. They both have the same initial point $(0,0)$ and the same terminal point $(\rk \ul{\wh\CF},\dim\ul{\wh\CF})$. 
We obtain the analogue of the theorem of Grothendieck--Katz \cite{GrothendieckCristaux, Katz}.

\begin{theorem} \label{ThmNewtonPolygon}
Let $\ul{\wh\CF}$ be a Dieudonn\'e $\BF_q\dbl z\dbr$-module of rank $r$ over $S$ and let $P$ be the graph of a continuous real valued function on $[0,r]$ which is linear between successive integers. Then the set of points in $S$ at which the Hodge (respectively Newton) polygon of $\ul{\wh\CF}$ lies above $P$ is Zariski closed.
\end{theorem}

Note that the stratification of $\AbSh^{r,d}_H$ considered in Section~\ref{SectAlg} is related to the stratification according to the Hodge-polygon. This relation comes from Construction~\ref{ConstrZDivGp} below.

\begin{definition}
Let $S\in\Nilp_{\BF_q\dbl\zeta\dbr}$ be the spectrum of a field. We say that a Dieudonn\'e $\BF_q\dbl z\dbr$-module over $S$ is \emph{isoclinic} if its Newton polygon has only a single slope.
\end{definition}

\begin{proposition} \label{PropIsoclinicCrystal}
Let $\ul{\wh\CF}$ be an isoclinic Dieudonn\'e $\BF_q\dbl z\dbr$-module of rank $r$ and dimension $d$ over a perfect field $K$. Then $\ul{\wh\CF}$ is isogenous over $K$ to a Dieudonn\'e $\BF_q\dbl z\dbr$-module 
satisfying $\im F^r=z^d\wh\CF$.
\end{proposition}

Later we will also need the analogue of Katz's Constancy Theorem.

\begin{theorem} \label{ThmConstancy}
Let $\Spec K\dbl\pi\dbr\in\Nilp_{\BF_q\dbl\zeta\dbr}$ be the spectrum of a power series ring over an algebraically closed field $K$ and let $\ul{\wh\CF}$ be a Dieudonn\'e $\BF_q\dbl z\dbr$-module over $\Spec K\dbl\pi\dbr$. Suppose that at the two points of $\Spec K\dbl\pi\dbr$ the Newton polygons coincide, and that this common Newton polygon has only a single slope. Then $\ul{\wh\CF}$ is isogenous to a constant Dieudonn\'e $\BF_q\dbl z\dbr$-module (i.e.\ one obtained by pullback under the morphism $\Spec K\dbl\pi\dbr\to\Spec K$).
\end{theorem}

\smallskip

We now want to apply this theory to \emph{special} $z$-divisible $\CO_{\!\Delta}$-modules. 
Since for every $z$-divisible group $E$ there is a canonical isomorphism $\Lie^\ast\! E=\coker F_E$ the property of being special reflects on $M_E$. We consider the following class of Dieudonn\'e $\BF_q\dbl z\dbr$-modules. Let $S$ be in $\Nilp_{\BF_{q^\ell}\dbl\zeta\dbr}$.

\begin{definition}
A \emph{formal abelian sheaf of rank $r$ and dimension $d$ over $S$} is a sheaf $\wh{\CF}$ of $\CO_{\!\Delta}\wh{\otimes} \CO_S$-modules on $S$ together with a morphism of $\CO_{\!\Delta}\wh{\otimes} \CO_S$-modules $F:\sigma^\ast\wh\CF\to\wh\CF$ such that
\begin{enumerate}
\item 
$(\wh{\CF},F)$ is a Dieudonn\'e $\BF_q\dbl z\dbr$-module over $S$ of rank $r\ell$,
\item
$\coker F$ is locally free of rank $d$ as an $\BF_{q^\ell}\otimes \CO_S$-module.
\end{enumerate}
\end{definition}

In the situation of special $z$-divisible $\CO_{\!\Delta}$-modules Theorem~\ref{ThmGroupsVSCrystals} takes the following form.

\begin{theorem} \label{ThmAntiEquiv}
The functor $E\mapsto M_E$ is an anti-equivalence between the category of special $z$-divisible $\CO_{\!\Delta}$-modules of height $r\ell$ and dimension $d\ell$ over $S$ and the category of formal abelian sheaves of rank $r$ and dimension $d$ over $S$.
\end{theorem}

\medskip

Next we want to describe the relation with abelian sheaves. We will obtain the analogue of the classical functor which assigns to every abelian variety its $p$-divisible group.

\begin{construction} \label{ConstrZDivGp}
Let $S$ be in $\Nilp_{\BF_q\dbl\zeta\dbr}$ and assume that there is a morphism $\beta:S\to\Spec \BF_{q^\ell}\dbl\zeta\dbr$.
Let $(\CF_i,\Pi_i,\tau_i)$ be an abelian sheaf of rank $r$ and dimension $d$ over $S$.
Consider the completions $\wh{\CF}_i = \CF_i \otimes_{\CO_{C_S}} \CO_S\dbl z\dbr$ of the sheaves $\CF_i$  at $\infty$. Using the periodicity $\CF_\ell\cong\CF_0(k\cdot\infty)$ we obtain morphisms
\begin{eqnarray}
\Pi_i: \wh{\CF}_i  \to \wh{\CF}_{i+1}  \quad \text{for all} \es i=0,\ldots,\ell-2\,, \quad  & z^k\,\Pi_{\ell-1}:\;&\wh{\CF}_{\ell-1}  \to  \wh{\CF}_0  \quad \text{and}  \hfill\label{EqFuerPi} \\[0.2cm]
{}\qquad\tau_i:\sigma^\ast\wh{\CF}_i  \to \wh{\CF}_{i+1}  \quad \text{for all} \es i=0,\ldots,\ell-2\,, \quad & \;z^k\,\tau_{\ell-1}:&\sigma^\ast\wh{\CF}_{\ell-1}  \to  \wh{\CF}_0 \label{EqFuerTau}
\end{eqnarray}
We set $\wh{\CF}=\wh{\CF}_0 \oplus \ldots \oplus \wh{\CF}_{\ell-1}$ and let the endomorphism $\Pi:\wh{\CF}\to \wh{\CF}$ be given by the morphisms (\ref{EqFuerPi}). We let $\lambda\in \BF_{q^\ell}$ act on $\wh{\CF}_i$ as the scalar $\beta^\ast\lambda^{q^i}$ and we let the $\sigma$-linear endomorphism $F:\sigma^\ast\wh{\CF}\to \wh{\CF}$ be given by the morphisms (\ref{EqFuerTau}). I.e.\ $\Pi$, $\lambda$ and $F$ are expressed by the block matrices
\begin{eqnarray*}
& \Pi = \left(
\begin{array}{cccc}
0\, &  \ldots&&\!\! z^k\Pi_{\ell-1} \\
\raisebox{1mm}{$\Pi_0$} & \ddots & & \raisebox{-1mm}{$\vdots$}\\[-3.4mm]
& \ddots & \hspace{-1.6mm}\ddots \\
& & \Pi_{\ell-2} \!\!& 0
\end{array}\right),\quad
\lambda = \left(
\begin{array}{cccc}
\beta^\ast\lambda \Id_\ell\\
&\beta^\ast\lambda^q\Id_\ell  \\[-2mm]
& & \ddots \!\!\\
& & & \beta^\ast\lambda^{q^{\ell-1}}\Id_\ell
\end{array}\right)\es\text{and}\\[2mm]
& F = \left(
\begin{array}{cccc}
0\, &  \ldots&&\!\! z^k\tau_{\ell-1} \\
\raisebox{1mm}{$\tau_0$} & \ddots & & \raisebox{-1mm}{$\vdots$} \\[-3.5mm]
& \ddots & \hspace{-1.1mm}\ddots \\
& & \tau_{\ell-2} \!\!& 0
\end{array}\right)\,.
\end{eqnarray*}
In this way $(\wh\CF,F)$ is a formal abelian sheaf of rank $r$ and dimension $d$ over $S$.
By Theorem~\ref{ThmAntiEquiv} we obtain a functor from $\AbSh^{r,d}(S)$ to the category of special $z$-divisible $\CO_{\!\Delta}$-modules of height $r\ell$ and dimension $d\ell$ over $S$.

Moreover, every isogeny of abelian sheaves over $S$ induces an isogeny of the associated formal abelian sheaves and an isogeny of the associated special $z$-divisible $\CO_{\!\Delta}$-modules.
\end{construction}

%
%

\section{The Serre-Tate-Theorem} \label{SectSerreTate}
\setcounter{equation}{0}

Classically the Serre-Tate-Theorem relates the deformation theory of an abelian variety in characteristic $p$ to the deformation theory of its $p$-divisible group. In the case of abelian sheaves the same principle prevails. We begin with the analogue of the following classical construction. Let $A$ be a fixed abelian variety over a scheme $S\in\Nilp_{\BZ_p}$ and let $A[p^\infty]$ be its $p$-divisible group. Then to every pair $(X,\wh{\alpha})$ consisting of a $p$-divisible group $X$ over $S$ and an isogeny $\wh{\alpha}:A[p^\infty]\to X$ of $p$-divisible groups there exists a uniquely determined abelian variety $\wh{\alpha}_\ast A$ and a $p$-isogeny $\alpha:A\to \wh{\alpha}_\ast A$ which induces $\wh{\alpha}$ on $p$-divisible groups. This construction can be applied to abelian sheaves as well.

\begin{proposition} \label{PropPullbackQIsog}
Let $S\in\Nilp_{\BF_{q^\ell}\dbl\zeta\dbr}$ and let $\ul{\CF}$ be an abelian sheaf of rank $r$ and dimension $d$ over $S$. Let $\ul{\wh\CF}$ be the associated formal abelian sheaf and let $\wh{\alpha}:\ul{\wh\CF}'\to\ul{\wh\CF}$ be a quasi-isogeny. Then there exists an abelian sheaf $\ul{\CF}'$ and a quasi-isogeny $\alpha:\ul{\CF}'\to\ul{\CF}$ over $S$ giving rise to $\wh{\alpha}$. If we require that $\alpha$ is an isomorphism over $C'$ then $(\ul{\CF}',\alpha)$ is unique up to canonical isomorphism. In this case we denote $\ul{\CF}'$ by $\wh{\alpha}^\ast\ul{\CF}$.
\end{proposition}

\begin{proof}
We denote by $\beta:S\to\Spf\BF_{q^\ell}\dbl\zeta\dbr$ the structure morphism of $S$. For each $i=0,\ldots,\ell-1$ we extract from the sheaf ${\wh{\CF}}'$ of $\CO_S\dbl z\dbr$-modules underlying $\ul{\wh\CF}'$ the locally free subsheaf $\wh{\CF}'_i\subset\wh{\CF}'$ of rank $r$ on which $\lambda\in\BF_{q^\ell}$ acts through the character $\lambda\mapsto\beta^\ast\lambda^{q^i}$.
The quasi-isogeny $\wh\alpha$ gives an inclusion $\wh\CF'_i\into\CF_i\otimes_{\CO_{C'_S}}\CO_S\dbl z\dbr[\frac{1}{z}]$. This permits us to glue $\wh{\CF}'_i$ with the corresponding sheaf $\CF_i|_{C'_S}$ 
to obtain a locally free sheaf $\CF'_i$ of rank $r$ on $C_S$. For arbitrary $i\in \BZ$ we take $n\in \BZ$ such that $0\leq i-n\ell <\ell$ and define
\[
\CF'_i \;:=\; \CF'_{i-n\ell} \otimes \CO_{C_S}(nk\cdot\infty)\,.
\]
From the morphisms $\Pi'$ and $F'$ of $\ul{\wh\CF}'$ and $\Pi_i$ and $\tau_i$ of $\ul\CF$ we obtain morphisms
\[
\Pi'_i:\CF'_i\to\CF'_{i+1}\quad \text{and}\quad\tau'_i:\sigma^\ast\CF'_i\to\CF'_{i+1}\,.
\]
The morphism $\Pi'_{\ell-1}:\CF'_{\ell-1}\to\CF'_\ell$ is induced from $z^{-k}\Pi':\wh{\CF}'_{\ell-1}\to\wh{\CF}'_{\ell}(k\cdot\infty)$. The same applies to $\tau'_{\ell-1}$. These morphisms make $\wh\alpha^\ast\ul\CF:=(\CF'_i,\Pi'_i,\tau'_i)$ into an abelian sheaf of rank $r$ and dimension $d$ over $S$ whose associated formal abelian sheaf is $\ul{\wh\CF}'$. 
By construction there is a quasi-isogeny $\alpha:\wh\alpha^\ast\ul\CF\to\ul\CF$ which is an isomorphism over $C'$ and induces the quasi-isogeny $\wh\alpha$ on formal abelian sheaves.
\qed
\end{proof}

\begin{proposition} \label{PropQIsogLiftAbSh}
Let $S\in\Nilp_{\BF_{q^\ell}\dbl\zeta\dbr}$ and let $j:\bar{S}\into S$ be a closed subscheme defined by a sheaf of ideals which is locally nilpotent. Let $\ul{\CF}$ and $\ul{\CF}'$ be two abelian sheaves of rank $r$ and dimension $d$ over $S$. Then every quasi-isogeny from $j^\ast\ul{\CF}$ to $j^\ast\ul{\CF}'$ lifts in a unique way to a quasi-isogeny from $\ul{\CF}$ to $\ul{\CF}'$.
\end{proposition}

\begin{proof}
Let $\bar{\alpha}:j^\ast\ul{\CF}\to j^\ast\ul{\CF}'$ be a quasi-isogeny.
It suffices to treat the case where the $q$-th power of the ideal sheaf defining $\bar{S}$ is zero. In this case the morphisms $\sigma_S$ and $\sigma_{\bar S}$ factor through $j$
\[
\sigma_S\es = \es j\circ\bar{\sigma}: \es S\to \bar{S}\to S\qquad\text{and}\qquad
\sigma_{\bar S}\es = \es \bar{\sigma}\circ j: \es \bar{S}\to S\to \bar{S}\,.
\]
Consider the quasi-isogeny $\bar{\sigma}^\ast\bar{\alpha}[1]:\bar{\sigma}^\ast j^\ast\ul{\CF}[1]\to\bar{\sigma}^\ast j^\ast\ul{\CF}'[1]$ where $[1]$ denotes the shift by $1$ (cf.\ \ref{DefStack2}). We view the morphisms $\tau_i$ as an isogeny $(\tau_i):\sigma^\ast\ul{\CF}[1]\to\ul{\CF}$ and obtain a commutative diagram
\[
\begin{CD}
\bar{\sigma}^\ast j^\ast\ul{\CF}[1] \es& @>\;\bar{\sigma}^\ast\bar{\alpha}[1]\;>> & \es\bar{\sigma}^\ast j^\ast\ul{\CF}'[1]\\
@V(\tau_i)VV & & @VV(\tau'_i)V \\
\ul{\CF} & @>\alpha>> & \ul{\CF}'
\end{CD}
\]
which defines a quasi-isogeny $\alpha$. Pulling back this diagram under $j$ we see that \mbox{$j^\ast\alpha=\bar{\alpha}$}. 
Moreover the diagram shows that $\alpha$ is uniquely determined by $\bar\alpha$. This proves the proposition.
\qed
\end{proof}

From the proof we even see the following.

\begin{corollary} \label{CorIsomOverC'}
Keep the situation of the proposition. 
\begin{enumerate}
\item 
If $\bar\alpha:j^\ast\ul\CF\to j^\ast\ul\CF'$ is an isogeny then the lift is an isogeny $\alpha:\ul\CF\to\ul\CF'(n\cdot\infty)$ for some integer $n\ge 0$.
\item 
If $\bar\alpha$ is an isomorphism over $C'$ then the same holds for $\alpha$.
\end{enumerate}
\end{corollary}

Next we come to the analogue of the Serre-Tate-Theorem. Let $S\in\Nilp_{\BF_{q^\ell}\dbl\zeta\dbr}$ and let $j:S'\into S$ be a closed subscheme defined by a sheaf of ideals which is locally nilpotent. Let $\ul{\CF}'$ be an abelian sheaf of rank $r$ and dimension $d$ over $S'$ and let $\ul{\wh\CF}'$ be the associated formal abelian sheaf. The \emph{category of lifts of $\ul\CF'$ to $S$} has 
\begin{description}
\item[as objects:\es] all pairs $(\ul\CF\,,\,\alpha:j^\ast\ul\CF\isoto\ul\CF')$ where $\ul\CF$ is an abelian sheaf over $S$ and $\alpha$ an isomorphism of abelian sheaves over $S'$,
\item[as morphisms:\es] isomorphisms between the $\ul\CF$'s that are compatible with the $\alpha$'s.
\end{description}
Similarly we define the \emph{category of lifts of $\ul{\wh\CF}'$ to $S$}. By Propositions~\ref{PropQIsogLiftAbSh} and \ref{PropQIsogLiftZDiv} all $\Hom$-sets in these categories contain at most one element.

\begin{theorem}[Analogue of the Serre-Tate-Theorem] \label{ThmSerreTate}
The category of lifts of $\ul{\CF}'$ to $S$ and  the category of lifts of $\ul{\wh\CF}'$ to $S$ are equivalent.
\end{theorem}

\begin{proof}
Let $\wh{\es}$ be the functor that assigns to a lift of $\ul{\CF}'$ the corresponding lift of $\ul{\wh\CF}'$. Full faithfulness of $\wh{\es}$ follows from Corollary~\ref{CorIsomOverC'}. It remains to show that $\wh{\es}$ is essentially surjective. So let $(\ul{\wh\CF}\,, \;\wh{\alpha}:j^\ast\ul{\wh\CF}\isoto\ul{\wh\CF}')$ be a lift of $\ul{\wh\CF}'$ to $S$.

It suffices to treat the case where the $q$-th power of the ideal sheaf defining $S'$ is zero. In this case the morphism $\sigma_S$ factors through $j$
\[
\sigma_S\es= \es j\circ\sigma':\es S\to S'\to S \,.
\]
Consider the abelian sheaf $\ul{\wt{\CF}}:=(\sigma'{}^\ast\ul{\CF}')[1]$ over $S$. I.e.
\[
(\wt{\CF}_i,\wt{\Pi}_i,\wt{\tau}_i) \es=\es (\sigma'{}^\ast\CF'_{i-1},\sigma'{}^\ast\Pi'_{i-1},\sigma'{}^\ast\tau'_{i-1})\,.
\]
The morphisms $\tau'_i$ constitute an isogeny $\tau':=(\tau'_i):j^\ast\ul{\wt{\CF}} \to\ul{\CF}'$ which is an isomorphism over $C'$. We let 
\[
\wh{\gamma}'\es = \es \wh{\tau}'{}^{-1} \circ \wh{\alpha}:\es j^\ast\ul{\wh\CF} \to \ul{\wh\CF}'\to j^\ast\ul{\wh{\wt{\CF}}}
\]
be the resulting quasi-isogeny of formal abelian sheaves. By Proposition~\ref{PropQIsogLiftZDiv} it lifts to a quasi-isogeny $\wh{\gamma}:\ul{\wh\CF} \to\ul{\wh{\wt{\CF}}}$. We put $\ul{\CF}:=\wh{\gamma}^\ast\ul{\wt{\CF}}$ (Prop.~\ref{PropPullbackQIsog}) and we let $\gamma:\ul{\CF}\to\ul{\wt{\CF}}$ be the induced quasi-isogeny of abelian sheaves. Then $(\ul{\CF},\tau'\circ j^\ast\gamma)$ is the desired lift of $\ul{\CF}$.
\qed
\end{proof}

In the remainder of this section we give an example for an abelian sheaf and we compute the associated formal abelian sheaf. This example will be crucial for the uniformization of $\AbSh^{r,d}_H$ in Part Three.

\begin{bigexample} \label{ExAbSheaf}
On the scheme $\bar S=\Spec\BF_q$ we set for $i=0,\ldots,\ell$
\[
\CM_i\es=\es\CO_{C_{\bar S}}(k\cdot\infty)^{\oplus i}\oplus \CO_{C_{\bar S}}^{\oplus \ell-i}\,.
\]
We let $\Pi_i:\CM_i\to\CM_{i+1}$ be the morphism coming from the natural inclusion $\CO_{C_{\bar S}}\subset\CO_{C_{\bar S}}(k\cdot\infty)$ in the $(i+1)$-st summand.
We define the $\sigma$-linear morphism
\[\tau_i= \left(
\begin{array}{cccc}
0\, &  \ldots&&\!\! 1 \\[-2mm]
\raisebox{1mm}{$1$} & \ddots & & \raisebox{-1mm}{$\vdots$}\\[-3.4mm]
& \ddots & \hspace{-1.6mm}\ddots \\[-1mm]
& & \!\!1 & 0
\end{array}\right)\!\cdot\sigma^\ast:\es\sigma^\ast\CM_i\to\CM_{i+1}\,.
\]
For arbitrary $i\in \BZ$ we take $n\in \BZ$ such that $0\leq i-n\ell <\ell$ and set
\[
\CM_i := \CM_{i-n\ell} \otimes \CO_{C_{\bar S}}(nk\cdot\infty)\,,\quad \Pi_i := \Pi_{i-n\ell} \otimes \id\,,\quad \tau_i :=\tau_{i-n\ell} \otimes \id\,.
\]
Then $\ul{\CM}:=(\CM_i,\Pi_i,\tau_i)$ is an abelian sheaf of rank $\ell$ and dimension $k$ over $\Spec\BF_q$. We let $e=\frac{r}{\ell}$ and set $\ol{\BM}=\ul{\CM}^{\oplus e}$. This is an abelian sheaf of rank $r$ and dimension $d$. We compute the formal abelian sheaf $\wh{\ol{\BM}}$ associated to $\ol{\BM}$.

In order to do this we have to extend the base scheme $\bar S$ to $\bar S'=\Spec\BF_{q^\ell}$. Since we afterwards want to lift $\wh{\ol{\BM}}$ to $S'=\Spf\BF_{q^\ell}\dbl\zeta\dbr$, we right away describe this lift. For $i=0,\ldots,\ell-1$ we consider the $\ell\times\ell$-matrices
\begin{equation} \label{EqPi_i}
\Pi_i=\left(
\begin{array}{ccccc}
1 &\\[-2mm]
& \ddots \\[-2mm]
& &z^k  \\[-2mm]
& & & \ddots \\[-2mm]
& & & & 1
\end{array}\right)\quad \text{and}\quad 
T= \left(
\begin{array}{cccc}
0\, &  \ldots&&\!\! (z-\zeta)^k \\[-2mm]
\raisebox{1mm}{$1$} & \ddots & & \raisebox{-1mm}{$\vdots$}\\[-3.4mm]
& \ddots & \hspace{-1.6mm}\ddots \\
& & 1 \!\!& 0
\end{array}\right)\,,
\end{equation}
where the $z^k$ in $\Pi_i$ sits in the $(i+1)$-st row. We let the formal abelian sheaf $\wh{\CM}$ of rank $\ell$ and dimension $k$ be the $\CO_{S'}\dbl z\dbr$-module
$\wh{\CM} = \CO_{S'}\dbl z\dbr^{\ell^2}$ together with the morphisms
\begin{eqnarray*}
& \Pi = \left(
\begin{array}{cccc}
0\, &  \ldots&&\!\!\! \Pi_{\ell-1} \\[-1mm]
\raisebox{1mm}{$\Pi_0$} & \ddots & & \!\raisebox{-1mm}{$\vdots$}\\[-3.4mm]
& \ddots & \hspace{-4mm}\ddots \\
& & \Pi_{\ell-2} \!\!\!& 0
\end{array}\right),\quad
\lambda = \left(
\begin{array}{cccc}
\lambda \Id_r\\
&\lambda^q\Id_r  \\[-2mm]
& & \ddots \!\!\\
& & & \lambda^{q^{\ell-1}}\Id_r
\end{array}\right):\;\wh{\CM}\to\wh{\CM}\\[2mm]
& \text{and}\qquad F = \left(
\begin{array}{cccc}
0\, &  \ldots&& T \\[-2mm]
\raisebox{1mm}{$T$} & \ddots & & \raisebox{-1mm}{$\vdots$} \\[-3.5mm]
& \ddots & \hspace{-1.1mm}\ddots \\
& & \!T & 0
\end{array}\right)\!\cdot\sigma^\ast:\es\sigma^\ast\wh{\CM}\to\wh{\CM}\,.
\end{eqnarray*}
We set $\wh{\BM}=\wh{\CM}^{\oplus e}$. If $j:\bar S'\into S'$ denotes the inclusion, then $j^\ast\wh{\BM}$ is the formal abelian sheaf associated to $\ol{\BM}$. We denote it by $\wh{\ol{\BM}}$. Via Theorem~\ref{ThmAntiEquiv} we obtain from $\wh{\ol{\BM}}$ a special $z$-divisible $\CO_{\!\Delta}$-module $\ol{\BE}$ over $\Spec\BF_{q^\ell}$ whose Dieudonn\'e $\BF_q\dbl z\dbr$-module is $\wh{\ol{\BM}}$. 
One easily checks $F^\ell=z^k\cdot\sigma^{\ell\ast}$ on $\wh{\ol{\BM}}$. In the terminology of \cite{Crystals} this means that $\ol\BE$ is \emph{descent}. The Newton-polygon of $\wh{\ol{\BM}}$ is a straight line between the end points $(0,0)$ and $(r\ell,d\ell)$.  

Via the Serre-Tate-Theorem~\ref{ThmSerreTate} we obtain from the lift $\wh\BM$ of the formal abelian sheaf of $\ol\BM$ an abelian sheaf $\BM$ over $\Spf\BF_{q^\ell}$ which lifts $\ol\BM$.
\end{bigexample}

Let us compute the group of quasi-isogenies of these (formal) abelian sheaves.

\begin{proposition} \label{PropQIsogMBar}
Fix a generator $\lambda\in \BF_{q^\ell}$ of the extension $\BF_{q^\ell}/\BF_q$. Let $V$ be the Vandermonde matrix
\begin{equation} \label{EqVandermonde}
V\es=\es\left(
\begin{array}{ccccc}
1 &  \lambda\es &\lambda^2\es &\ldots &\lambda^{\ell-1}\es  \\
1 &  \;\lambda^q\es &\;\lambda^{2q}\es &\ldots &\;\lambda^{(\ell-1)q}\es  \\
\vdots &\vdots &\vdots & &\vdots \\
1 &  \es \lambda^{q^{\ell-1}}&\es \lambda^{2q^{\ell-1}}&\ldots &\es \lambda^{(\ell-1)q^{\ell-1}} 
\end{array}\right)\quad\text{and}\quad V_e\es=\es\left(
\begin{array}{ccc}
V \\[-3mm]
& \!\!\ddots\! & \\[-1mm]
& & V
\end{array}
\right)
\end{equation}
be the block diagonal matrix of dimension $r$.
Then the group of quasi-isogenies of $\ol{\BM}$ over $\BF_q^{\,\alg}$ is the group of $Q$-valued points of the algebraic group  $J=V_e\GL_r V_e^{-1}$ over $Q$.
\end{proposition}

\begin{proof}
Let $g\in J(Q)$ and let $D\subset C$ be an effective divisor  satisfying $\di(g_{\mu\nu})\geq -D$ for all entries of $g$. Then multiplication with $g$ defines morphisms $\iota^\ast\CM_i^{\oplus e}\to\iota^\ast\CM_i^{\oplus e}(D)$ which form a quasi-isogeny of $\ol{\BM}$ (defined over $\Spec\BF_{q^\ell}$). Conversely one computes that every quasi-isogeny of $\ol{\BM}$ arises in this way.
\qed
\end{proof}

If $T\in\Nilp_{\BF_{q^\ell}\dbl\zeta\dbr}$ is a scheme we let $\beta:T\to\Spf\BF_{q^\ell}\dbl\zeta\dbr$ be its structure morphism and $\bar\beta:\bar T\to\Spec\BF_{q^\ell}$ be the reduction modulo $\zeta$.

\begin{proposition} \label{PropQIsogMS}
Let $T\in\Nilp_{\BF_{q^\ell}\dbl\zeta\dbr}$ be connected. Then there is a canonical isomorphism $g\mapsto g_T$ of $J(Q)$ to the group ${\rm QIsog}_T(\beta^\ast\BM)$ of quasi-isogenies of $\beta^\ast\BM$ over $T$.
\end{proposition}

\begin{proof}
The map $g\mapsto g_T$ is defined as
\[
\begin{array}{ccccc}
J(Q)&\es\longto\es& {\rm QIsog}_{\bar T}(\bar\beta^\ast\ol{\BM}) &\es\isoto\es &{\rm QIsog}_T(\beta^\ast\BM)\\[1mm]
g & \mapsto & \bar\beta^\ast g & \mapsto & g_T\,,
\end{array}
\]
the last isomorphism coming from Proposition~\ref{PropQIsogLiftAbSh}.
Clearly this defines a monomorphism of groups.

We show that $g\mapsto \bar\beta^\ast g$ is surjective. Let $\alpha:\ol{\BM}_{\bar T}\to\ol{\BM}_{\bar T}(D)$ be an isogeny for some effective divisor $D\subset C$. There exists an $a\in A$ whose divisor is $\geq D$ on $C'$. We view $\alpha|_{C'_{\bar T}}$ as a matrix $U\in M_r(\CO_{\bar T}\otimes A[\frac{1}{a}])$. Then the matrix $V_e^{-1}UV_e$ satisfies
\[
{}^{\sigma\!}(V_e^{-1}UV_e)\es=\es V_e^{-1}UV_e
\]
and hence lies in $\GL_r(Q)$. We conclude that $U\in J(Q)$ and $\bar\beta^\ast U=\alpha$.
\qed
\end{proof}

\begin{proposition}\label{PropQIsogMHat}
The group of quasi-isogenies of the formal abelian sheaf $\wh\BM$ is isomorphic to $J(Q_\infty)$. 
\end{proposition}

\begin{proof}
A straightforward calculation shows that this isomorphism can be described as follows. Let
\[
W_i\es=\es \left(
\begin{array}{cc}
z^k\Id_i & 0 \\
0 & \Id_{\ell-i}
\end{array}\right)
\qquad \text{and} \qquad
W\es=\es\left(
\begin{array}{cccc}
g \\
& W_1 gW_1^{-1} \\[-2mm]
& & \!\!\ddots\!\! & \\
& & & W_{\ell-1}gW_{\ell-1}^{-1}
\end{array}
\right)
\,.
\]
Then an element $\in (Q_\infty)$ is mapped to the quasi-isogeny $W^{\oplus e}$ of $\wh{\CM}$.
\qed 
\end{proof}

%
%

\section{Moduli Spaces for $z$-Divisible Groups}  \label{SectModSpZDiv}
\setcounter{equation}{0}

Consider the formal abelian sheaf $\wh{\BM}$ from Example~\ref{ExAbSheaf}. We will define a moduli problem for formal abelian sheaves which are quasi-isogenous to $\wh{\BM}$. This is a higher dimensional variant of a moduli problem studied by Drinfeld~\cite{Drinfeld2}. At the same time it is a close analogue of a moduli problem for $p$-divisible groups considered by Rapoport--Zink~\cite{RZ}. Like these two problems our moduli problem too will be solved by a formal scheme over $\Spf\BF_{q^\ell}\dbl\zeta\dbr$. Following \cite{Drinfeld2, RZ} we will use this formal moduli scheme in Section~\ref{SectUnifOfAbSh} to (partly) uniformize the stacks of abelian sheaves.

For a scheme $S$ in $\Nilp_{\BF_q\dbl\zeta\dbr}$ we denote by $\bar S$ the closed subscheme defined by the sheaf of ideals $\zeta\CO_S$. We call $\bar S$ the special fiber of $S$. If $E$ is a $z$-divisible group over $S$ we denote by $E_{\bar S}=E\times_S{\bar S}$ the base change to the special fiber. A similar notation will be applied to Dieudonn\'e $\BF_q\dbl z\dbr$-modules, etc.
If $\beta:S\to\Spf\BF_{q^\ell}\dbl\zeta\dbr$ is a morphism of formal schemes we denote by $\bar\beta:{\bar S}\to\Spec\BF_{q^\ell}$ its restriction to special fibers. 
We define the following moduli problem for formal abelian sheaves. 

\begin{definition}
Let $G$ be the contravariant functor $\Nilp_{\BF_q\dbl\zeta\dbr} \longto \CS ets$
\begin{eqnarray*}
S &\longmapsto \Bigl\{ &\text{Isomorphism classes of triples}\es(\beta,\ul{\wh\CF},\wh\alpha) \es\text{where} \\
&&\es \bullet\es \beta:S\to\Spf\BF_{q^\ell}\dbl\zeta\dbr \text{ is a morphism of formal schemes,}\\
&&\es \bullet\es \ul{\wh\CF} \text{ is a formal abelian sheaf of rank $r$ and dimension $d$ over $S$}\,,\\
&&\es \bullet\es \wh\alpha:\ul{\wh{\CF}}_{\bar S}\to \bar\beta^\ast \wh{\ol\BM}\text{ is a quasi-isogeny of formal abelian sheaves.}\es\Bigr\}
\end{eqnarray*}
\end{definition}

Thereby two  triples $(\beta,\ul{\wh\CF},\wh\alpha)$ and $(\beta',\ul{\wh\CF}',\wh\alpha')$ are isomorphic if $\beta=\beta'$ and if there is an isomorphism between $\ul{\wh\CF}$ and $\ul{\wh\CF}'$ over $S$ which is compatible with $\wh\alpha$ and $\wh\alpha'$.

Using the equivalence between special $z$-divisible $\CO_{\!\Delta}$-modules and formal abelian sheaves from Theorem~\ref{ThmAntiEquiv} we see that 

\begin{proposition}
$G$ can be described as the functor $\Nilp_{\BF_q\dbl\zeta\dbr} \longto \CS ets$
\begin{eqnarray*}
&S &\longmapsto \Bigl\{\text{ Isomorphism classes of triples}\es(\beta,E,\rho) \es\text{where} \\[2mm]
&&\es \bullet\es \beta:S\to\Spf\BF_{q^\ell}\dbl\zeta\dbr \text{ is a morphism of formal schemes,}\\
&&\es \bullet\es E \text{ is a special $z$-divisible $\CO_{\!\Delta}$-module of height $r\ell$ and dimension $d\ell$ over $S$}\,,\\
&&\es \bullet\es \rho:\bar\beta^\ast\ol{\BE}\to E_{\bar S} \text{ is a quasi-isogeny of $z$-divisible $\CO_{\!\Delta}$-modules.}\es\Bigr\}
\end{eqnarray*}
\end{proposition}

Sending $(\beta,\ul{\wh\CF},\wh\alpha)$ to $\beta$ gives a morphism $G\to \Spf\BF_{q^\ell}\dbl\zeta\dbr$.
We define an action of the Galois group $\Gal(\BF_{q^\ell}/\BF_q)$ on $G$ over $\Spf\BF_{q^\ell}\dbl\zeta\dbr$. Namely for each $\pi^\ast\in \Gal(\BF_{q^\ell}/\BF_q)=\Gal\bigl(\BF_{q^\ell}\dbl\zeta\dbr/\BF_q\dbl\zeta\dbr\bigr)$ consider the cartesian square
\[
\begin{CD}
S^\pi & @>\pi_S>> & S \\
@VV\pi^\ast\beta V & & @VV\beta V \\
\Spf\BF_{q^\ell}\dbl\zeta\dbr & @>\pi>> &\Spf\BF_{q^\ell}\dbl\zeta\dbr\es.
\end{CD}
\]
Then we let $\pi^\ast$ act by mapping the element $(\beta,\ul{\wh\CF},\wh\alpha)\in G(S)$ to
\[
\bigl(\pi^\ast\beta,\pi_S^\ast\ul{\wh\CF},\pi_S^\ast\wh\alpha\bigr)\es\in\es G(S^\pi)\,.
\]

\medskip

Following the arguments given in Rapoport--Zink~\cite{RZ} in the case of $p$-divisible groups, one can prove the following representability theorem. Recall that an adic formal scheme $G$ over $\Spf\BF_{q^\ell}\dbl\zeta\dbr$ is called \emph{locally formally of finite type} if $G_\red$ is locally of finite type over $\Spec\BF_{q^\ell}$.

\begin{theorem}\label{ThmGTildeRepresentable}
The functor that assigns to a scheme $S\in\Nilp_{\BF_q\dbl\zeta\dbr}$ the set of isomorphism classes of triples $(\beta,E,\rho)$ where
\begin{eqnarray*}
&&\bullet\es \beta:S\to\Spf\BF_{q^\ell}\dbl\zeta\dbr \text{ is a morphism of formal schemes,}\\
&&\bullet\es E\text{ is a $z$-divisible group of dimension $d\ell$ and height $r\ell$ over $S$,\qquad\qquad}\\ 
&&\bullet\es \rho:\bar\beta^\ast\ol{\BE}\to E_{\bar S} \text{ is a quasi-isogeny of $z$-divisible groups,}
\end{eqnarray*}
is representable by a quasi-separated, locally noetherian, adic formal scheme which is locally formally of finite type over $\Spf\BF_{q^\ell}\dbl\zeta\dbr$.
\end{theorem}

The proof of this theorem is given in \cite{Crystals}. It makes use of Dieudonn\'e $\BF_q\dbl z\dbr$-modules which replace the crystals of the $p$-divisible groups in \cite{RZ}.

\begin{corollary} \label{ThmGRepresentable}
The functor $G$ is representable by a quasi-separated, locally noetherian, adic formal scheme which is locally formally of finite type over $\Spf\BF_{q^\ell}\dbl\zeta\dbr$.
\end{corollary}

\begin{proof}
Let $\wt G$ be the formal scheme whose existence is stated in Theorem~\ref{ThmGTildeRepresentable}. Let $E$ be the universal $z$-divisible group over $\wt G$. We transport the $\CO_{\!\Delta}$-action from $\BE$ to $E$ via $\rho$. Then $G$ is the closed formal subscheme of $\wt G$ on which  $E$ is special and $\CO_{\!\Delta}$ acts through isogenies (Propositions~\ref{PropQIsogClosed} and \ref{PropSpecialClosed}).
\qed
\end{proof}

\begin{point} \label{PointActionOnGAndG'}
We define an action of the group $J(Q_\infty)$ on $G$.
By Proposition~\ref{PropQIsogMHat} there is an isomorphism $\epsilon_\infty$ from the group $J(Q_\infty)$ to the group of quasi-isogenies of $\wh{\ol\BM}$. We let $g\in J(Q_\infty)$ act on $G$ through $\epsilon_\infty$
\[
(\beta,\,\ul{\wh\CF},\wh\alpha) \es\mapsto\es (\beta,\,\ul{\wh\CF},\,\bar\beta^\ast\epsilon_\infty(g)\circ\wh\alpha) \,.
\]
This action commutes with the Galois action on $G$.
\end{point}

Let now $\Gamma\subset J(Q_\infty)$ be a discrete subgroup. We say that $\Gamma$ is \emph{separated} if it is separated in the profininte topology. This means that for every $g\in \Gamma$ there is a normal subgroup $\Gamma'\subset\Gamma$ of finite index that does not contain $g$.

Again following the arguments of Rapoport--Zink we prove in \cite{Crystals}

\begin{theorem} \label{ThmQuotientsOfG}
Let $\Gamma\subset J(Q_\infty)$ be a separated discrete subgroup. Then the quotient $\Gamma\backslash G$ is a locally noetherian, adic formal algebraic $\Spf\BF_{q^\ell}\dbl\zeta\dbr$-stack locally formally of finite type over $\Spf\BF_{q^\ell}$. Moreover, the 1-morphism $G\to\Gamma\backslash G$ is adic.
\end{theorem}

See the appendix \ref{DefFormalAlgebraicStacks} -- \ref{DefAdicMorphismOfStacks} for an explanation of this statement.

%
%

\bigskip

\section*{{\Large Part Three \es Uniformization}}  \label{PartIII}
\setcounter{equation}{0}

\mtaddtocont{\protect\contentsline {mtchap}{\protect\large\protect\rule{0mm}{2.5ex}Part Three \es Uniformization}{\thepage}}

We now turn towards the uniformization of the stacks $\AbSh^{r,d}_H$ at $\infty$. Therefore we again consider their base change
\[
\AbSh^{r,d}_H\times_C\;\Spf\BF_q\dbl\zeta\dbr\,.
\]
As remarked in Part Two these are \emph{formal algebraic $\Spf\BF_q\dbl\zeta\dbr$-stacks} (\ref{DefFormalAlgebraicStacks}). In analogy with the work of Rapoport--Zink~\cite{RZ} and Drinfeld~\cite{Drinfeld2} the space used to uniformize these stacks will be the formal scheme $G$. We view it as a formal algebraic $\Spf\BF_{q^\ell}\dbl\zeta\dbr$-stack (\ref{PointFormalSchemesAreStacks}). The uniformization we obtain will only be partial. To be precise, we find a closed subset $Z$ in $\AbSh^{r,d}_H\times_C\,\infty$ and we consider the formal completion $\AbSh^{r,d}_H{}_{\!\!/Z}$ of $\AbSh^{r,d}_H$ along $Z$. It is a formal algebraic $\Spf\BF_q\dbl\zeta\dbr$-stack (\ref{PropComplIsAlgStack}). We will uniformize $\AbSh^{r,d}_H{}_{\!\!/Z}$. This uniformization therefore takes place in the 2-category of formal algebraic $\Spf\BF_q\dbl\zeta\dbr$-stacks.

%
%

\section{Algebraizations}  \label{SectAlgebr}
\setcounter{equation}{0}

We first give still another interpretation of the moduli space $G$. Namely since the formal abelian sheaf $\wh{\ol{\BM}}$ comes from the abelian sheaf $\ol\BM$, the universal formal abelian sheaf on $G$ and its quasi-isogeny $\wh\alpha$ to $\wh{\ol{\BM}}$ can be algebraized. I.e.\ they too come from an abelian sheaf, namely from $\wh\alpha^\ast\ol\BM$. Now consider a scheme $S\in\Nilp_{\BF_q\dbl\zeta\dbr}$ and denote by $\bar S$ the closed subscheme of $S$ defined by $\zeta=0$.

\begin{definition}
  Let $G'$ be the contravariant functor $\Nilp_{\BF_q\dbl\zeta\dbr} \longto \CS ets$
\begin{eqnarray*}
S &\longmapsto \Bigl\{ &\text{Isomorphism classes of pairs}\es(\ul{\CF},\alpha) \es\text{where} \\
&&\es \bullet\es \ul{\CF}\text{ is an abelian sheaf of rank $r$ and dimension $d$ over }S\\
&&\es \bullet\es \alpha:\ul\CF_{\bar S}\to\ol\BM_{\bar S}\text{ is a quasi-isogeny which is an isomorphism over }C'\es\Bigr\}
\end{eqnarray*}
\end{definition}

Thereby two such pairs $(\ul{\CF},\alpha)$ and $(\ul{\CF}',\alpha')$ are isomorphic if there is an isomorphism between $\ul{\CF}$ and $\ul{\CF}'$ over $S$ which is compatible with $\alpha$ and $\alpha'$.

\begin{theorem} \label{ThmIsomOfGAndG'}
The functors $G$ and $G'\times_{\Spf \BF_q\dbl\zeta\dbr}\Spf \BF_{q^\ell}\dbl\zeta\dbr$ are canonically isomorphic as $\Gal(\BF_{q^\ell}/\BF_q)$-modules (where $\Gal(\BF_{q^\ell}/\BF_q)$ acts trivially on $G'$). 
\end{theorem}

\begin{proof}
We will exhibit two mutually inverse maps between these two functors. We start by describing the morphism $G'\times_{\Spf \BF_q\dbl\zeta\dbr}\Spf \BF_{q^\ell}\dbl\zeta\dbr\to G$. 

So let $(\ul{\CF},\alpha)\in G'(S)$ and $\beta:S\to\Spf \BF_{q^\ell}\dbl\zeta\dbr$.
By Construction~\ref{ConstrZDivGp} we obtain from $\ul{\CF}$ and $\beta$ a formal abelian sheaf $\ul{\wh\CF}$ of rank $r$ and dimension $d$ over $S$.
The quasi-isogeny $\alpha$ induces a quasi-isogeny of formal abelian sheaves $\wh\alpha: \es\ul{\wh\CF}_{\bar S}\to \bar\beta^\ast \wh{\ol\BM}$\,.
The triple $(\beta,\ul{\wh\CF},\wh\alpha)$ defines an $S$-valued point of the functor $G$. One easily checks that the morphism just constructed is $\Gal(\BF_{q^\ell}/\BF_q)$-equivariant.

Conversely let $(\beta,\ul{\wh\CF},\wh{\ol\alpha})\in G(S)$. By Proposition~\ref{PropQIsogLiftZDiv} there is a unique lift of $\wh{\ol\alpha}$ to a quasi-isogeny $\wh\alpha:\ul{\wh\CF}\to\beta^\ast\wh\BM$. From Proposition~\ref{PropPullbackQIsog} we obtain an abelian sheaf $\ul\CF:=\wh\alpha^\ast(\beta^\ast\BM)$ whose formal abelian sheaf is $\ul{\wh\CF}$, and a quasi-isogeny $\alpha:\ul\CF\to\beta^\ast\BM$ which is an isomorphism over $C'$ and which induces $\wh\alpha$ on the formal abelian sheaves. We have $(\ul\CF,\alpha_{\bar S})\in G'(S)$. Thus we have constructed a morphism $G\to G'$. Again one checks that this morphism is $\Gal(\BF_{q^\ell}/\BF_q)$-invariant.
The two morphisms just described are mutually inverse and yield the desired isomorphism between $G$ and $G'\times_{\Spf \BF_q\dbl\zeta\dbr}\Spf \BF_{q^\ell}\dbl\zeta\dbr$. 
\qed
\end{proof}

The action of $J(Q_\infty)$ on $G$ from \ref{PointActionOnGAndG'} induces an action of $J(Q_\infty)$ on $G'$, since it is compatible with the Galois action on $G$.

\begin{definition} \label{DefAlgebraization}
The pair $(\ul\CF,\alpha)\in G'(S)$ which is associated to an element $(\beta,\ul{\wh\CF},\wh\alpha)\in G(S)$ by Theorem~\ref{ThmIsomOfGAndG'} will be called the \emph{algebraization of $(\beta,\ul{\wh\CF},\wh\alpha)$}.
\end{definition}

\begin{bigexample} \label{ExTrivialAbSh}
We want to explain how a quasi-isogeny $\alpha:\ul\CF_{\bar S}\to\ol\BM_{\bar S}$ which is an isomorphism over $C'$ induces $H$-level structures on $\ul\CF$ for compact open subgroups $H\subset\GL_r(\BA_f)$.

Consider the abelian sheaf $\ol{\BM}=(\CM_i^{\oplus e},\Pi_i^{\oplus e},\tau_i^{\oplus e})$ from Example~\ref{ExAbSheaf} pulled back to ${\bar S}=\Spec \BF_{q^\ell}$. The restrictions $\CM_i^{\oplus e}|_{C'_{\bar S}}$ are all isomorphic via the morphisms $\Pi_i$. We denote this restriction by $\ol{\BM}|_{C'_{\bar S}}$. The same holds for the morphisms $\tau_i$. So we obtain a morphism
\[
\tau|_{C'_{\bar S}}:\sigma^\ast\ol{\BM}|_{C'_{\bar S}}\to\ol{\BM}|_{C'_{\bar S}}\,.
\]
Via the canonical identification $\ol{\BM}|_{C'_{\bar S}}=\CO_{C'_{\bar S}}^r$ this morphism is expressed by the block diagonal matrix
\[
\tau|_{C'_{\bar S}} \es = \es \left(
\begin{array}{ccc}
T \\[-3mm]
& \!\!\ddots\! & \\[-1mm]
& & T
\end{array}
\right)\!\cdot\sigma^\ast\, \qquad\text{where}\qquad 
T\es=\es\left(
\begin{array}{cccc}
0\, &  \ldots&&\! 1 \\[-1mm]
\raisebox{1mm}{$1$} & \ddots & & \raisebox{-1mm}{$\vdots$}\\[-3.4mm]
& \ddots & \hspace{-1.6mm}\ddots \\
& & 1 \!& 0
\end{array}\right)\;\in\;\GL_\ell(\BF_q)\,.
\]
Therefore multiplication with the matrix $V_e$ from (\ref{EqVandermonde}) defines an isomorphism 
\[
\psi:(\ol{\BM}|_{C'_{\bar S}},\tau|_{C'_{\bar S}})\es\isoto\es(\CO_{C'_{\bar S}}^r,\Id_r\cdot\sigma^\ast)
\]
that commutes with the $\sigma$-linear endomorphisms $\tau|_{C'_{\bar S}}$ on one and $\Id_r\cdot\sigma^\ast$ on the other side. 
This $\psi$ induces an isomorphism
\[
\gamma =(\psi|_{\BA_f})^\tau: \es (\ol{\BM}|_{\BA_f})^\tau \es\isoto\es \BA_f^r
\]
which gives rise to an $H$-level structure on $\ol{\BM}$.

Now let $S\in\Nilp_{\BF_q\dbl\zeta\dbr}$ be an arbitrary scheme equipped with a fixed morphism $\beta:S\to\Spf\BF_{q^\ell}\dbl\zeta\dbr$. Let $\ul\CF$ be an abelian sheaf over $S$ and let $\alpha:\ul\CF_{\bar S}\to\ol\BM_{\bar S}$ be a quasi-isogeny which is an isomorphism over $C'$. For an algebraically closed base point $\iota:s\to S$ the composition of $\alpha$ with $\beta^\ast\psi$ gives rise to an isomorphism 
\[
\gamma\circ(\alpha|_{\BA_f})^\tau\;:=\;\bigl((\beta^\ast\psi\circ\alpha)|_{\BA_f}\bigr)^\tau: (\iota^\ast\ul{\CF}|_{\BA_f})^\tau(s) \isoto \BA_f^r
\]
which is fixed under the action of $\pi_1(S,s)$ on the source. This induces an $H$-level structure on $\ul\CF$.
\end{bigexample}

%
%

\section{Closedness of the Uniformizable Locus}  \label{SectUnifClosed}
\setcounter{equation}{0}

In this section we show that an abelian sheaf of rank $r$ and dimension $d$ over an algebraically closed field in $\Nilp_{\BF_{q^\ell}\dbl\zeta\dbr}$ is isogenous to $\BM$ if and only if its formal abelian sheaf is isoclinic. The locus of these points inside $\AbSh^{r,d}_H\times_C\,\Spf\BF_{q^\ell}\dbl\zeta\dbr$ will be the one uniformized by $G$. We prove that this locus is formally closed. A couple of lemmas are needed beforehand.

\begin{lemma} \label{LemTEqualsPi}
Let $\ul\CF$ be an abelian sheaf of rank $r$ and dimension $d$ over a finite field in $\Nilp_{\BF_{q^\ell}\dbl\zeta\dbr}$. Assume that the formal abelian sheaf associated to $\ul\CF$ is isoclinic. Then for some integer $n>0$ divisible by $r$ there is an isomorphism $(\sigma^n)^\ast\ul\CF\cong\ul\CF$ which maps the isogeny
\[
T\;:=\;(\tau_i)\circ\sigma^\ast(\tau_i)\circ\ldots\circ(\sigma^{n-1})^\ast(\tau_i):\es(\sigma^n)^\ast\ul\CF\,[n] \to \ul\CF
\]
to the isogeny $(\Pi_i)^n:\ul\CF\,[n] \to \ul\CF$.
\end{lemma}

\begin{proof}
It suffices to prove the assertion for an abelian sheaf isogenous to $\ul\CF$.
By Proposition~\ref{PropIsoclinicCrystal} the formal abelian sheaf associated to $\ul\CF$ is isogenous to a formal abelian sheaf $(\wh\CF,F)$ which satisfies $\im F^r=z^d\wh\CF$. Pulling back $\ul\CF$ along this isogeny (Proposition~\ref{PropPullbackQIsog}) we can assume that $(\wh\CF,F)$ is the formal abelian sheaf associated to $\ul\CF$. Since $\ul\CF$ is defined over a finite field we certainly obtain $(\sigma^n)^\ast\ul\CF\cong\ul\CF$ for a suitable $n>0$. We may even assume that $r$ divides $n$. The isogeny $T$ is an isomorphism over $C'$. Now the equation $\im F^r=z^d\wh\CF$ implies that the image of $T$ is $\ul\CF(-\frac{nk}{\ell}\cdot\infty)$. Hence $T$ factors as $T=\alpha\circ(\Pi_i)^n$ for an automorphism $\alpha$ of $\ul\CF$. Since the automorphism group of $\ul\CF$ is finite by Proposition~\ref{PropAutomGpFinite} the lemma follows.
\qed
\end{proof}

\begin{definition}
Let $\ul\CF$ and $\ul\CF'$ be abelian sheaves of rank $r$ and dimension $d$ over a field $K$ in $\Nilp_{\BF_{q^\ell}\dbl\zeta\dbr}$. We define the $Q$-vector space
\[
\Hom_K^0(\ul\CF,\ul\CF') \es =\es \dirlim[D]\;\Hom_K\bigl(\ul\CF,\ul\CF'(D)\bigr)
\]
where the limit is taken over all effective divisors $D\subset C$. Furthermore for formal abelian sheaves $\ul{\wh\CF}$ and $\ul{\wh\CF}'$ over $K$ we define the $Q_\infty$-vector space
\[
\Hom_K^0(\ul{\wh\CF},\ul{\wh\CF}') \es =\es \Hom_K(\ul{\wh\CF},\ul{\wh\CF}')\otimes_{\BF_q\dbl z\dbr}Q_\infty\,.
\]
\end{definition}
Note that each time the invertible elements in $\Hom^0$ are precisely the quasi-isogenies.

\begin{corollary} \label{CorIsogExistFinFld}
Let $\ul\CF$ and $\ul\CF'$ be abelian sheaves of rank $r$ and dimension $d$ over a finite field $K$ in $\Nilp_{\BF_{q^\ell}\dbl\zeta\dbr}$. Assume that their associated formal abelian sheaves $\ul{\wh\CF}$ and $\ul{\wh\CF}'$ are isoclinic. Then for a suitable finite extension $K'/K$ we have an isomorphism
\[
\Hom^0_{K'}(\ul\CF,\ul\CF') \otimes_Q Q_\infty \isoto \Hom^0_{K^\alg}(\ul{\wh\CF},\ul{\wh\CF}')\,.
\]
\end{corollary}

\begin{proof}
We first assume that $C=\BP^1_{\BF_q}$. Let $\Spec K[z]\subset\BP^1_K$ be a neighborhood of $\infty=\Var(z)$. Restricted to this neighborhood all the sheaves $\CF_i$ and $\CF'_i$ are free of rank $r$. We fix bases for $i=0,\ldots,\ell-1$ and consider for the other $i$ the bases induced by the periodicity $\CF_i=\CF_{i-\ell}(k\cdot\infty)$. Then the morphism
\[
\bigoplus_{i=0}^{\ell-1}\tau_i:\es\bigoplus_{i=0}^{\ell-1}\sigma^\ast\CF_i\longto\bigoplus_{i=1}^{\ell}\CF_i
\]
is represented by a matrix $U\in M_{r\ell}\bigl(K[z]\bigr)$. We endow the formal abelian sheaf $(\wh\CF,F)$ of $\ul\CF$ with the induced basis. With respect to this basis $F$ is of the form $F\;=\;\wh U\cdot\sigma^\ast$ for the matrix
\[
\wh U\es= \es
\left(\begin{array}{cccc}
0\, &  \ldots&&\!\!\! \Id_r \\[-1mm]
\raisebox{1mm}{$\Id_r$} & \ddots & & \!\raisebox{-1mm}{$\vdots$}\\[-3.4mm]
& \ddots & \hspace{-1mm}\ddots \\
& & \!\!\Id_r \!& 0
\end{array}\right)\!\cdot U
\es=\es\sum_{\nu=0}^N\wh U_\nu z^\nu\es\in M_{r\ell}\bigl(K[z]\bigr)\,.
\]
The same holds for $\ul\CF'$ where we denote the corresponding matrix by $\wh U'$. Let us for a moment forget the structure of the (formal) abelian sheaves that is given by the $\Pi$'s. Let $n=rm$ be the integer from Lemma~\ref{LemTEqualsPi} and let $K'$ be the compositum of $\BF_{q^n}$ and $K$ inside $K^\alg$. We claim that there is an isomorphism of $Q_\infty$-vector spaces
\begin{eqnarray} \label{EqIsomOfHom0}
\qquad\bigl\{\,\Phi\in M_{r\ell}\bigl(Q\otimes_{\BF_q} K'\bigr): \;\Phi\,\wh U \,=\, \wh U'\,{}^{\sigma\!}\Phi\,\bigr\} \otimes_Q Q_\infty \es\isoto \qquad\qquad\qquad&&\\ 
\bigl\{\,\Phi\in M_{r\ell}\bigl(K^\alg\dpl z\dpr\bigr):\; \Phi\,\wh U \,=\, \wh U'\,{}^{\sigma\!}\Phi\,\bigr\}&&\nonumber
\end{eqnarray}
where the superscript ${}^{\sigma\!}\Phi$ denotes the application of $\sigma^\ast$ to the entries of the matrix $\Phi$. The injectivity is obvious. We have to prove the surjectivity. So let an element of the right hand side be given. After multiplying it with a power of $z$ it is represented by a matrix
\[
\Phi =\sum_{\mu=0}^\infty\Phi_\mu z^\mu\quad\in M_{r\ell}\bigl(K^\alg\dbl z\dbr\bigr)\,.
\]
We expand the equation $\Phi\,\wh U=\wh U'\,{}^{\sigma\!}\Phi$ into powers of $z$ and get for all $\mu$
\begin{equation} \label{EqForPhiMu}
\Phi_\mu\,\wh U_0 - \wh U'_0\,{}^{\sigma\!}\Phi_\mu\es = \es \sum_{\nu=1}^N\bigl(\Phi_{\mu-\nu}\,\wh U_\nu - \wh U'_\nu\,{}^{\sigma\!}\Phi_{\mu-\nu}\bigr)\,.
\end{equation}
Now by Lemma~\ref{LemTEqualsPi} we have 
\[
\wh U\,{}^{\sigma\!}\wh U\cdots\,{}^{\sigma^{rm-1}\!}\wh U\es = \es z^{dm}\cdot\Id_{r\ell} \qquad \text{and}\qquad
\wh U'\,{}^{\sigma\!}\wh U'\cdots\,{}^{\sigma^{rm-1}\!}\wh U'\es = \es z^{dm}\cdot\Id_{r\ell}\,.
\]
This implies the equation $z^{dm}\cdot\Phi=z^{dm}\cdot{}^{\sigma^{rm}\!}\Phi$, whence $\Phi_\mu = {}^{\sigma^{rm}\!}\Phi_\mu$ for all $\mu$. We find $\Phi_\mu\in M_{r\ell}(K')$. Now we consider for an integer $i$ the sequence of matrices $(\Phi_i,\ldots,\Phi_{i+N})$. As $i$ varies, these sequences run through a finite set. Therefore there are infinitely many $i$ giving rise to the same sequence. Let $j$ be the difference of two such $i$ and consider the matrix
\[
\Phi -z^j\Phi\es=:\es\sum_{\mu=0}^\infty \wt\Phi_\mu z^\mu\,.
\]
In this matrix we find a sequence with $\wt\Phi_i=\ldots=\wt\Phi_{i+N}=0$. Looking at equation (\ref{EqForPhiMu}) we see that we may set all $\wt\Phi_\mu=0$ for $\mu>i+N$ to obtain a matrix
\[
\wt\Phi\es=\es\sum_{\mu=0}^{i-1}\wt\Phi_\mu z^\mu\es\in \;M_{r\ell}\bigl(K'[z]\bigr)
\]
which satisfies $\wt\Phi\,\wh U=\wh U'\,{}^{\sigma\!}\wt\Phi$ and is congruent to $\Phi$ modulo $z^j$. 
As $j$ can be chosen arbitrarily large the surjectivity of (\ref{EqIsomOfHom0}) is established.

Now note that the $\Pi$'s on the (formal) abelian sheaves induce endomorphisms of the $Q_\infty$-vector spaces in (\ref{EqIsomOfHom0}). So the compatibility with the $\Pi$'s is a condition that cuts out isomorphic linear subspaces on both sides of (\ref{EqIsomOfHom0}). From this the corollary follows in the case $C=\BP^1$. For arbitrary $C$ consider a finite flat morphism $\pi:C\to\BP^1$ mapping $\infty_C$ to $\infty$. We have just proved the assertion for $\pi_\ast\ul\CF$ and $\pi_\ast\ul\CF'$. Since again the elements of $\CO_C$ induce endomorphisms of the $Q_\infty$-vector spaces in (\ref{EqIsomOfHom0}) we may deduce the assertion for $\ul\CF$ and $\ul\CF'$.
\qed
\end{proof}

Recall the abelian sheaf $\BM$ over $\Spf\BF_{q^\ell}\dbl\zeta\dbr$ and the level structures on $\BM$ from Example~\ref{ExTrivialAbSh}.

\begin{proposition} \label{PropAltDescrOfZ}
Let $S\in\Nilp_{\BF_{q^\ell}\dbl\zeta\dbr}$ be locally of finite type over $\BF_{q^\ell}\dbl\zeta\dbr$ and let $H\subset\GL_r(\BA_f)$ be a compact open subgroup. Let $\ul\CF$ be an abelian sheaf of rank $r$ and dimension $d$ with a rational $H$-level structure over $S$. Then for a point $s\in S$ the following assertions are equivalent:
\begin{enumerate}
\item 
The formal abelian sheaf associated to $\ul\CF_s$ is isoclinic.
\item
Over a finite extension of the residue field $\kappa(s)$ of $s$ there is a quasi-isogeny between $\ul\CF_s$ and $\BM_s$ which is compatible with the $H$-level structures on both sides.
\item
There is an abelian sheaf $\ul\CF'$ over a finite field $\BF\subset\BF_{q^\ell}$, a finite quasi-isogeny $\alpha:\ol\BM_\BF\to\ul\CF'$ over $\BF$ compatible with the $H$-level structures, and a quasi-isogeny $\phi_s:\ul\CF'_s\to\ul\CF_s$ which is an isomorphism over $C'$ and which is defined over a finite extension of $\kappa(s)$.
\end{enumerate}
\end{proposition}

\begin{proof}
Note that the formal abelian sheaf of $\ul\CF_s$ is isoclinic if and only if it is isogenous over an algebraically closed extension of $\kappa(s)$ to the formal abelian sheaf $\wh\BM$ of $\BM$. Therefore 2 implies 1. Clearly 2 follows from 3.

To prove that 1 implies 3 we proceed by induction on the transcendence degree of the residue field $\kappa(s)$ of $s$ over $\BF_q$. Observe that $\kappa(s)$ is finitely generated since $S$ is locally of finite type over $\BF_{q^\ell}\dbl\zeta\dbr$. If $\kappa(s)\subset\BF_q^{\,\alg}$ we obtain from Corollary~\ref{CorIsogExistFinFld} a quasi-isogeny $\alpha:\BM_\BF\to\ul\CF_\BF$ over a finite field $\BF\supset \kappa(s)$. Altering $\alpha$ by a quasi-isogeny of $\BM_\BF$ we can achieve that $\alpha$ is compatible with the $H$-level structures. Let $\wh\alpha$ be the induced quasi-isogeny on formal abelian sheaves. We set $\ul\CF'=\wh\alpha^\ast\ul\CF_\BF$. Then $\alpha$ factors as
\[
\ol\BM_\BF\es\xrightarrow{\es\alpha'}\es\ul\CF'\es\xrightarrow{\es\alpha''}\es\ul\CF_\BF\,,
\]
where $\alpha'$ is a finite quasi-isogeny compatible with the $H$-level structures, and $\alpha''$ is a quasi-isogeny which is an isomorphism over $C'$.

If $\kappa(s)\not\subset\BF_q^{\,\alg}$ we choose a point $s'$ of codimension 1 in the closure of $s$ inside $S$. 
By Theorem~\ref{ThmNewtonPolygon} the formal abelian sheaf of $\ul\CF_{s'}$ is also isoclinic. So the induction hypothesis asserts that there is an abelian sheaf $\ul\CF'$ over a finite field $\BF\subset\BF_{q^\ell}$, a finite quasi-isogeny $\alpha:\ol\BM_\BF\to\ul\CF'$ over $\BF$ compatible with the $H$-level structures, and a quasi-isogeny $\phi_K:\ul\CF'_K\to\ul\CF_K$ over an algebraic closure $K$ of $\kappa(s')$, which is an isomorphism over $C'$. A suitable extension of $\CO_{\ol{\{s\}},s'}$ is a power series ring in one variable $K\dbl\pi\dbr$ over $K$. 

Now consider the abelian sheaf $\ul\CF_{S'}$ over $S':=\Spec K\dbl\pi\dbr$ which is the pullback of $\ul\CF$ under the morphism $S'\to S$. Let $\ul{\wh\CF}_{S'}$ be the formal abelian sheaf associated to $\ul\CF_{S'}$. The Newton polygon of $\ul{\wh\CF}_{S'}$ is constant over $S'$. So by Theorem~\ref{ThmConstancy} there is a quasi-isogeny of formal abelian sheaves $\wh\phi:\ul{\wh\CF}'_{S'}\to\ul{\wh\CF}_{S'}$. 
After changing $\wh\phi$ by a quasi-isogeny of $\ul{\wh\CF}'$ we may assume that $\phi_K$ gives rise to $\wh\phi|_K$. Due to the Serre-Tate-Theorem~\ref{ThmSerreTate} $\wh\phi$ induces for all $n$ a lift of $\phi_K$ to a quasi-isogeny $\phi_n$ over $S'_n:=\Spec K\dbl\pi\dbr/(\pi^{n+1})$. By \ref{CorIsomOverC'} $\phi_n$ is an isomorphism over $C'$. Let $D$ be the divisor of the pole of $\wh\phi$ at $\infty$. Then $\phi_n:\ul\CF'_{S'_n}\to\ul\CF(D)_{S'_n}$ is a true isogeny. So in the limit we obtain a quasi-isogeny between $\ul\CF'_{S'}$ and $\ul\CF_{S'}$ over $\Spf K\dbl\pi\dbr$ which is an isomorphism over $C'$. By Grothendieck's existence theorem it comes from a quasi-isogeny over $\Spec K\dbl\pi\dbr$ which gives us a quasi-isogeny $\phi_s:\ul\CF'_s\to\ul\CF_s$ over $K\dpl\pi\dpr$. From the following lemma we obtain the desired quasi-isogeny over a finite extension of $\kappa(s)$.
\qed
\end{proof}

\begin{lemma}
Let $K\supset\BF_q$ be an arbitrary field. Consider two abelian sheaves $\ul\CF$ and $\ul\CF'$ over $K$ and a quasi-isogeny $\phi:\ul\CF'_L\to\ul\CF_L$ defined over some extension $L$ of $K$. Then there is a quasi-isogeny $\phi':\ul\CF'\to\ul\CF$ defined over a finite extension of $K$. If $\phi$ is an isomorphism over $C'$ we can find a $\phi'$ which is also an isomorphism over $C'$.
\end{lemma}

\begin{proof}
Let $D\subset C$ be an effective divisor such that $\phi:\ul\CF'_L\to\ul\CF_L(D)$ is an isogeny. In the description of the morphisms $\phi_i:\CF'_i\to\CF_i(D)$ there are only finitely many coefficients from $L$ involved due to the periodicity. Let $R$ be the $K$-subalgebra of $L$ generated by these coefficients. Now consider the locally closed subscheme $S$ of $\Spec R$ defined by the conditions that the $\phi_i$ give an isogeny. These are the equations $\Pi_i\circ\phi_i=\phi_{i+1}\circ\Pi'_i$, and $\tau_i\circ\sigma^\ast\phi_i=\phi_{i+1}\circ\tau'_i$, and the conditions that $\phi_i$ is injective, and that $\coker\phi_i$ is locally free of finite rank and supported on $\wt D\times S$ for an effective divisor $\wt D\subset C$. As $S$ is of finite type over $K$ we find a $K'$-rational point on it for a finite extension $K'/K$. The data over this point defines a quasi-isogeny $\phi':\ul\CF'_{K'}\to\ul\CF_{K'}$. If $\phi$ is an isomorphism over $C'$ it is clear that we can achieve the same for $\phi'$.
\qed
\end{proof}

\begin{corollary} \label{CorrZIsClosed}
Let $S\in\Nilp_{\BF_{q^\ell}\dbl\zeta\dbr}$ and let $H\subset\GL_r(\BA_f)$ be a compact open subgroup. Let $\ul\CF$ be an abelian sheaf of rank $r$ and dimension $d$ with a rational $H$-level structure over $S$. Then the set of points in $S$ over which there is a quasi-isogeny between $\ul\CF$ and $\BM$ which is compatible with the $H$-level structures, is closed.
\end{corollary}

\begin{proof}
By Theorem~\ref{ThmNewtonPolygon} the set of points over which the formal abelian sheaf associated to $\ul\CF$ is isoclinic is closed in $S$. If $S$ is locally of finite type over $\BF_{q^\ell}\dbl\zeta\dbr$ then the corollary follows from Proposition~\ref{PropAltDescrOfZ}.

Let $S$ be arbitrary. We only need to treat the case where $S$ is reduced. Then the abelian sheaf $\ul\CF$ induces a 1-morphism $f:S\to\AbSh^{r,d}_H\times_C \;\infty$ of algebraic $\Spec\BF_{q^\ell}$-stacks. As the question is local on $S$ we can assume that $S$ is quasi-compact. Since $\AbSh^{r,d}_H$ is locally of finite type over $C$ we may further assume that $f$ factors through a presentation $X\to\AbSh^{r,d}_H\times_C \;\infty$, where $X$ is a scheme of finite type over $\BF_{q^\ell}\dbl\zeta\dbr$. Then the closedness of the set on $X$ implies the closedness of the set on $S$.
\qed
\end{proof}

\begin{example} \label{ExPinksExample}
Consider the universal abelian sheaf over $M^{2,2}_I$ from the example in Section~\ref{SectExample}. Let $S=M^{2,2}_I\times_C\,\infty$. Pink has computed that the closed set from Corollary~\ref{CorrZIsClosed} is the proper subset
\[
\bigcup \,g\Var(\zeta,a_{11},a_{22},a_{21}) \es\subset\es S\,,
\]
where the union runs over all $g\in\GL_2(\BF_q)$ which act on the points $(a_{\mu\nu})\in S$ by conjugation $(a_{\mu\nu})\mapsto g(a_{\mu\nu})g^{-1}$\,; cf.\ \cite[{\S}7]{BH}. (Note that $\AbSh^{2,2}_{H_I}$ is a $\GL_2(\BF_q)$-torsor over $\AbSh^{2,2}$.)
\end{example}

\begin{remark}
There is an interesting consequence for the uniformizability of $t$-motives. Namely by work of Gardeyn~\cite{Gardeyn} the $t$-motive associated to an abelian sheaf over a complete extension $K$ of $\BF_q\dpl\zeta\dpr$ is uniformizable in the sense of Anderson~\cite{Anderson}, if and only if firstly the abelian sheaf extends to an abelian sheaf over the valuation ring $R$ of a finite extension of $K$, and secondly its reduction modulo the maximal ideal of $R$ is isogenous to $\BM$. 

Let $X$ be an admissible formal scheme in the sense of Raynaud~\cite{Raynaud,FRG} and let $(\ul{\CF},\bar\eta)$ be an abelian sheaf with $H$-level structure over $X$. Then we deduce that the set of points on the associated rigid-analytic space $X^\rig$ over which the abelian sheaf $\ul{\CF}$ is uniformizable, is formally closed. In particular if $X^\rig$ is quasi-compact then the complement of this set is also quasi-compact. For more elaboration on this issue see B{\"o}ckle--Hartl~\cite{BH}.
\end{remark}

\smallskip

In the remainder of this section we prove a weak result on the uniform existence of the quasi-isogeny between $\ul\CF$ and $\BM$ which above has been studied point-wise. It will suffice for our purposes in this article. Undoubtedly there should be much stronger results in this direction.

\begin{lemma}
Let $\BF\supset\BF_q$ be a finite field and let $S=\Spec R\in\Nilp_{\BF_q\dbl\zeta\dbr}$ be a reduced noetherian affine scheme. Consider abelian sheaves $\ul\CF$ and $\ul\CF'$ of rank $r$ and dimension $d$ over $S$ and $\BF$ respectively. Let $s:\Spec L\to S$ be an algebraically closed point over which a quasi-isogeny $\phi_s:\ul\CF'_s\to\ul\CF_s$ is given. Assume that $\phi_s$ is an isomorphism over $C'$. Then there exists a quasi-compact reduced scheme $S'$ of finite type over $S$ containing a lift of $s$, such that $\phi_s$ extends over all of $S'$ to a quasi-isogeny $\phi:\ul\CF'_{S'}\to\ul\CF_{S'}$ which is an isomorphism over $C'$.
\end{lemma}

\begin{proof}
We first treat the case where $C=\BP^1_{\BF_q}$ and $C'=\Spec\BF_q[t]$. Then the projective $R[t]$-module $\CF|_{C'_S}$ underlying $\ul\CF$ is a direct summand of a free $R[t]$-module
\[
\CF|_{C'_S}\oplus \CF^{nil} \es=\es\wt\CF
\]
of rank $\wt r$. We extend $\tau_{\ul\CF}$ by zero to the endomorphism $\wt\tau:=\tau_{\ul\CF}|_{C'_S}\oplus 0 : \sigma^\ast\wt\CF\to\wt\CF$. 
After choosing bases of $\wt\CF$ and $\ul\CF'|_{C'_S}$ the coherent sheaves on $C'_S$ with $\sigma$-linear endomorphism $(\wt\CF,\wt\tau)$ and $\ul\CF'|_{C'_S}$ are isomorphic to
\[
(R[t]^{\wt r}, U\cdot\sigma^\ast) \qquad\text{and}\qquad (R[t]^r, U'\cdot\sigma^\ast)
\]
for suitable matrices $U=\sum_\nu U_\nu t^\nu\in M_{\wt r}(R[t])$ and $U'=\sum_\nu U'_\nu t^\nu\in\GL_r(\BF[t])$.
The quasi-isogeny $\phi:\ul\CF'_{S'}\to\ul\CF_{S'}$ we are looking for, then corresponds to a matrix $\Phi=\sum_\mu \Phi_\mu \,t^\mu\in M_{\wt r\times r}(\CO_{S'}[t])$ satisfying the condition $\Phi\,U'= U\,{}^{\sigma\!}\Phi$. We expand this condition according to powers of $t$ to get
\begin{equation} \label{EqForMatrixPhi}
\Phi_\mu \,U'_0 \;- \;U_0\, {}^{\sigma\!}\Phi_\mu\es = \es \sum_{\nu=1}^\mu\bigl( U_\nu\,{}^{\sigma\!}\Phi_{\mu-\nu} \;- \;\Phi_{\mu-\nu}\,U'_\nu\bigr)\,.
\end{equation}
The quasi-isogeny $\phi_s$ over $\Spec L$ corresponds to a matrix for some integer $N$
\[
\ol\Phi\es=\es\sum_{\mu=0}^N \ol\Phi_\mu t^\mu\es\in\es M_{\wt r\times r}(L[t])
\]
satisfying (\ref{EqForMatrixPhi}).
In order to extend $\ol\Phi$ to a neighborhood of $s$ we simply adjoin the entries of indeterminant matrices $\Phi_0,\ldots,\Phi_N$ to the ring $R$ and divide out the relations (\ref{EqForMatrixPhi}) for $\mu=0,\ldots,N$. Thus we obtain an $R$-algebra $R''$ and an \'etale morphism $S''=\Spec R''\to S$. The point $s$ lifts to a point $s'':\Spec L\to S''$ by mapping $\Phi_\mu$ to $\ol\Phi_\mu$. 

We let $\wt S\subset S''$ be the closed subset on which the right hand side of (\ref{EqForMatrixPhi}) and the matrices $\Phi_\mu$ are zero for all $\mu>N$. Clearly $s''$ lies in $\wt S$. Over $\wt S$ the matrix $\Phi$ defines a morphism $f:\ul\CF'|_{C'_{\wt S}}\to(\wt\CF,\wt\tau)_{\wt S}$. Since $\tau_\ul{\CF'}|_{C'}$ is an isomorphism and $\wt\tau|_{\CF^{nil}}=0$ we see that $f$ factors through a morphism $\phi:\ul\CF'|_{C'_{\wt S}}\to\ul\CF|_{C'_{\wt S}}$. The coherent sheaves $\ker\phi$ and $\coker\phi$ are supported on closed subschemes of $C'_{\wt S}$. Observe that $\wt S$ is noetherian. So the projection of these closed subschemes to $\wt S$ are constructible subsets of $\wt S$ by Chevalley's theorem \cite[Corollaire IV.1.8.5]{EGA}. Their complement $S'$ is also constructible. Hence $S'$ is a finite union of locally close reduced subschemes of $\wt S$. We replace $S'$ with the \emph{disjoint} union of these subschemes. Then $S'$ is a quasi-compact reduced scheme of finite type over $S$. Now $\phi$ defines a quasi-isogeny $\ul\CF'_{S'}\to\ul\CF_{S'}$ which is an isomorphism over $C'_{S'}$. Moreover $s''\in S'$ and $\phi$ specializes to $\phi_s$ at $s''$.

If $C$ is arbitrary we consider a finite flat morphism $\pi:C\to\BP^1_{\BF_q}$ mapping $\infty_C$ to $\infty$. Then we have just proved the existence of a quasi-isogeny $\pi_\ast\ul\CF_{S'}\to\pi_\ast\ul\CF'_{S'}$ over $S'$. In order for this to give a quasi-isogeny $\ul\CF_{S'}\to\ul\CF'_{S'}$ it must commute with the elements of $\CO_C$. Clearly this condition cuts out a reduced closed subscheme of $S'$ which contains $s''$. We replace $S'$ by this subset. This concludes the proof of the lemma.
\qed
\end{proof}

\begin{proposition} \label{PropTauInvExist}
Let $S\in\Nilp_{\BF_{q^\ell}\dbl\zeta\dbr}$ be a quasi-compact and reduced scheme and let $H\subset\GL_r(\BA_f)$ be a compact open subgroup. Consider an abelian sheaf $(\ul\CF,\bar\gamma)$ of rank $r$ and dimension $d$ with rational $H$-level structure over $S$. Assume that for every point $s\in S$ the formal abelian sheaf associated to $\ul\CF_s$ is isoclinic. Then there exists a surjective morphism of $\BF_{q^\ell}$-schemes $S'\to S$ with $S'$ quasi-compact and reduced, and a quasi-isogeny over $S'$ between $\ul\CF_{S'}$ and $\BM_{S'}$ compatible with the $H$-level structures.
\end{proposition}

\begin{proof}
Since the data $(\ul\CF,\bar\gamma)$ involves only finitely many coefficients from $\CO_S$ we may assume that $S$ is of finite type over $\BF_{q^\ell}$. In particular $S$ is noetherian. (We could also use the fact that $\AbSh^{r,d}_H$ is locally of finite type over $C$.)

Let $s\in S$ be a point. From Proposition~\ref{PropAltDescrOfZ} we obtain an abelian sheaf $\ul\CF'$ over a finite field $\BF$, a finite quasi-isogeny $\alpha:\ol\BM_\BF\to\ul\CF'$ over $\BF$ compatible with the $H$-level structures, and a quasi-isogeny $\phi_s:\ul\CF'_s\to\ul\CF_s$ which is an isomorphism over $C'$. The quasi-isogeny $\phi_s$ is defined over an algebraically closed extension of $\kappa(s)$. So by the previous lemma there is a morphism $S'_s\to S$ of finite type from a quasi-compact reduced scheme $S'_s$, such that $s$ lifts to a point of $S'_s$, and $\phi_s$ extends to a quasi-isogeny $\phi$ over all of $S'_s$ which is an isomorphism over $C'$. Clearly the quasi-isogeny $\phi\circ\alpha:\ol\BM_{S'_s}\to\ul\CF_{S'_s}$ over $S'_s$ is compatible with the $H$-level structures. By Chevalley's theorem the image of the morphism $S'_s\to S$ is a constructible subset $S_s$ of $S$ containing $s$. Since $S$ is of finite type over $\BF_{q^\ell}$ the $S_s$ form a countable covering of $S$ by constructible subsets. By \cite[Corollaire 0.9.2.4]{EGA} finitely many of the $S_s$ suffice to cover $S$. We let $S'$ be the finite disjoint union of the corresponding $S'_s$. 
\qed
\end{proof}

%
%

\section{The Uniformization Theorem}  \label{SectUnifOfAbSh}
\setcounter{equation}{0}

Let $H\subset\GL_r(\BA_f)$ be a compact open subgroup. We will define 1-morphisms of formal algebraic $\Spf\BF_{q^\ell}\dbl\zeta\dbr$-stacks from $G$ to the stacks $\AbSh^{r,d}_H\times_C\;\Spf\BF_{q^\ell}\dbl\zeta\dbr$ of abelian sheaves with rational $H$-level structure (Definition~\ref{DefStack3}). For this purpose recall the $H$-level structures on $\ol\BM$ constructed in Example~\ref{ExTrivialAbSh}.

\begin{point}
Let $S$ be in $\Nilp_{\BF_q\dbl\zeta\dbr}$ and denote by $\bar S$ its special fiber. The morphism $G\to G'$ from Theorem~\ref{ThmIsomOfGAndG'} associates to an element $(\beta,\ul{\wh\CF},\wh\alpha)\in G(S)$ an abelian sheaf $\ul{\CF}$ and a quasi-isogeny $\alpha:\ul\CF_{\bar S}\to\bar\beta^\ast\ol\BM$ which is an isomorphism over $C'$. For an algebraically closed base point $\iota:s\to S$ the quasi-isogeny $\alpha$ has lead in \ref{ExTrivialAbSh} to an isomorphism 
\[
\gamma\circ(\alpha|_{\BA_f})^\tau: (\iota^\ast\ul{\CF}|_{\BA_f})^\tau(s) \isoto \BA_f^r
\]
which is fixed under the action of $\pi_1(S,s)$ on the source.
Then we define a 1-morphism $\Theta$ of formal algebraic $\Spf\BF_{q^\ell}\dbl\zeta\dbr$-stacks by the following map on $S$-valued points
\begin{eqnarray} \label{EqMorphismTheta}
\Theta:\es G\times\GL_r(\BA_f)/H &\es\longto\es&\AbSh^{r,d}_H\;\times_C\;\Spf\BF_{q^\ell}\dbl\zeta\dbr\\[2mm]
(\beta,\ul{\wh\CF},\wh\alpha) \times h H & \mapsto & \bigl(\ul{\CF},\,h^{-1}\gamma\,(\alpha|_{\BA_f})^\tau\bigr) \times \beta\,. \nonumber
\end{eqnarray}
It is equivariant with respect to the right $\GL_r(\BA_f)$-action on the projective systems on both sides of (\ref{EqMorphismTheta}).
\end{point}


\begin{point} \label{PointActionOnSource}
There is an action of $J(Q)$ on the source which we describe next. Recall that we have defined in \ref{PointActionOnGAndG'} an action of $J(Q_\infty)$ on $G$ through the isomorphism $\epsilon_\infty$ from $J(Q_\infty)$ to the group of quasi-isogenies of $\wh{\ol\BM}$.
We let $J(Q)$ act on $G$ via the inclusion $J(Q)\subset J(Q_\infty)$.

On the other hand we have a morphism
\[
\epsilon^\infty: J(Q) \into \GL_r(\BA_f)
\]
which is defined by the commutative diagram
\[
\begin{CD}
(\ol{\BM}|_{\BA_f})^\tau & @>(g|_{\BA_f})^\tau>> & (\ol{\BM}|_{\BA_f})^\tau \\
@V\gamma VV & & @VV\gamma V \\
\BA_f^r & @>\epsilon^\infty(g)>> & \BA_f^r
\end{CD}
\]
for $g\in J(Q)$. Note that $\epsilon^\infty$ identifies $J(Q)$ with the diagonal embedding of $\GL_r(Q)$ into $\GL_r(\BA_f)$.

We define a left action of the group $J(Q)$ on $G\times\GL_r(\BA_f)/H$ 
\[
(\beta,\ul{\wh\CF},\wh\alpha) \times h H \es\longmapsto\es \bigl(\beta,\,\ul{\wh\CF},\,\bar\beta^\ast\epsilon_\infty(g)\circ\wh\alpha\bigr) \;\times \;\epsilon^\infty(g)h H \,.
\]
\end{point}

\begin{proposition} \label{PropAction}
We abbreviate $Y=G\times\GL_r(\BA_f)/H$. The action of $J(Q)$ induces a 1-isomorphism of formal algebraic $\Spf\BF_{q^\ell}\dbl\zeta\dbr$-stacks
\[
Y\times J(Q) \es\isoto\es Y\times_{\AbSh^{r,d}_H\times\Spf\BF_{q^\ell}\dbl\zeta\dbr}\;Y\,.
\]
\end{proposition}

\begin{proof}
We have to show that the functor between the categories of $S$-valued points on both sides is an equivalence. Note that in the stack $Y$ the only morphisms are the identities. From this full faithfulness follows. For essential surjectivity 
we have to show that two $S$-valued points of $Y$ lie in the same orbit if and only if they are mapped to the same point by $\Theta$. So let $(\beta,\ul{\wh\CF},\wh\alpha,h H)$ and $(\beta,\ul{\wh\CF}',\wh\alpha',h'H)$ be in $Y(S)$ and let $(\ul{\CF},\alpha,h H)$ and $(\ul{\CF}',\alpha,h'H)$ be their algebraizations. 

First we assume that 
\[
\bigl(\,\beta,\,\ul{\wh\CF}',\wh\alpha', h'H\,\bigr)\es =\es\bigl(\,\beta,\,\ul{\wh\CF},\,\bar\beta^\ast\epsilon_\infty(g)\circ\wh\alpha,\, \epsilon^\infty(g)h H \,\bigr)
\]
in $Y(S)$ for a $g\in J(Q)$.
Consider the following diagram of quasi-isogenies over $\bar S$
\begin{equation} \label{DiagPropAction}
\begin{CD}
\ul{\CF}_{\bar S} & @>\alpha>> & \ol{\BM}_{\bar S} \\
@V\bar\phi VV & & @VVg_{\bar S} V \\
\ul{\CF}'_{\bar S} & @>\alpha'>> & \ol{\BM}_{\bar S} \,,
\end{CD}  
\end{equation}
where $g_{\bar S}$ is the quasi-isogeny obtained from $g$ by Proposition~\ref{PropQIsogMS}. The quasi-isogeny $\bar\phi$ defined by this diagram is finite by construction. By Proposition~\ref{PropQIsogLiftAbSh} it lifts uniquely to a quasi-isogeny $\phi:\ul\CF\to\ul\CF'$. We claim that $\phi$ is finite. Namely the induced quasi-isogeny $\wh\phi$ on formal abelian sheaves satisfies $\wh\phi_{\bar S}=\id_{\ul{\wh\CF}_{\bar S}}$. Hence we find $\wh\phi=\id_{\ul{\wh\CF}}$ by the uniqueness of the lift. This shows that $\phi$ is a finite quasi-isogeny.

Now $\bar\phi$ induces on $\ul{\CF}$ the $H$-level structure 
\[
h^{-1}\,\epsilon^\infty(g)^{-1}\,\gamma\;(\alpha'|_{\BA_f})^\tau\,(\bar\phi|_{\BA_f})^\tau\es=\es h^{-1}\,\gamma\;(\alpha|_{\BA_f})^\tau\,.
\]
Therefore the two points have the same image under $\Theta$.

Conversely let $\phi:\ul{\CF}\to\ul{\CF'}$ be a finite quasi-isogeny inducing an equality of $H$-level structures on $\ul{\CF}$
\[
h'{}^{-1}\,\gamma\;(\alpha'|_{\BA_f})^\tau \,(\phi|_{\BA_f})^\tau \es =\es h^{-1}\,\gamma\;(\alpha|_{\BA_f})^\tau
\]
This time diagram (\ref{DiagPropAction}) defines a quasi-isogeny $g_{\bar S}$ from $\ol{\BM}_{\bar S}$ to itself. By Proposition~\ref{PropQIsogMS} it comes from an element $g\in J(Q)$. Then we have
\[
\bigl(\,\beta,\,\ul{\wh\CF}',\wh\alpha', h'H\,\bigr)\es =\es\bigl(\,\beta,\,\ul{\wh\CF},\,\bar\beta^\ast\epsilon_\infty(g)\circ\wh\alpha,\, \epsilon^\infty(g)h H \,\bigr)\,.
\]
\qed
\end{proof}

Hence the map $\Theta$ factors through a 1-morphism of formal algebraic $\Spf\BF_{q^\ell}\dbl\zeta\dbr$-stacks
\[
\raisebox{-0.2mm}{$J(Q)$}\backslash \raisebox{0mm}{$G\times\GL_r(\BA_f)$}/\raisebox{-0.2mm}{$H$} \es\longto\es\AbSh^{r,d}_H\;\times_C\;\Spf\BF_{q^\ell}\dbl\zeta\dbr\,.
\]
Note that the quotient $\raisebox{-0.2mm}{$J(Q)$}\backslash \raisebox{0mm}{$G\times\GL_r(\BA_f)$}/\raisebox{-0.2mm}{$H$}$ is a formal algebraic $\Spf\BF_{q^\ell}\dbl\zeta\dbr$-stack due to Theorem~\ref{ThmQuotientsOfG}. Indeed, the subgroup $J(Q) \into J(Q_\infty)\times\GL_r(\BA_f)$ is discrete. Hence
\[
\raisebox{-0.2mm}{$J(Q)$}\backslash \raisebox{0mm}{$G\times\GL_r(\BA_f)$}/\raisebox{-0.2mm}{$H$} \es\cong\es\coprod_\Gamma \Gamma\backslash G
\]
where $\Gamma$ runs through a countable set of subgroups of $J(Q_\infty)$ of the form 
\[
\bigl(J(Q_\infty)\times g H g^{-1}\bigr)\cap J(Q) \es\subset\es J(Q_\infty)\,.
\]
These are separated discrete subgroups. 


\begin{proposition} \label{PropDefTheta'}
The 1-morphism $\Theta$ defines a 1-morphism of formal algebraic $\Spf\BF_q\dbl\zeta\dbr$-stacks
\[
\Theta':\es G'\times\GL_r(\BA_f)/H \es\longto\es\AbSh^{r,d}_H\;\times_C\;\Spf\BF_q\dbl\zeta\dbr
\]
which is invariant with respect to the action of $J(Q)$ on the source defined in Section~\ref{SectAlgebr}.
\end{proposition}

\begin{proof}
We have to show that $\Theta$ and the action of $J(Q)$ commute with the Galois-descent data on the source and the target of the 1-morphism $\Theta$. For the $J(Q)$-action this was already observed in \ref{PointActionOnGAndG'}. For $\Theta$ we must check that
\[
\bigl(\ul{\CF},\,h^{-1}\gamma\;(\alpha|_{\BA_f})^\tau\bigr) \es = \es\bigl(\ul{\CF},\,h^{-1}{}^{\sigma\!}\gamma\;(\alpha|_{\BA_f})^\tau\bigr)
\]
in $\AbSh^{r,d}_H$, where ${}^{\sigma\!}\gamma$ is obtained from $\gamma$ by the action of the generator $Frob_q$ of $\Gal(\BF_{q^\ell}/\BF_q)$. Now ${}^{\sigma\!}\gamma^{-1}\circ\gamma$ comes in fact from an automorphism of $\ol{\BM}$ (cf.\ \ref{ExTrivialAbSh}). This induces an automorphism of $\ul{\CF}$ which carries the $H$-level structure $h^{-1}\gamma\;(\alpha|_{\BA_f})^\tau$ to $h^{-1}{}^{\sigma\!}\gamma\;(\alpha|_{\BA_f})^\tau$.
\qed
\end{proof}

Let $Z$ be the set of points $s:\Spec L\to\AbSh^{r,d}_H\;\times_C\;\infty$ such that the universal abelian sheaf over $s$ is isogenous to $\ol{\BM}$ over an algebraic closure of $L$. Consider the preimage $Z'\subset\AbSh^{r,d}_H\times_{\BF_q}\BF_{q^\ell}$ of $Z$ under the base change morphism coming from $\BF_q\subset\BF_{q^\ell}$.

\begin{lemma}
The set $Z'$ can also be described as the set of points $s$ over which the associated formal abelian sheaf is isoclinic. In particular $Z$ and $Z'$ are closed subsets.
\end{lemma}

\begin{proof}
The formal abelian sheaf $\wh{\ol\BM}$ of $\ol\BM$ is isoclinic. Therefore $Z'$ is contained in the set of the lemma. 
Conversely let $s$ belong to this later set. Since $\AbSh^{r,d}_H$ is locally of finite type over $C$ we can assume that $s$ comes from a point on a local presentation $X\to\AbSh^{r,d}_H\times_C\;\infty$ where $X$ is a scheme of finite type over $\BF_q$. From Proposition~\ref{PropAltDescrOfZ} we obtain a quasi-isogeny $\alpha:\ul\CF_s\to\BM_s$ over an algebraic closure of $L$. Hence $s$ belongs to $Z'$. 

Now Theorem~\ref{ThmNewtonPolygon} implies that the subset $Z'$ is closed. Namely $Z'$ is the complement of the open substack on which the associated Newton polygon lies strictly below the Newton polygon of $\wh{\ol\BM}$. Therefore also the image $Z$ of $Z'$ is closed.
\qed
\end{proof}

We denote by $\AbSh^{r,d}_H{}_{\!\!/Z}$ the formal completion of $\AbSh^{r,d}_H$ along $Z$ (\ref{DefFormalComplOfStacks}). It is a formal algebraic $\Spf\BF_q\dbl\zeta\dbr$-stack. By its definition the 1-morphism $\Theta'$ factors through $\AbSh^{r,d}_H{}_{\!\!/Z}$. Indeed, if a point $s\in\AbSh^{r,d}_H\;\times_C\;\infty$ lies in the image of $\Theta'$ the formal abelian sheaf associated to $s$ is isogenous to $\wh\BM$ by definition and hence isoclinic. We can now formulate our

\begin{uniftheorem} \label{ThmUnifOfAbSh}
There are $\GL_r(\BA_f)$-equivariant 1-isomorphisms of formal algebraic $\Spf\BF_q\dbl\zeta\dbr$-stacks
\begin{eqnarray*}
\bar\Theta:\es\raisebox{-0.2mm}{$J(Q)$}\backslash \raisebox{0mm}{$G\times\GL_r(\BA_f)$}/\raisebox{-0.2mm}{$H$} &\es\isoto\es&\AbSh^{r,d}_H{}_{\!\!/Z}\times_{\BF_q}\BF_{q^\ell}\,,\\[2mm]
\bar\Theta':\es\raisebox{-0.2mm}{$J(Q)$}\backslash \raisebox{0mm}{$G'\times\GL_r(\BA_f)$}/\raisebox{-0.2mm}{$H$} &\es\isoto\es&\AbSh^{r,d}_H{}_{\!\!/Z}\,.
\end{eqnarray*}
\end{uniftheorem}

\medskip

\begin{example}
In the case of elliptic sheaves 
\[
\AbSh^{r,1}_H{}_{\!\!/Z}\times_{\BF_q}\;\,\BF_{q^\ell}\es=\es\AbSh^{r,1}_H\times_C \;\Spf\BF_{q^\ell}\dbl\zeta\dbr\,.
\]
This follows from the fact that there is only one polygon between the points $(0,0)$ and $(r,1)$ with non-negative slopes and integral break points, namely the straight line. So all formal abelian sheaves are isoclinic and $Z$ is all of $\AbSh^{r,d}_H\times_C\;\infty$.
In this case the open and closed subscheme of $G$ on which the universal quasi-isogeny has height zero is the formal scheme $\Omega^{(r)}$ used by Drinfeld~\cite{Drinfeld2}. A detailed account on this can be found in Genestier~\cite{Genestier}. Therefore we have an isomorphism of formal schemes $G\cong\BZ\times\Omega^{(r)}$. This decomposition of $G$ is compatible with the decomposition of $\AbSh^{r,1}_{H_I} \cong\BZ\times\DrMod^r_I$ from Example~\ref{ExDrinfeldModules}. So Drinfeld's uniformization theorem which announces an isomorphism of formal schemes
\[
\raisebox{-0.2mm}{$\GL_r(Q)$}\backslash \raisebox{0mm}{$\Omega^{(r)}\times\GL_r(\BA_f)$}/\raisebox{-0.2mm}{$H_I$} \es\isoto\es \DrMod^r_I\times_C \;\Spf\BF_{q^\ell}\dbl\zeta\dbr
\]
is equivalent to ours.
\end{example}

\begin{example}
Consider the algebraic stack $\AbSh^{2,2}_H$ from Section~\ref{SectExample}. We describe $Z(0):=Z\cap\AbSh^{2,2}_H(0)$. In \ref{ExPinksExample} we have remarked that
\[
Z\,\cap\, M^{2,2}_I \times_C\infty \es=\es\bigcup \,g\Var(\zeta,a_{11},a_{22},a_{21})\,.
\]
Now we claim that
\[
Z(0)\es=\es\bigl(Z\cap M^{2,2}_I \times_C\infty\bigr) \;\cup\; \bigl(\AbSh^{2,2}_H(0)\,\setminus\,M^{2,2}_I\bigr)\times_C\infty\,.
\]
Indeed let $s:\Spec L\to\bigl(\AbSh^{2,2}_H(0)\,\setminus\,M^{2,2}_I\bigr)\times_C\infty$ be a point. We must show that the abelian sheaf $\ul\CF$ over $L$ is isogenous to $\ol\BM$. By what was said in Section~\ref{SectExample} we have $\CF_0=\CO_{\BP^1_L}(m\cdot\infty)\oplus\CO_{\BP^1_L}(-m\cdot\infty)$ for an integer $m\ge 1$. With respect to this basis, $\tau_0$ is described by a matrix
\[
\tau_0\es=\es\left(
\begin{array}{c@{\quad\es}c}
a_0+a_1 t & b_0+\ldots+b_{2m+1} \,t^{2m+1}\\[2mm]
0 & d_0+d_1 t\qquad\quad
\end{array}\right)\cdot\sigma^\ast\;.
\]
Due to the presence of the $I$-level structure we may assume $a_0=d_0=1$ and $b_0=0$. Since $\coker\tau_0$ is supported at $\infty$ we must have $a_1=d_1=0$. Now let $u_i\in L^\alg$ be solutions of the equations $u_i^q-u_i+b_i=0$ for $i=1,\ldots,2m+1$. Then the isomorphism
\[
\ol\BM|_{C'}\;\to\; \ul\CF|_{C'}\,,\quad \left(\begin{array}{c}x\\y\end{array}\right)\es\mapsto 
\left(
\begin{array}{c@{\quad}c}
1 & u_0 t+\ldots+u_{2m+1} \,t^{2m+1}\\[2mm]
0 & 1\qquad\qquad
\end{array}
\right)\cdot\left(\begin{array}{c}x\\y\end{array}\right)
\]
extends to a quasi-isogeny $\ol\BM\to\ul\CF$ over $L^\alg$. Hence $s$ belongs to $Z(0)$.
\end{example}

\begin{remark}
In \cite{RZ} Rapoport--Zink study the uniformization of Shimura varieties of EL- and PEL-type. One of their theorems yields the uniformization of the formal completion of a Shimura variety along the most supersingular isogeny class \cite[6.30]{RZ}. The Newton polygon in this isogeny class is maximal. In this sense our Uniformization Theorem is closely analogous to theirs. 
Beyond this,
Rapoport--Zink also obtain uniformization theorems of other isogeny classes. There is no doubt that these theorems too have counterparts for abelian sheaves.
\end{remark}

The remainder of this article is devoted to the proof of the Uniformization Theorem.

%
%

\section{Proof of the Uniformization Theorem}  \label{SectProofUnifThm}
\setcounter{equation}{0}

\begin{point}
By Proposition~\ref{PropDefTheta'} it suffices to prove the assertion for the 1-morphism $\bar\Theta$. We fix the following notation. On the formal scheme $G$ we let $\CJ$ be the largest ideal of definition of $G$; cf.\ \cite[I$_{\rm new}$, 10.5.4]{EGA}. For an integer $n\geq0$ we denote by $G_n$ the scheme $(G,\CO_G/\CJ^{n+1})$. We set
\begin{eqnarray*}
Y &\es:=\es& G\times\GL_r(\BA_f)/H \,,\\[2mm]
Y_n &\es:=\es& G_n\times\GL_r(\BA_f)/H \,,\\[2mm]
\CY & \es:=\es& \raisebox{-0.2mm}{$J(Q)$}\backslash \raisebox{0mm}{$G\times\GL_r(\BA_f)$}/\raisebox{-0.2mm}{$H$}\,,\\[2mm]
\CX & \es:=\es & \AbSh^{r,d}_H{}_{\!\!/Z}\times_{\Spf\BF_q\dbl\zeta\dbr}\;\,\Spf\BF_{q^\ell}\dbl\zeta\dbr\,.
\end{eqnarray*}
Let $S\in\Nilp_{\BF_q\dbl\zeta\dbr}$ and let $(\ul{\CF},{\bar\gamma},\beta)\in\CX(S)$ which we consider as a 1-morphism $S\to\CX$. We have to show that the stack
\[
\CY\times_\CX S
\]
is a scheme mapping isomorphically to $S$. The proof relies on several intermediate lemmas.
\end{point}

\begin{lemma}
The 1-morphism $\bar\Theta:\CY\to\CX$ is a 1-monomorphism of formal algebraic stacks, i.e.\ for every $S\in \Nilp_{\BF_{q^\ell}\dbl\zeta\dbr}$ the functor $\bar\Theta(S):\CY(S)\to\CX(S)$ is fully faithful.
\end{lemma}

\begin{proof}
We can view the formal algebraic $\Spf\BF_{q^\ell}\dbl\zeta\dbr$-stacks $\CX$ and $\CY$ as $\Spec\BF_{q^\ell}\dbl\zeta\dbr$-stacks $\wt{\CX}$ and $\wt{\CY}$ by setting for an $\Spec\BF_{q^\ell}\dbl\zeta\dbr$-scheme $S$
\[
\wt{\CX}(S)\es=\es\left\{
\begin{array}{cl}
\CX(S) \es& \text{if } S\in\Nilp_{\BF_{q^\ell}\dbl\zeta\dbr} \\[1mm]
\emptyset & \text{if } S\notin\Nilp_{\BF_{q^\ell}\dbl\zeta\dbr} \,.
\end{array}\right.
\]
Then the assertion follows from Proposition~\ref{PropAction}, Lemma~\ref{LemIsomFiberProducts} and \cite[Proposition 3.8]{LaumonMB}.
\qed
\end{proof}

\begin{lemma} \label{LemThetaRepresentable}
The 1-morphism of algebraic $\Spec\BF_{q^\ell}$-stacks $\bar\Theta_\red:\CY_\red\to\CX_\red$ (\ref{PointReducedStacks}) is representable by a morphism of schemes.
\end{lemma}

\begin{proof}
This follows from the fact that every 1-monomorphism of algebraic stacks is representable by a morphism of schemes; cf.\ \cite[Th{\'e}or{\`e}me A.2 and Corollaire 8.1.3]{LaumonMB}.
\qed
\end{proof}

\begin{lemma} \label{LemThetaRedSurj}
$\bar\Theta_\red:\CY_\red\to\CX_\red$ is surjective.
\end{lemma}

\begin{proof}
Note that $\CX_\red$ is the closed substack $Z\subset\AbSh^{r,d}_H\;\times_C\;\infty$ with its induced reduced structure. Let $s\in\CX_\red$ be a point. By definition of $Z$ there is a quasi-isogeny $\alpha:\ul\CF_s\to\BM_s$. We can multiply it with a quasi-isogeny of $\BM$ and thus assume that $\alpha$ is compatible with the $H$-level structures. The induced quasi-isogeny $\ul\CF_s\to\wh\alpha^\ast\BM_s$ is finite and also compatible with the $H$-level structures. Therefore $s$ lies in the image of $\bar\Theta_\red$.
\qed
\end{proof}

\begin{lemma} \label{LemTauInvExist}
Let $S\in\Nilp_{\BF_{q^\ell}\dbl\zeta\dbr}$ be quasi-compact and reduced, and consider a 1-morphism $S\to\CX_\red$. Then there exists a surjective morphism of $\BF_{q^\ell}$-schemes $S'\to S$ with $S'$ quasi-compact and reduced, and a 2-commutative diagram
\[
\begin{CD}
S' & @>>> & Y_\red \\
@VVV & & @VVV \\
S & @>>> & \CX_\red\,.
\end{CD}
\]
\end{lemma}

\begin{proof}
This is just a reformulation of Proposition~\ref{PropTauInvExist}.
\qed
\end{proof}

\begin{lemma} \label{LemThetaRedQC}
$\bar\Theta_\red:\CY_\red\to\CX_\red$ is quasi-compact.
\end{lemma}

\begin{proof}
Let $S\in\Nilp_{\BF_{q^\ell}\dbl\zeta\dbr}$ be quasi-compact and reduced, and let $S\to\CX_\red$ be a 1-morphism. Let $S'\to S$ be the surjective morphism from Lemma~\ref{LemTauInvExist}. It gives rise to a surjective morphism of schemes over $\BF_{q^\ell}$
\[
\CY_\red\times_{\CX_\red}S \es\lsurj\es \CY_\red\times_{\CX_\red}S'\es \lsurj \es Y_\red\times_{\CX_\red} S'\,.
\]
Since $Y_\red\times_{\CX_\red} S' \cong S'\times J(Q)$ by Proposition~\ref{PropAction}, we obtain an epimorphism $S'\rsurj\CY_\red\times_{\CX_\red}S$. Now the lemma follows from the quasi-compactness of $S'$.
\qed
\end{proof}

\begin{lemma}
$\bar\Theta_\red:\CY_\red\to\CX_\red$ is proper.
\end{lemma}

\begin{proof}
Since $\CY_\red$ and $\CX_\red$ are locally of finite type and $\bar\Theta_\red$ is quasi-compact, we see that $\bar\Theta_\red$ is of finite type. Being a 1-monomorphism it is also separated. So it remains to prove that it is universally closed. For this we use the valuative criterion \cite[Th{\'e}or{\`e}me 7.3]{LaumonMB}. Let $R$ be a valuation ring with $\Spec R\in\Nilp_{\BF_{q^\ell}\dbl\zeta\dbr}$ and let $K$ be its field of fractions. We have to show that for every 2-commutative diagram
\[
\begin{CD}
\Spec K & @>>> & \CY_\red \\
@VVV & & @VVV \\
\Spec R & @>>> & \CX_\red
\end{CD}
\]
there exists a finite extension $R'/R$ of valuation rings and a 1-morphism of $\Spec\BF_{q^\ell}$-stacks $\Spec R'\to\CY_\red$  which 2-commutes with the above diagram. However, the existence of this data follows from Lemma~\ref{LemTauInvExist}. 
\qed
\end{proof}

\begin{lemma} \label{LemThetaIsAdic}
The 1-morphism $\bar\Theta:\CY\to\CX$ is adic (\ref{DefAdicMorphismOfStacks}).
\end{lemma}

\begin{proof}
We will show that $\bar\Theta_\red$ is a 1-isomorphism. Since both $\CX$ and $\CY$ are adic (\ref{PointIdealOfDefinition}) this suffices.

Let $P:X\to\CX_\red$ be a presentation. Then $P$ is an epimorphism since $\CX_\red$ is an algebraic $\Spec\BF_{q^\ell}$-stack. Thus it suffices to show that $\CY_\red\times_{\CX_\red}X\to X$ is a 1-isomorphism. From the previous lemmas we know that it is a proper monomorphism of schemes, hence a closed immersion. Since it is also surjective and $X$ is reduced it is an isomorphism as desired. 
\qed
\end{proof}

\begin{lemma} \label{LemThetaEtale}
The 1-morphism $\bar\Theta:\CY\to\CX$ is {\'e}tale.
\end{lemma}

\begin{proof}
Since quasi-isogenies of $z$-divisible groups lift to infinitesimal neighborhoods we first see that $Y\to\CX$ is {\'e}tale.

Now let $\CJ$ be an ideal of definition of $\CX$ and let $\CZ$ be the closed substack defined by $\CJ$. Since $\bar\Theta$ and also the presentation $Y\to\CY$ are adic (Theorem~\ref{ThmQuotientsOfG}), we obtain 1-morphisms
\[
Y\times_\CX\CZ \es\to\es \CY\times_\CX \CZ \es\to\es \CZ\,.
\]
of algebraic $\Spec\BF_{q^\ell}\dbl\zeta\dbr$-stacks. Since these 1-morphisms are representable by morphisms of schemes, we can apply \cite{EGA} to see that $\CY\times_\CX \CZ \es\to\es \CZ$ is {\'e}tale. This proves the lemma.
\qed
\end{proof}

\smallskip

We can now finish the

\begin{proof}[of the Uniformization Theorem~\ref{ThmUnifOfAbSh}]
Keep the notation of the proof of Lemma~\ref{LemThetaEtale}. We have to show that
\[
\CY\times_\CX \CZ \es\to\es \CZ
\]
is a 1-isomorphism of algebraic $\Spec\BF_q\dbl\zeta\dbr$-stacks.
Let $S$ be a $\Spec\BF_q\dbl\zeta\dbr$-scheme and let $S\to\CZ$ be a 1-morphism. Then from the previous lemmas we conclude that
\[
\CY\times_\CX S \es\to\es S
\]
is an {\'e}tale monomorphism of schemes, hence an open immersion. Being also surjective, it is indeed an isomorphism. This completes the proof of the Uniformization Theorem.
\qed
\end{proof}

%
%

\begin{appendix}

\section{Background on Formal Algebraic Stacks} \label{AppBackgroundStacks}
\setcounter{equation}{0}

For a general introduction to the theory of algebraic stacks we refer to Laumon--Moret-Bailly~\cite{LaumonMB} or Deligne--Mumford~\cite{DM}. In this appendix we propose the notion of formal algebraic stacks which generalizes the notion of algebraic stacks in the same way as formal schemes are a generalization of usual schemes. In fact much of the theory of algebraic stacks can be developed also for formal algebraic stacks. See Hartl~\cite{FormalStacks} for details.

For a scheme $S$ let $\Sch_S$ be the category of $S$-schemes equipped with the {\'e}tale topology. An \emph{algebraic stack} $\CX$ over $S$ is defined as a category fibered in groupoids over $\Sch_S$ satisfying further conditions \cite[Definition 4.1]{LaumonMB}. We transfer this concept to the formal category.

For the rest of this appendix we let $S$ be a formal scheme. We denote by $\Nilp_S$ the category of schemes over $S$ on which an ideal of definition of $S$ is locally nilpotent. We remind the reader that every scheme may be considered as a formal scheme having $(0)$ as an ideal of definition. In this sense every $U\in\Nilp_S$ is itself a formal scheme. We equip $\Nilp_S$ with the {\'e}tale topology. 
We make the following definitions (compare \cite{Laumon,Knutson}).

\begin{definition}
A \emph{formal $S$-space} is a sheaf of sets on the site $\Nilp_S$.
\end{definition}

\begin{definition} \label{DefFormalAlgebraicSpaces}
A \emph{(quasi-separated) formal algebraic $S$-space} is a formal $S$-space $X$ such that 
\begin{enumerate}
\item 
the diagonal morphism $X\to X\times_S X$ is relatively representable by a quasi-compact morphism of formal schemes and
\item
there is a formal scheme $X'$ over $S$ and a morphism of formal $S$-spaces $X'\to X$ which is representable (automatic because of 1) by an {\'e}tale surjective morphism of formal schemes.
\end{enumerate}
\end{definition}

\begin{definition} \label{DefFormalStacks}
A \emph{formal $S$-stack $\CX$} is a category $\CX$ fibered in groupoids over $\Nilp_S$ such that 
\begin{enumerate}
\item 
for every $U\in\Nilp_S$ and every $x,y\in \CX(U)$ the presheaf
\begin{eqnarray*}
\CI som(x,y):\es \Nilp_U &\es\to\es& \CS ets \\
(V\to U) &\mapsto& \Hom_{\CX(V)}(x_V,y_V) 
\end{eqnarray*}
is in fact a sheaf on $\Nilp_U$,
\item
for every covering $U_i\to U$ in $\Nilp_S$ all descent data for this covering are effective.
\end{enumerate}
\end{definition}

\begin{definition} \label{DefRepresentableStacks}
A formal $S$-stack is called \emph{representable} if it is 1-isomorphic to a formal algebraic space.

A 1-morphism $\CX\to\CY$ of formal $S$-stacks is called \emph{representable} if for every $U\in\Nilp_S$ and every $y\in\CY(U)$ viewed as a 1-morphism $U\to\CY$ of formal $S$-stacks the fiber product $\CX\times_\CY U$ (in the 2-category of formal $S$-stacks) is representable.
\end{definition}

\begin{definition} \label{DefFormalAlgebraicStacks}
A \emph{(quasi-separated) formal algebraic $S$-stack} is a formal $S$-stack $\CX$ such that
\begin{enumerate}
\item 
the diagonal 1-morphism of formal $S$-stacks
\[
\CX \es\to\es\CX\times_S\CX
\]
is representable, separated, and quasi-compact,
\item
there exists a formal algebraic $S$-space $X$ and a 1-morphism of formal $S$-stacks
\[
P:\es X\to\CX
\]
which is representable (automatic because of 1) by a smooth and surjective morphism of formal algebraic $S$-spaces.
\end{enumerate}
The 1-morphism $P$ is called a \emph{presentation of $\CX$}. We say that $\CX$ is of \emph{DM-type} if the presentation $P$ in 2 can be chosen {\'e}tale.
\end{definition}

\begin{point} \label{PointReducedStacks}
Let $\CX$ be a formal algebraic $S$-stack and let $P:X\to\CX$ be a presentation. We define the underlying reduced stack $\CX_\red$ as follows: For every $U\in\Nilp_S$ we let $\CX_\red(U)$ be the full subcategory of $\CX(U)$ whose objects are the $x\in \CX(U)$ such that there is a covering $U'\to U$ in $\Nilp_S$, an element $x'\in X_\red(U')$, and an isomorphism in $\CX(U')$ between $x_{U'}$ and $P(x')$. Then $\CX_\red$ is an algebraic $S_\red$-stack. If moreover $\CX$ is of DM-type then $\CX_\red$ is an algebraic $S_\red$-stack in the sense of Deligne--Mumford. In this way we obtain from every 1-morphism $f:\CY\to\CX$ of formal algebraic $S$-stacks a 1-morphism $f_\red:\CY_\red\to\CX_\red$ of algebraic $S_\red$-stacks.
We say that $f$ is \emph{locally formally of finite type} if $f_\red$ is locally of finite type.
\end{point}

\begin{point} \label{PointIdealOfDefinition}
For an algebraic stack $\CX$ one can define its \emph{structure sheaf $\CO_\CX$} which is a sheaf on the lisse-{\'e}tale site of $\CX$; cf. \cite[{\S} 12]{LaumonMB}. Then one has the usual bijection between closed substacks of $\CX$ and quasi-coherent sheaves of ideals of $\CO_\CX$.

The same can be done for formal algebraic $S$-stacks $\CX$. In this setting we say that a sheaf of ideals $\CJ$ of $\CO_\CX$ is an \emph{ideal of definition of $\CX$} if for some (any) presentation $P:X\to\CX$ of $\CX$ the ideal sheaf $P^\ast\CJ$ is an ideal of definition of $X$. The formal algebraic $S$-stack is called \emph{$\CJ$-adic} if $\CJ^n$ is an ideal of definition for every $n$. It is called \emph{adic} if it is $\CJ$-adic for some $\CJ$. If $\CX$ (i.e. $X$) is locally noetherian then there exists a unique \emph{largest ideal of definition $\CK$ of $\CX$}, namely the one defining the closed substack $\CX_\red$ of $\CX$. Note that if $\CX$ is $\CJ$-adic for some $\CJ$ then it is also $\CK$-adic.
One easily verifies the following proposition which generalizes \ref{PointReducedStacks}.
\end{point}

\begin{proposition} \label{PropIdealOfDefinition}
Let $\CI$ be an ideal of definition of $S$ and let $\CJ$ be an ideal of definition of a formal algebraic $S$-stack $\CX$ with $\CI\cdot\CO_\CX\subset\CJ$. Then the closed substack of $\CX$ which is defined by the ideal $\CJ$ is an algebraic $(S,\CO_S/\CI)$-stack.
\end{proposition}

\begin{definition} \label{DefAdicMorphismOfStacks}
A 1-morphism $f:\CY\to\CX$ of locally noetherian formal algebraic $S$-stacks is called \emph{adic} if for some (any) ideal of definition $\CJ$ of $\CX$ the ideal $f^\ast\CJ$ is an ideal of definition of $\CY$. 
\end{definition}

\smallskip

We discuss some examples.

\begin{bigexample} \label{PointFormalSchemesAreStacks}
Every formal scheme $G$ over $S$ can be viewed as a sheaf of sets on the category $\Nilp_S$. Moreover, every sheaf on $\Nilp_S$ is a formal $S$-stack. Therefore we can view every quasi-separated formal scheme $G$ over $S$ as a formal algebraic $S$-stack of DM-type.
\end{bigexample}

\begin{bigexample} \label{PointBaseChangeStacks}
Let $S^\alg$ be a scheme and let $S_0$ be a closed subscheme. (We allow the case $S_0=S^\alg$.) We let the formal scheme $S$ be the formal completion of $S^\alg$ along $S_0$. 
If $\CX$ is an (algebraic) $S^\alg$-stack (in the sense of Deligne--Mumford) then $\CX\times_{S^\alg}S$ is an adic formal (algebraic) $S$-stack (of DM-type).
\end{bigexample}

We generalize this example as follows.

\begin{definition} \label{DefFormalComplOfStacks}
Let $S^\alg$ be a scheme and let $S_0$ be a closed subscheme. We let the formal scheme $S$ be the formal completion of $S^\alg$ along $S_0$. Let $\CX$ be an algebraic $S^\alg$-stack and let $\CZ\subset\CX$ be a closed substack, contained in $\CX\times_{S^\alg} S_0$. We view the objects of $\CX$ as 1-morphisms $U\to\CX$ for varying $U\in\Sch_{S^\alg}$. We define the \emph{formal completion $\wh{\CX}_\CZ$ of $\CX$ along $\CZ$} as the full subcategory of $\CX$ consisting of those objects $U\to\CX$ such that $U_\red\to\CX$ factors through $\CZ$, i.e.\ such that there exists a 2-commutative diagram of 1-morphisms
\[
\begin{CD}
U & @>>> & \CX \\
@AAA & & @AAA \\
U_\red & @>>> & \CZ\,.
\end{CD}
\]
\end{definition}

Note that if $\CX$ is an $S^\alg$-scheme this definition coincides with the usual definition of formal completion along a closed subscheme. 

One verifies directly that $\wh{\CX}_\CZ$ is an $S^\alg$-stack. The embedding $\wh{\CX}_\CZ\to\CX$ is a 1-morphism of $S^\alg$-stacks which is a 1-monomorphism, i.e.\ the functors $\wh{\CX}_\CZ(U)\to\CX(U)$ are fully faithful. This implies that the diagonal 1-morphism
\[
\wh{\CX}_\CZ\es\to\es\wh{\CX}_\CZ\times_\CX \wh{\CX}_\CZ
\]
is a 1-isomorphism of $S^\alg$-stacks \cite[2.3]{LaumonMB}. As a consequence we obtain

\begin{lemma} \label{LemIsomFiberProducts}
Let $\CY\to\wh{\CX}_\CZ$ and $\CY'\to\wh{\CX}_\CZ$ be two 1-morphisms of $S^\alg$-stacks. Then there is a 1-isomorphism of $S^\alg$-stacks
\[
\CY\times_{\wh{\CX}_\CZ}\CY'\es\isoto\es\CY\times_\CX \CY'\,.
\]
\end{lemma}

The fibration $\wh{\CX}_\CZ\to\Sch_{S^\alg}$ factors through $\Nilp_S\into\Sch_{S^\alg}$. Note the fact that for an {\'e}tale covering $U_i\to U$ in $\Sch_{S^\alg}$ we have $U\in\Nilp_S$ if and only if $U_i\in\Nilp_S$ for all $i$. This implies that $\wh{\CX}_\CZ$ is a formal $S$-stack. We show that $\wh{\CX}_\CZ$ is even a formal algebraic $S$-stack. Namely, condition 1 of Definition~\ref{DefFormalAlgebraicStacks} can be read off from the following 2-cartesian diagram of $S^\alg$-stacks
\[
\begin{CD}
\wh{\CX}_\CZ & @>\sim>> & \wh{\CX}_\CZ\times_\CX \wh{\CX}_\CZ & @>>> & \wh{\CX}_\CZ\times_{S^\alg} \wh{\CX}_\CZ & @<\sim<< & \wh{\CX}_\CZ\times_S \wh{\CX}_\CZ\\
& & & & @VVV & & @VVV \\
& & & & \CX & @>>> & \CX \times_{S^\alg}\CX\,.
\end{CD}
\]
Condition 2 results from the fact that every presentation $X\to\CX$ of $\CX$ induces a presentation $\wh{X}_Z\to\wh{\CX}_\CZ$ of $\wh{\CX}_\CZ$ by the formal completion $\wh{X}_Z$ of $X$ along the closed subscheme $Z=X\times_\CX \CZ$. Thus we have proved

\begin{proposition} \label{PropComplIsAlgStack}
Keep the notation of Definition~\ref{DefFormalComplOfStacks}. Then the formal completion $\wh{\CX}_\CZ$ of $\CX$ along $\CZ$ is a formal algebraic $S$-stack. If $\CX$ is an algebraic $S^\alg$-stack in the sense of Deligne--Mumford then $\wh{\CX}_\CZ$ is of DM-type. If $\CJ$ is the ideal sheaf on $\CX$ defining the closed substack $\CZ$ then $\wh{\CX}_\CZ$ is $\CJ\cdot\CO_{\wh{\CX}_\CZ}$-adic.
\end{proposition}

\begin{point} \label{PointClosedSubset}
Let $S^\alg$ be a scheme and let $\CX$ be an algebraic $S^\alg$-stack.
There is a notion of \emph{points} of $\CX$. Namely a point of $\CX$ is given by a 1-morphisms $\Spec K\to\CX$ for an $S^\alg$-field $K$. The set $|\CX|$ of points of $\CX$ forms a topological space; cf.\ \cite[{\S} 5]{LaumonMB}. Let $Z\subset|\CX|$ be a closed subset, i.e.\ there is an open substack $\CU\subset\CX$ such that $Z=|\CX|\setminus|\CU|$. We can equip $Z$ in a unique way with a structure of reduced closed substack \cite[4.10]{LaumonMB}. By \ref{DefFormalComplOfStacks} we can consider the formal completion of $\CX$ along $Z$.
\end{point}

\smallskip

\begin{bigexample}[Quotients] \label{PointQuotientsAreStacks}
Let $U\in\Nilp_S$ and let $G$ be a formal $U$-group space (i.e.\ a group object in the category of formal $U$-spaces). A (left) \emph{$G$-torsor} is a formal $U$-space $P$ with an action of $G$ (from the left) such that there is a covering $U'\to U$ in $\Nilp_S$ for which $P\times_U U'$ is $G\times_U U'$-isomorph to $G\times_U U'$ which acts on itself by left translation.

Let $X$ be a formal $S$-space, $Y$ an $X$-space (i.e.\ a formal $S$-space equipped with a morphism $Y\to X$) and $G$ an $X$-group space which acts on $Y$ from the left. We define the quotient stack $G\backslash Y$ as the following category fibered in groupoids over $\Nilp_S$: For every $U\in\Nilp_S$ the category $(G\backslash Y)(U)$ consists of all triples $(x,P,\alpha)$ where $x\in X(U)$, $P$ is a $G\times_{X,x}U$-torsor and $\alpha:P\to Y_{X,x}U$ is a $G\times_{X,x}U$-equivariant morphism of formal $U$-spaces.
One easily verifies that the quotient $G\backslash Y$ is a formal $S$-stack.

In particular, if $X$ is an adic formal algebraic $S$-space and $G$ a finite {\'e}tale $S$-group scheme then the quotient $G\backslash X$ is even an adic formal algebraic $S$-stack of DM-type. In this case, the canonical projection $X\to G\backslash X$ is an {\'e}tale presentation of $G\backslash X$.
\end{bigexample}


\end{appendix}


\end{document}